\documentclass[12pt]{amsart}

\usepackage{amssymb}
\usepackage{latexsym} % to use some math symbols
\def\nc{\newcommand}

\def\lam{\lambda}\def\ep{\epsilon}

\def\om{\omega}

\def\Ph{\Phi} \def\Om{\Omega}
\def\ra{\rightarrow}
\def\Ra{\Rightarrow}

\def\x{\times}
\def\R{\mathbf R}
\def\DD{\mathcal D}
\nc\pa{\partial}

\nc\CC{\mathbf{C}}
\nc\RR{\mathbf{R}}
\nc\QQ{\mathbf{Q}}
\nc\ZZ{\mathbf{Z}}
\nc\NN{\mathbf{N}}

\nc\m[1]{\left| #1\right|}
\nc\norm[1]{\left\| #1\right\|}
\newtheorem{theorem}{Theorem}[section] 
\newtheorem{lemma}[theorem]{Lemma}
\newtheorem{corollary}[theorem]{Corollary}
\newtheorem{proposition}[theorem]{Proposition}
\newtheorem{definition}[theorem]{Definition}% Use {\rm ...}
\newtheorem{remark}[theorem]{Remark}        % Use {\rm ...}
\numberwithin{equation}{section}

\begin{document}

\title[Quasilinear and Hessian equations]
{Quasilinear and Hessian equations\\ of Lane--Emden type}

\author[Nguyen Cong Phuc]
{Nguyen Cong Phuc}
\address{Department of Mathematics,
University of Missouri,
Columbia, MO 65211, USA}
\email{nguyencp@math.missouri.edu}

\author[Igor E. Verbitsky]
{Igor E. Verbitsky$^*$}
\address{Department of Mathematics,
University of Missouri,
Columbia, MO 65211, USA}
\email{igor@math.missouri.edu}

\thanks{$^*$Supported in part 
by NSF Grant DMS-0070623.}

\begin{abstract} 

The existence problem is solved, and global pointwise estimates 
of solutions are obtained for  quasilinear and Hessian equations of Lane--Emden type, 
including the following 
two  model problems: 
$$ -\Delta_p u  = u^q + \mu, \qquad F_k[-u] = u^q + \mu, \qquad u \ge 0,  
$$
on  $\RR^n$, or on a bounded domain $\Omega \subset \RR^n$. 
Here $\Delta_p$ is the $p$-Laplacian defined by $\Delta_p u = {\rm div} \, ( \nabla u |\nabla u|^{p-2})$, 
and $F_k[u]$ is the $k$-Hessian defined as the sum of $k\times k$ principal minors of the 
Hessian matrix $D^2 u$ ($k=1,2, \ldots, n$); 
$\mu$ is a nonnegative  measurable function (or measure) 
on $\Omega$.  

The solvability of these classes of equations  in the renormalized (entropy) or viscosity sense   
has been an open problem even for good data $\mu \in L^s (\Omega)$, $s>1$. Such results are deduced 
from our existence criteria with the sharp exponents $s = \frac {n(q-p+1)} {pq}$ for the first equation, and 
$s = \frac{n(q-k)}{2kq}$ for the second one. Furthermore, 
a complete characterization of removable singularities is given. 

Our methods are based on systematic use of Wolff's potentials, dyadic models,  and nonlinear trace inequalities. 
We make use of recent advances in potential theory and PDE due to Kilpel\"ainen and Mal\'y, Trudinger 
and Wang, and Labutin. This enables us to treat  singular solutions, 
nonlocal operators,  and distributed singularities, and 
develop the theory simultaneously for quasilinear equations and equations of Monge-Amp\`ere type.

\end{abstract}

\maketitle

\section{Introduction}\label{Introduction}

We study a class of quasilinear and fully nonlinear equations and inequalities  
with nonlinear source terms, 
which appear in such diverse areas as quasi-regular mappings,   
non-Newtonian fluids,  reaction-diffusion problems, and stochastic control. 
In particular, the following two model equations are of substantial interest: 
 \begin{equation}\label{I0}  
  -\Delta_p u  = f(x, u), \qquad   F_k[-u] = f(x, u) ,  
 \end{equation} 
on   $\RR^n$, or on a bounded domain $\Om \subset \RR^n$, where $f(x, u)$ is a nonnegative function,  
convex and nondecreasing in $u$ for $u\ge 0$.   
Here $\Delta_p u = {\rm div} \, (\nabla u \, |\nabla u|^{p-2})$ is the $p$-Laplacian ($p>1$), 
and $F_k[u]$ is the $k$-Hessian ($k=1, 2, \ldots, n$) defined by 
 \begin{equation}\label{I01} 
 F_k[u] =\sum_{1\leq i_{1}<\cdots<i_{k}\leq n}\lambda_{i_{1}}\cdots\lambda_{i_{k}},
  \end{equation} 
where $ \lambda_{1},...,\lambda_{n}$ are the eigenvalues of the Hessian matrix 
$D^{2}u$. In other words, $F_{k}[u]$ is the sum of the $k\times k$ principal minors of $D^{2}u$, which coincides   
with the Laplacian $F_1[u] = \Delta u$ if $k=1$, and the Monge--Amp\`ere operator 
$F_n[u] = {\rm det} \, (D^2 u)$ if $k=n$. 

The form in which we write the second equation in  (\ref{I0})   
is chosen only for the sake of convenience, in order to emphasize the profound analogy 
between the quasilinear and Hessian equations. Obviously,  it may be stated  as  
$(-1)^k \, F_k[u] =  f(x, u)$, $u\geq 0$, or $F_k[u] = f(x, -u)$, $u\leq 0$.   

The existence and regularity theory, local and global estimates of 
sub- and super-solutions, the Wiener criterion, 
and Harnack  inequalities associated with  the $p$-Laplacian, as well as more general quasilinear 
operators,  can be found in  \cite{HKM}, \cite{IM}, \cite{KM2}, \cite{M1}, \cite{MZ}, 
 \cite{S1}, \cite{S2},  \cite{SZ}, \cite{TW4} where many 
fundamental results, and relations to other areas of analysis and geometry are presented. 

The theory  of fully 
nonlinear equations of Monge-Amp\`ere type which involve the $k$-Hessian operator  
$F_k[u]$ was originally developed  by Caffarelli, Nirenberg and 
Spruck,  Ivochkina, and Krylov in the classical setting. We refer to \cite{CNS}, \cite{GT}, \cite{Gu}, 
\cite{Iv}, \cite{Kr}, \cite{Ur}, \cite{Tru2}, \cite{TW1} for these and further results. Recent developments concerning the notion 
of the $k$-Hessian measure, weak convergence, and pointwise potential estimates 
due to Trudinger and Wang \cite{TW2}--\cite{TW4}, and Labutin \cite{L}  are used extensively in this paper.  

We are specifically interested in quasilinear and fully nonlinear equations of  
Lane--Emden type:
\begin{equation}\label{I02}  
  -\Delta_p u  = u^q, \quad {\rm and} \qquad  F_k[-u] = u^q,  \qquad u \ge 0 \quad {\rm in} \, \, \Om,
 \end{equation}
  where $p>1$, $q>0$, $k=1, 2, \ldots, n$, and the corresponding nonlinear inequalities:  
  \begin{equation}\label{I03}  
  -\Delta_p u  \ge  u^q, \quad {\rm and} \qquad F_k[-u] \ge u^q,  \qquad u \ge 0 \quad {\rm in} \, \, \Om. 
 \end{equation}
 The latter can be stated 
 in the form of the inhomogeneous equations with measure data, 
  \begin{equation}\label{I04}  
  -\Delta_p u  =  u^q +\mu, \qquad   F_k[-u] = u^q + \mu,  \qquad u \ge 0 \quad {\rm in} \, \, \Om,
 \end{equation}
 where $\mu$ is a nonnegative Borel measure on $\Om$.

The difficulties arising in studies of such equations and inequalities with competing nonlinearities are well known. 
In particular,  (\ref{I02})   may have singular solutions \cite{SZ}.  The existence problem for (\ref{I04}) has 
been open (\cite{BV2}, Problems 1 and 2; see also \cite{BV1}, \cite{BV3}, \cite{Gre}) even for the quasilinear equation 
$-\Delta_p u = u^q +f$ with good data $f \in L^s(\Om)$, $s>1$. 
Here solutions are generally understood in the renormalized (entropy) sense
 for quasilinear equations,  and viscosity, or  $k$-convexity   
 sense, for fully nonlinear equations of Hessian type (see \cite{BMMP}, 
 \cite{DMOP}, \cite{JLM}, \cite{TW1}--\cite{TW3}, \cite{Ur}). Precise definitions of these 
classes of {\it admissible solutions} are given  in Sec.~3, Sec.~6, and Sec.~7 below.

In this paper, we present a unified approach to   (\ref{I02})--(\ref{I04})  which  
 makes it possible to attack a number of open problems. 
It is  based on  global pointwise estimates, nonlinear integral inequalities in Sobolev 
spaces of fractional order,  and analysis of dyadic 
models, along with the weak convergence  and  Hessian measure results \cite{TW2}--\cite{TW4}. The latter   
are used to bridge the gap between the dyadic models and 
 partial differential equations. Some of these techniques were developed in  
 the linear case, in the framework 
of Schr\"odinger operators and harmonic analysis 
  \cite{ChWW}, \cite{Fef}, 
\cite{KS}, \cite{NTV},   \cite{V1}, \cite{V2}, 
and applications to semilinear equations   \cite{KV},  \cite{VW}, \cite{V3}.

Our goal is to establish  {\it necessary and sufficient\/}  conditions for 
the existence of  solutions to (\ref{I04}), sharp pointwise and integral estimates for solutions to (\ref{I03}), and 
 a complete characterization of  
removable singularities for (\ref{I02}). 
We are mostly concerned with admissible solutions to the corresponding equations and inequalities. 
However, even for locally bounded solutions, as in \cite{SZ},  our results 
yield new pointwise and integral estimates, and Liouville-type theorems.

 In the ``linear case'' $p=2$ and $k=1$,  problems (\ref{I02})--(\ref{I04})  
 with nonlinear sources 
 are associated with the names of 
 Lane and Emden, as well as Fowler.   Authoritative 
  historical and bibliographical comments can be found in \cite{SZ}. An up-to-date   survey  
  of the vast literature on nonlinear elliptic equations with measure data is given in 
  \cite{Ver}, including a thorough discussion of related work due to D. Adams and Pierre \cite{AP}, 
 Baras and Pierre \cite{BP},   Berestycki,
  Capuzzo-Dolcetta,  and  Nirenberg \cite{BCDN},  Brezis and Cabr\'e \cite{BC}, Kalton and Verbitsky \cite{KV}.

  It is worth mentioning that related  equations with absorption, 
    \begin{equation}\label{I04d}  
  - \Delta u  +  u^q = \mu, \qquad u \ge 0 \quad  {\rm in} \, \, \Om,
 \end{equation}
 were studied in detail  by  B\'enilan and Brezis, Baras and Pierre, 
 and Marcus and V\'eron analytically for $1<q<\infty$,    and by Le Gall, and Dynkin and Kuznetsov using
 probabilistic methods when $1<q\le 2$ (see \cite{D}, \cite{Ver}). For a general 
class of semilinear equations 
  \begin{equation}\label{I04db}  
  - \Delta u  +  g(u) = \mu, \qquad u \ge 0 \quad  {\rm in} \, \, \Om,
 \end{equation}
where $g$ belongs to the class of continuous nondecreasing  functions such that $g(0)=0$, 
sharp existence results have been obtained quite recently 
 by Brezis, Marcus, and 
Ponce \cite{BMP}. It is well known that equations with absorption 
generally require ``softer'' methods of analysis, and the conditions on $\mu$ which ensure the existence  of solutions 
are less stringent than in the case of equations with source terms. 
 
 Quasilinear problems of Lane--Emden type (\ref{I02})--(\ref{I04})  
   have been  studied extensively over the past 15  years. 
   Universal estimates for solutions, Liouville-type theorems,  and  analysis of removable singularities are  due to 
  Bidaut-V\'eron,  Mitidieri and Pohozaev  \cite{BV1}--\cite{BV3}, 
  \cite{BVP}, \cite{MP}, and Serrin and Zou \cite{SZ}.  (See also 
\cite{BiD},  \cite{Gre}, \cite{Ver}, and the literature cited there.) The profound difficulties in this theory 
are highlighted by the presence of the two  critical exponents, 
 \begin{equation}\label{I4} 
 q_* = \tfrac{n(p-1)}{n-p}, \qquad q^*= \tfrac{n (p-1) + p} {n-p},    
 \end{equation}
 where $1<p< n$. 
 As was shown in   \cite{BVP}, \cite{MP}, and \cite{SZ},  the quasilinear inequality (\ref{I04})
   does not 
 have nontrivial  weak solutions  
 on $\RR^n$, or exterior domains, 
 if $q\le  q_*$. For  $q > q_*$ , there exist 
 $u \in W^{1, \, p}_{\rm loc}\cap L^\infty_{\rm loc} $   which obey  (\ref{I03}), as well 
as singular solutions to 
(\ref{I02})  on $\RR^n$. However, for the existence of nontrivial solutions  
$u \in W^{1, p}_{\rm loc}\cap L^\infty_{\rm loc} $ to 
 (\ref{I02}) on $\RR^n$, it is necessary and sufficient that $q\geq q^*$ \cite{SZ}. In the 
``linear case'' $p=2$, this is classical  \cite{GS},  \cite{BP}, \cite{BCDN}. 

The following local 
estimates of solutions to quasilinear inequalities are used extensively in the studies mentioned above (see, e.g., \cite{SZ}, Lemma 2.4). Let  $B_R$ denote a ball of radius $R$ such that 
 $B_{2R} \subset \Omega$. Then, for every solution $u \in W^{1, p}_{\rm loc}\cap L^\infty_{\rm loc} $ 
 to the inequality $-\Delta_p u\ge u^q$  in $\Omega$, 
\begin{align}\label{I6} 
& \int_{B_R} u^\gamma \, dx  \le C \, R^{n- \tfrac{\gamma p}{q-p+1}}, & \qquad 0<\gamma< q,  \\ 
 & \int_{B_R} |\nabla u|^{\tfrac{\gamma p}{q+1}} \, dx  \le C \, R^{n-\tfrac {\gamma p}{q-p+1}}, 
&  \qquad 0<\gamma < q,  
 \label{I7} 
 \end{align}
 where  the constants $C$ in (\ref{I6}) and (\ref{I7}) 
 depend only on 
 $p, q, n, \gamma$. Note that  (\ref{I6})
   holds even for $\gamma=q$ (cf. \cite{MP}), while 
  (\ref{I7}) generally fails in this case.  In what follows, we 
 will substantially strengthen (\ref{I6}) in the end-point case $\gamma=q$, and 
 obtain global pointwise estimates of solutions.

  In \cite{PV}, we proved that 
 all compact sets $E\subset \Omega$ of zero Hausdorff measure, $H^{n-\tfrac{pq}{q-p+1}}(E) =0$, are removable 
 singularities for the equation $-\Delta_p u=u^q$, $q>q_*$, and a more general class of nonlinear equations.
 Earlier results of this kind,  under a stronger restriction 
 ${\rm cap}_{1, \, \tfrac {pq}{q-p+1}+\epsilon} (E) =0$ for some $\epsilon>0$,
  are due to Bidaut-V\'eron   \cite{BV3}. 
  Here ${\rm cap}_{1, \, s}(\cdot)$ is the capacity associated with the Sobolev space $W^{1, \, s}$.

  In fact, much more is true.  We will show below that a compact set $E\subset \Omega$ is a removable 
  singularity for $-\Delta_p u=u^q$ 
 if and only if it has  zero  fractional capacity: $\hbox{cap}_{p, \, \tfrac q {q-p+1}} \, (E) = 0$.  
 Here $\hbox{cap}_{\alpha, \, s}$ stands for  the Bessel capacity 
 associated with the Sobolev space $W^{\alpha, \, s}$ which is  defined in 
 Sec.~2.  We 
  observe that the usual $p$-capacity 
 $\hbox{cap}_{1,\,  p}$  used in the studies of the $p$-Laplacian \cite{HKM}, \cite{KM2} plays a secondary role 
 in the theory of equations of Lane--Emden type. Relations between these and other capacities 
 used in nonlinear PDE are discussed in \cite{AH}, \cite{M2}, and [V4].

Our characterization of removable singularities is based on the solution of the existence problem for the equation 
\begin{equation}\label{I5} 
 -\Delta_p u = u^q + \mu, \qquad u \ge 0,   
 \end{equation}
 with nonnegative measure $\mu$ obtained in Sec.~\ref{Om}. Main existence theorems 
 for quasilinear equations are stated below (Theorems \ref{T1} and \ref{T2}).
 Here we only 
  mention the following corollary in the case $\Omega=\RR^n$: 
  If (\ref{I5}) has an admissible solution $u$, then  
  \begin{equation}\label{I09a} 
 \int_{B_R} d \mu  \le C \, R^{n - \tfrac {pq}{q-p+1}},
 \end{equation}
 for every ball $B_R$ in $\RR^n$, where $C=C(p, q, n)$, provided 
 $1<p <n$ and $q>q_*$;  if $p \ge n$ or $q\le q_*$, then  $\mu=0$.
  
  Conversely, suppose that $1<p<n$, $q>q_*$, and $ d \mu = f\, dx$, $f \ge 0$, where 
 \begin{equation}\label{I10} 
\int_{B_R} f^{1+\epsilon} \, dx  \le C \, R^{n - \tfrac {(1+\epsilon) pq}{q-p+1}}, 
\end{equation}
 for some $\epsilon>0$. 
Then there exists a constant $C_0(p,q,n)$ such that (\ref{I5}) has 
an admissible solution on $\RR^n$  if 
 $C\le C_0(p,q,n)$. 
 
  The preceding inequality  is an analogue 
of the classical Fefferman--Phong condition \cite{Fef}, which appeared in applications to 
Schr\"odinger operators.  In particular,  (\ref{I10}) holds if $f \in L^{n(q-p+1)/pq, \, \infty}(\R^n)$. 
Here $L^{s, \, \infty}$ stands for the weak $L^s$ space. This sufficiency result, which  to the best of our 
knowledge is new even in the 
$L^s$ scale, provides a comprehensive solution to Problem 1 in \cite{BV2}. Notice 
that the exponent $s=\tfrac{n(q-p+1)}{pq}$  is sharp. Broader classes of  measures $\mu$ (possibly singular with respect 
to Lebesgue measure) which guarantee the existence of admissible solutions to (\ref{I5}) will be discussed in the sequel.

A substantial part of our work is concerned with  integral inequalities for 
nonlinear potential operators, which are at the heart of our approach. We employ the notion 
of Wolff's potential  introduced originally in  \cite{HW} 
in relation to the spectral synthesis 
problem for Sobolev spaces. For a nonnegative Borel measure $\mu$ on 
$\RR^n$, $s \in (1, \, +\infty)$, and $\alpha >0$, the Wolff potential 
${\rm\bf W}_{\alpha, \, s} \, \mu$ is 
defined by
 \begin{equation}\label{I05w}  
{\rm\bf W}_{\alpha, \, s} \, \mu(x)=\int_{0}^{\infty}\Big[\frac{\mu(B_{t}(x))}{t^{n-\alpha s}
}\Big]^{\frac{1}{s-1}} \frac{dt}{t}, \qquad x \in \RR^n. 
\end{equation}
We 
write  ${\rm\bf W}_{\alpha, \, s} \, f$ in place of ${\rm\bf W}_{\alpha, \, s} \, \mu$ if  $d \mu=f dx$, 
where  $f \in L^1_{\rm loc}(\RR^n)$, $f \ge 0$.   
When dealing with equations in a bounded domain $\Om\subset\RR^n$, a truncated version is useful: 
 \begin{equation}\label{I05t}
{\rm\bf W}_{\alpha, \, s}^{r} \, \mu(x)=\int_{0}^{r}\Big[\frac{\mu(B_{t}(x))}{t^{n-\alpha s}
}\Big]
^{\frac{1}{s-1}}\frac{dt}{t}, \qquad x \in \Om,
\end{equation}
where $0<r\leq 2 {\rm diam}  (\Omega)$. In many instances, it is more convenient to work with the dyadic version, 
also introduced in \cite{HW}:
 \begin{equation}\label{I05d}  
{\mathcal W}_{\alpha, \, s} \, \mu(x)=\sum_{Q \in \DD} \, \Big[\frac{\mu(Q)}{\ell(Q)^{n-\alpha s}
}\Big]^{\frac{1}{s-1}} \, \chi_Q (x), \qquad x \in \RR^n,
\end{equation}
where $\DD = \{ Q\}$ is the collection of the dyadic cubes 
$Q = 2^{i} (k + [0, \, 1)^n)$,  $i \in \ZZ, \, k \in \ZZ^n$, and $\ell(Q)$ is the side length of $Q$.  

An indispensable source on nonlinear potential theory is provided by \cite{AH},  
where   
the fundamental  Wolff inequality and its applications are discussed. Very recently, an
 analogue of Wolff's inequality for general dyadic and radially decreasing kernels was obtained in \cite{COV};  
 some of the tools developed there are employed below.

The dyadic Wolff potentials appear in the following  discrete model of  (\ref{I04}) studied in Sec.~3:  
\begin{equation}\label{Imod}  
u  = {\mathcal W}_{\alpha, \, s} \, u^q  + f,  \qquad u \ge 0. 
\end{equation}
 As it turns out, this nonlinear integral equation with $f = {\mathcal W}_{\alpha, \, s} \, \mu$
 is  intimately connected to 
the quasilinear differential equation (\ref{I5})  in the case  
$\alpha =1$, $s=p$, and to its $k$-Hessian counterpart  in the case $\alpha = \frac {2k} {k+1}$, $s= k+1$.  
Similar discrete models are used extensively 
in harmonic analysis and function spaces (see, e.g., \cite{NTV}, \cite{St2}, \cite{V1}). 

The profound role of Wolff's potentials in the theory of quasilinear equations was discovered 
by Kilpel\"ainen and Mal\'y \cite{KM1}. They established local 
pointwise estimates for nonnegative $p$-superharmonic functions in terms of Wolff's potentials 
of the associated $p$-Laplacian measure $\mu$. 
More precisely,  if $u\ge 0$ is a  $p$-superharmonic function in $B(x,3r)$ such that $-\Delta_p u = \mu $, then
\begin{equation}\label{IKM}
C_{1} \, {\rm\bf W}_{1, \, p}^{r} \, \mu(x)\leq u(x)\leq C_{2}\, \inf_{B(x,r)} \, u + C_{3} \, {\rm\bf W}_{1, \, p}^{2r} \, \mu(x),
\end{equation} 
where $C_{1}, C_{2}$ and $C_{3}$ are positive constants which depend only on $n$ and $p$.

 In \cite{TW1}, \cite{TW2}, Trudinger and Wang introduced the notion of 
the Hessian measure $\mu[u]$ associated with  $F_k [u]$ for a $k$-convex function $u$. 
Very recently, Labutin \cite{L}  
proved local pointwise estimates for Hessian equations 
analogous to (\ref{IKM}),  where the Wolff potential ${\rm\bf W}_{\frac {2k}{k+1}, \, k+1}^r \, \mu$ 
is used in place  of  ${\rm\bf W}_{1, \, p}^{r} \, \mu$.

In what follows, we will need   {\it global} pointwise estimates of this type. In the case of 
a $k$-convex solution to the equation $F_k[u] = \mu$ on $\RR^n$ such that $\inf_{ x \in \RR^n} \,(- u(x)) =0$, 
 one has 
\begin{equation}\label{glob}
  C_{1} \, {\rm\bf W}_{\frac {2k}{k+1}, \, k+1} \, \mu (x)  \le  - u(x) \le C_{2} \, 
  {\rm\bf W}_{\frac {2k}{k+1}, \, k+1} \, \mu (x),
 \end{equation} 
where $C_{1}$ and  $C_{2}$ are positive constants which depend only on $n$ and $k$. Analogous global 
estimates are obtained below for admissible solutions of the Dirichlet problem for $-\Delta_p u = \mu$ and 
$F_k[-u] = \mu$ in a bounded domain $\Omega \subset \RR^n$.

In the special case  $\Omega = \RR^n$, 
our criterion for the 
  solvability of (\ref{I5}) can be stated in the form of the pointwise condition 
  involving Wolff's potentials:
 \begin{equation}\label{I8} 
 {\rm\bf W}_{1, \, p}\left ( {\rm\bf W}_{1, \, p} \, \mu  \,   \right)^q(x) 
 \le C \, {\rm\bf W}_{1, \, p} \, \mu (x)<+\infty \quad 
{\rm a.e.},
 \end{equation}
 which is necessary with $C=C_1(p, q, n)$, and sufficient with another constant $C=C_2(p, q, n)$. Moreover, 
 in the latter case there exists an admissible solution $u$ to (\ref{I5})  such that 
 \begin{equation}\label{Ipointwise} 
 c_1 \, {\rm\bf W}_{1, \, p} \, \mu (x) \le u(x) \le c_2 \,  {\rm\bf W}_{1, \, p} \, \mu (x), \qquad x \in \RR^n,
  \end{equation}
  where $c_1$ and $c_2$ are positive constants which depend only on $p, q, n$, provided 
  $1<p<n$ and $q>q_*$; if $p \ge n$ or $q \le q_*$ 
  then $u=0$ and $\mu =0$.

The iterated Wolff potential condition (\ref{I8}) plays a crucial role 
in our approach. 
As we will demonstrate in Sec.~\ref{R^n}, it turns out to be equivalent to the fractional Riesz capacity condition 
\begin{equation}\label{I9} 
\mu(E) \le C \, \hbox{Cap}_{p, \frac q {q-p+1}} \, (E),
 \end{equation}
where $C$ does not depend on a compact set $E\subset \RR^n$. Such classes of measures $\mu$ 
were introduced by V. Maz'ya in the early 60-s in the framework of linear problems.

It follows that every admissible solution $u$ 
to (\ref{I5}) on $\RR^n$ obeys the inequality 
\begin{equation}\label{ineq} 
\int_E u^q \, dx  \le C \, \hbox{Cap}_{p, \frac q {q-p+1}} \, (E),
 \end{equation}
 for all compact sets $E\subset \RR^n$. We also prove an analogous estimate in a bounded domain 
 $\Omega$ (Sec.~\ref{Om}). Obviously, this yields  (\ref{I6}) in the end-point 
 case $\gamma=q$:   
   \begin{equation}\label{I6sharp} 
\int_{B_R} u^q \, dx  \le C \, R^{n- \tfrac{q p}{q-p+1}},    
 \end{equation}
 where $B_{2R} \subset \Omega$. In the critical case $q=q_*$, 
  we obtain an improved estimate: 
   \begin{equation}\label{I6crit} 
\int_{B_r} u^{q_*} \, dx  \le C \,  \left(\log (\tfrac {2R} {r})\right)^{\frac{1-p}{q-p+1}}, 
 \end{equation}
 for every ball $B_r$ of radius $r$ such that $B_r\subset B_R$, and $B_{2R} \subset \Omega$. 
 Certain Carleson measure inequalities are employed in the proof of (\ref{I6crit}).  
 We observe that (\ref{I6sharp}) and (\ref{I6crit})  
 yield Liouville-type theorems   for all admissible solutions to  
 (\ref{I5}) on $\RR^n$, or in exterior domains,  provided $q\le q_*$ (cf. \cite{BVP}, \cite{SZ}). 

Analogous results will be established  in Sec.~\ref{hessianequation} 
for  equations of Lane--Emden type involving the $k$-Hessian operator $F_k[u]$. 
We will prove that there exists a constant $C(k, q, n)$ 
such that, if 
\begin{equation}\label{I8H} 
 {\rm\bf W}_{\frac {2k}{k+1},\,  k+1}  ( {\rm\bf W}_{\frac {2k}{k+1}, \, k+1}   \mu )^q(x) \le C \, 
 {\rm\bf W}_{\frac {2k}{k+1}, \, k+1}  \mu (x)< +\infty \, \, 
 {\rm a.e.},
 \end{equation}
where $0\le C\le C(k, q, n)$, then the equation 
\begin{equation}\label{I12} 
 F_k[-u] = u^q + \mu, \qquad u \ge 0, 
\end{equation} 
has a  solution $u$ so that $-u$ is $k$-convex on $\RR^n$, and 
\begin{equation}\label{I13} 
c_1 \, {\rm\bf W}_{\frac {2k}{k+1}, \, k+1} \, \mu (x) \le u(x) \le c_2 \,  {\rm\bf W}_{\frac {2k}{k+1}, \, k+1} \, 
\mu (x), \qquad x \in \RR^n,
\end{equation} 
where $c_1$, $c_2$ are positive constants which depend only on $k, q, n$, 
for $1\le k < \frac n 2$. Conversely, (\ref{I8H}) is necessary 
in order that (\ref{I12}) have a  solution $u$ such that $-u$ is $k$-convex on $\RR^n$ provided  $1\le k < \frac n 2$ 
and $q>q_* =  \frac {nk}{n-2k}$; 
if $k\ge \frac n 2$ or $q \le q_*$ then  $u=0$ and $\mu =0$. 

In particular, 
(\ref{I8H}) holds if $d \mu = f \, dx$, where $f\ge 0$ and $f \in L^{n(q-k)/2kq, \, \infty}(\R^n)$;  
the exponent $\frac{n(q-k)}{2kq}$ is sharp. 

In Sec.~\ref{hessianequation}, we 
will obtain precise existence theorems for equation (\ref{I12}) in a bounded domain $\Omega$ 
with  the  Dirichlet 
boundary condition   $u = \phi$, $\phi \ge 0$, on $\partial \Omega$, for $1 \le k \le n$.  
Furthermore, removable singularities $E\subset \Omega$ for the homogeneous equation $F_k[-u] = u^q$, 
$u\ge 0$,  will be characterized  as the sets  of 
zero Bessel  capacity 
$\hbox{cap}_{2k, \, \frac q {q-k}} \, (E) = 0$,  in the most interesting case  $q> k$.

The notion of the {\it $k$-Hessian} capacity  introduced  by Trudinger and Wang proved 
to be very useful in studies of the uniqueness 
problem for $k$-Hessian equations  \cite{TW3}, as well as associated $k$-polar sets \cite{L}. Comparison theorems 
for  this capacity and the corresponding Hausdorff measures were obtained by Labutin in \cite{L} where it is 
proved that the $(n-2k)$-Hausdorff dimension is critical in this respect. 
We will enhance this result (see Theorem~\ref{capequiv} below)
 by showing that the $k$-Hessian capacity is in fact 
locally equivalent to the fractional 
Bessel capacity ${\rm cap}_{{\frac {2k}{k+1}}, \, k+1}$.

In conclusion, we remark that our methods 
provide a promising approach for a wide class of nonlinear problems,  
including curvature and subelliptic equations, and more general nonlinearities.

%**************************************************************
%**************************************************************

\section{Main results}\label{Main results}

Let $\Om$ be a bounded domain in $\RR^n$, $n\geq 2$. 
We study the 
existence problem for the quasilinear equation 
\begin{eqnarray}\label{Asuperharmonicom}
\left\{\begin{array}{c}
-{\rm div}\mathcal{A}(x,\nabla u)=u^{q}+\om,\\
u\geq 0 \quad {\rm in} \quad  \Om,\\
\hspace*{.15in}u=0 \quad {\rm on} \quad \partial \Om,
\end{array}
\right.
\end{eqnarray}
where $p>1$, $q>p-1$ and 
\begin{equation}\label{CA}
\mathcal{A}(x,\xi)\cdot\xi\geq \alpha\m{\xi}^p,\qquad
\m{\mathcal{A}(x,\xi)}\leq\beta\m{\xi}^{p-1}.
\end{equation}
for some $\alpha,\beta >0$.
The precise structural conditions imposed on  $\mathcal{A}(x,\xi)$ are stated in Sec.~4, 
formulae (\ref{2.1})--(\ref{2.6}). This includes the  principal model problem
\begin{eqnarray}\label{psuperharmonicom}
\left\{\begin{array}{c}
-\Delta_{p}u=u^{q}+\om,\\
u\geq 0 \quad {\rm in}\quad  \Om,\\
\hspace*{.15in} u=0 \quad {\rm on} \quad \partial \Om.
\end{array}
\right.
\end{eqnarray}
Here $\Delta_{p}$ is the $p$-Laplacian defined by $\Delta_{p}u= {\rm div}(\m{\nabla u}^{p-2}\nabla u).$ 
We observe that in the well-studied case $q \le p-1$ hard analysis techniques are not needed, and 
  many of our results simplify. We refer to \cite{Gre}, \cite{SZ} for further comments 
and references, especially in the classical case $q=p-1$. 

Our approach also applies to the following class of fully nonlinear equations 
\begin{equation}\label{ksubharmonicom}
\left\{\begin{array}{c}
F_{k}[-u]=u^{q}+\om , \\
u\geq 0 \quad{\rm in}\quad\Om,\\
\hspace*{.15in}u=\varphi \quad {\rm on} \quad\partial \Om,
\end{array}
\right.
\end{equation} 
 where $k=1,2,\dots,n$, and 
$F_{k}$ is  the $k$-Hessian operator,
$$F_{k}[u]=\sum_{1\leq i_{1}<\cdots<i_{k}\leq n}\lambda_{i_{1}}\cdots\lambda_{i_{k}}.$$ 
Here $(\lambda_{1},\dots,\lambda_{n})$ are the eigenvalues of the Hessian matrix 
$D^{2}u$, and $-u$ belongs to the class of $k$-subharmonic (or $k$-convex) functions on $\Om$ introduced by Trudinger and 
Wang in \cite{TW1}--\cite{TW2}.  Analogues of equations (\ref{Asuperharmonicom}) and (\ref{ksubharmonicom}) on the entire 
space $\RR^n$ are studied as well. \\
\indent To state our results, let us introduce some necessary definitions and notations.
Let $\mathcal{M}_{B}^{+}(\Om)$ (resp. ${\mathcal M}^{+}(\Om)$) denote the class of all 
nonnegative finite (respectively locally finite) Borel measures on $\Om$. For $\mu\in
{\mathcal M}^{+}(\RR^n)$ and a Borel set $E\subset\RR^n$, we denote by $\mu_{E}$ the restriction of 
$\mu$ to $E$: $d\mu_{E}=\chi_{E}d\mu$ where $\chi_{E}$ is the characteristic function 
of $E$. We define the Riesz potential ${\rm\bf I}_{\alpha}$ of order $\alpha$, 
$0<\alpha<n$, on $\RR^n$ by 
\begin{eqnarray*}
{\rm\bf I}_{\alpha}\mu(x)=c(n,\alpha)\int_{\RR^n}\m{x-y}^{\alpha-n}d\mu(y), \qquad x \in \RR^n, 
\end{eqnarray*}
where $\mu\in {\mathcal M}^{+}(\RR^n)$ and $c(n,\alpha)$ is a normalized constant.
For $\alpha>0$, $p>1$, such that $\alpha p<n$, the Wolff potential ${\rm\bf W}_{\alpha,\, p} \mu$ is 
defined by
$${\rm\bf W}_{\alpha,\, p}\mu(x)=\int_{0}^{\infty}\Big[\frac{\mu(B_{t}(x))}{t^{n-\alpha p}
}\Big]
^{\frac{1}{p-1}}\frac{dt}{t}, \qquad x \in \RR^n. 
$$ 
When dealing with equations in a bounded domain $\Om\subset
\RR^n$, it is convenient to use the truncated versions of Riesz and Wolff potentials.  For  
$0<r \le \infty$, $\alpha>0$ and $p>1$, we set 
$$
{\rm\bf I}_{\alpha}^{r}\mu(x)=\int_{0}^{r}\frac{\mu(B_{t}(x))}{t^{n-\alpha}}\frac{dt}{t}, \qquad  
{\rm\bf W}_{\alpha,\, p}^{r}\mu(x)=\int_{0}^{r}\Big[\frac{\mu(B_{t}(x))}{t^{n-\alpha p}
}\Big]
^{\frac{1}{p-1}}\frac{dt}{t}.$$
Here ${\rm\bf I}_{\alpha}^{\infty}$ and  ${\rm\bf W}_{\alpha,\, p}^{\infty}$ are understood as 
${\rm\bf I}_{\alpha}$ and  ${\rm\bf W}_{\alpha,\, p}$ respectively. 
For $\alpha>0$, we denote by 
${\rm\bf G}_{\alpha}$ the Bessel kernel of order $\alpha$ (see \cite{AH}, Sec. 1.2.4). The Bessel potential 
of a measure $\mu\in {\mathcal M}^{+}(\RR^n)$ is defined by 
\begin{eqnarray*}
{\rm\bf G}_{\alpha}\mu(x)=\int_{\RR^n}{\rm\bf G}_{\alpha}(x-y)d\mu(y), \qquad x \in \RR^n. 
\end{eqnarray*}
Various  capacities will be used throughout the paper. Among them are the Riesz 
and Bessel capacities defined respectively by
\begin{eqnarray*}
{\rm Cap}_{{\rm\bf I}_{\alpha},\, s}(E)=\inf\{\norm{f}_{L^{s}
(\RR^n)}^{s}: {\rm\bf I}_{\alpha}f\geq \chi_{E},~
 0\leq f\in L^{s}(\RR^n)\},
\end{eqnarray*}  
\begin{eqnarray*}
{\rm Cap}_{{\rm\bf G}_{\alpha},\, s}(E)=\inf\{\norm{f}_{L^{s}
(\RR^{n})}^{s}: {\rm\bf G}_{\alpha}f\geq \chi_{E},~
 0\leq f\in L^{s}(\RR^n)\}, 
\end{eqnarray*}  
for any $E\subset\RR^n$.\\
\indent Our first two theorems are concerned with {\it global} pointwise potential estimates for 
quasilinear  and Hessian equations on 
a bounded domain $\Om$ in $\RR^n$. 
\begin{theorem}\label{om-estimate-intro} Suppose that $u$ is a renormalized solution to
the equation 
\begin{eqnarray}\label{renor-intro}
\left\{\begin{array}{c}
-{\rm div}\mathcal{A}(x,\nabla u)=\om \quad {\rm in}\quad \Om,\\
\hspace{.18in} u=0 \quad {\rm on} \quad \partial \Om,
\end{array}
\right.
\end{eqnarray}
with data $\om\in {\mathcal M}_{B}^{+}(\Om)$. Then there is a positive constant $K$ which does not depend on 
$u$ and $\Om$ such that
\begin{eqnarray}\label{uestimate-intro}
\frac{1}{K} \,  {\rm\bf W}_{1,\, p}^{\frac{{\rm dist}(x,\partial\Om)}{3}}\om(x) \leq u(x)\leq K \, {\rm\bf W}_{1,\, p}^
{2{\rm diam}(\Om)}\om(x),  
\end{eqnarray}
for all $x$ in $\Om$.
\end{theorem}
\begin{theorem}\label{globalH-intro} Let $\om$ be a nonnegative finite 
measure on $\Om$ such that $\om\in L^{s}(\Om\setminus E)$ for 
a compact set $E\subset\Om$. Here 
$s>\frac{n}{2k}$ if $1\leq k\leq \frac n 2$, and $s=1$ if $\frac n 2<k\leq n$. 
Suppose that $-u$ is a nonpositive $k$-subharmonic function in  $\Om$ such that
$u$ is continuous near $\partial\Om$, and solves the equation
\begin{eqnarray*}
\left\{\begin{array}{c}
F_{k}[-u]=\om \quad {\rm in}\quad \Om,\\
u=0 \quad {\rm on}\quad \partial \Om,
\end{array}
\right.
\end{eqnarray*}
Then, for all $x\in\Om$,
\begin{equation}\label{pointwiseH-intro} \frac{1}{K}\, {\rm\bf W}_{\frac{2k}{k+1},\, k+1}^{\frac{{\rm dist}
(x,\partial\Om)}{8}}\om(x) 
\leq  u(x) \leq K \, {\rm\bf W}_{\frac{2k}{k+1},\, k+1}^{2{\rm diam}(\Om)}\om(x),
\end{equation}
where $K$ is a constant which does not depend on $x$, $u$, and $\Om$.
\end{theorem}

 We remark that the upper estimate 
in (\ref{uestimate-intro}) does not hold in general if $u$ is merely a weak solution 
of (\ref{renor-intro}) in the sense of \cite{KM1}. For a counter example, 
see \cite[Sec. 2]{Kil}. Upper estimates similar to the one 
in (\ref{pointwiseH-intro}) hold also for  $k$-subharmonic
functions with non-homogeneous boundary condition as well (see Sec. \ref{hessianequation}). Equivalent definitions of renormalized solutions
to the problem (\ref{renor-intro}) are given in Sec. \ref{Om}. For definitions
of $k$-subharmonic  functions, see Sec. \ref{hessianequation}.  

Note also that in the  case of the entire space $\Om=\RR^n$, 
if $-u$ is a non-positive $k$-subharmonic function such that 
$F_{k}[-u]=\mu$ and $\inf_{x\in \RR^n} u(x)=0$, then  
 \begin{equation*}
 \frac 1 K \, {\rm\bf W}_{\frac{2k}{k+1},\, k+1}\mu(x) 
 \le u(x) \le K \, {\rm\bf W}_{\frac{2k}{k+1},\, k+1}\mu(x).
\end{equation*} 
An analogous two-sided estimate holds for  $\mathcal{A}$-superharmonic 
functions as well, with ${\rm\bf W}_{1,\, p}\mu$ in place of ${\rm\bf W}_{\frac{2k}
{k+1}, \, k+1}\mu$. These global estimates  are deduced from  the local ones given
in \cite{L}, \cite{KM2}.  

In the next two theorems we give criteria for the solvability of quasilinear 
and Hessian equations on the entire space $\RR^n$.    
\begin{theorem}\label{T1}
Let $\om$ be a measure in $\mathcal{M}^{+}(\RR^n)$, $1<p<n$ and $q>p-1$. 
Then the following statements are equivalent.\\
{\rm(i)} There exists a nonnegative $\mathcal{A}$-superharmonic solution $u\in 
L^{q}_{\rm loc}(\RR^n)$
 to the equation 
\begin{equation}
\label{T1ep-equation}
\left\{\begin{array}{c}
\inf_{x\in\RR^{n}}u(x)=0\\ 
-{\rm div}\mathcal{A}(x, \nabla u)=u^{q}+\ep \, \om \quad{\rm in}\quad \RR^n, 
\end{array}
\right.
\end{equation} 
for some $\ep>0$.\\
{\rm(ii)} For all compact sets $E\subset\RR^n$, 
\begin{equation}\label{E1}
\om(E)\leq C \, {\rm Cap}_{{\rm\bf I}_{p},\, \frac{q}{q-p+1}}(E).
\end{equation}
{\rm(iii)} The testing inequality 
\begin{equation}\label{E2}
\int_{B}\Big[{\rm\bf W}_{1,\, p}\om_{B}(x)\Big]^{q}dx\leq C \, \om(B)
\end{equation}
holds for all balls $B$ in $\RR^{n}$ .\\
{\rm(iv)} There exists a constant C such that
\begin{equation}
\label{T1pointwise}
{\rm\bf W}_{1,\, p}({\rm\bf W}_{1,\, p}\om)^{q} (x) \leq C \, {\rm\bf W}_{1,\, p}\om (x)< \infty \quad a.e.
\end{equation}
Moreover, there is a constant $C_{0}=C_{0}(n,p,q,\alpha,\beta)$ such that if any one of the conditions 
(\ref{E1})--(\ref{T1pointwise}) holds with 
$C\leq C_{0}$, then equation (\ref{T1ep-equation}) has a solution $u$ with $\ep=1$ which satisfies the 
two-sided estimate
\begin{equation*}
\frac 1 K \, {\bf W}_{1,\, p}\om(x)\leq u(x)\leq K \, {\bf W}_{1,\, p}\om(x), \qquad x\in \RR^n, 
\end{equation*}  
where $K$ depends only on $n, p, q, \alpha,\beta$. Conversely, if (\ref{T1ep-equation}) has a solution $u$ as in statement {\rm(i)} with $\ep=1$, then 
conditions (\ref{E1})--(\ref{T1pointwise}) hold with $C=C_1(n, p, q,\alpha, \beta)$. Here 
$\alpha$ and $\beta$ are the structural constants of $\mathcal{A}$ defined in 
(\ref{CA}).
\end{theorem}
Using  condition (\ref{E1}) in the above theorem, we can now deduce from the 
isoperimetric inequality: 
$$\m{E}^{1-\frac{pq}{q-p+1}}\leq C {\rm Cap}_{{\rm \bf I}_{p},\, 
\frac{q}{q-p+1}}(E), $$
(see \cite{AH} or \cite{M2}), a simple sufficient condition for the solvability of (\ref{T1ep-equation}).  
\begin{corollary}\label{iso} Suppose that $f\in L^{\frac{n(q-p+1)}{pq},\, \infty}(\RR^n)$ and 
$d\om=fdx$. If $q>p-1$ and $\frac{pq}{q-p+1}< n$, then equation (\ref{T1ep-equation}) has a nonnegative solution for some 
$\epsilon>0$.
\end{corollary}
\begin{remark}{\rm The condition $f\in L^{\frac{n(q-p+1)}{pq},\, \infty}(\RR^n)$ in 
 Corollary \ref{iso} can  be relaxed by using the Fefferman--Phong condition 
\cite{Fef}:
$$\int_{B_{R}}f^{1+\delta}dx\leq C R^{n-\frac{(1+\delta)pq}{q-p+1}},$$
for some $\delta>0$, which is known to be sufficient for the validity of (\ref{E1}); see, e.g., 
\cite{KS}, \cite{V2}. 
}\end{remark}
\begin{theorem}\label{T3}
Let $\om$ be a measure in $\mathcal{M}^{+}(\RR^n)$, $1\leq k< \frac n 2$, 
and  $q>k$. Then the 
following statements are equivalent.\\
{\rm(i)} There exists a  solution $u\geq 0$, $-u\in\Phi^{k}(\Om)\cap 
L^{q}_{\rm loc}(\RR^n)$, to the equation 
\begin{equation}
\label{ep-equationHT3}
\left\{\begin{array}{c}
\inf_{x\in\RR^{n}}u(x)=0\\ 
F_{k}[-u]=u^{q}+\ep \, \om \quad{\rm in}\quad \RR^n, 
\end{array}
\right.
\end{equation} 
for some $\ep>0$.\\
{\rm(ii)} For all compact sets $E\subset\RR^n$, 
\begin{equation}\label{TE1}
\om(E)\leq C \, {\rm Cap}_{{\rm\bf I}_{2k},\, \frac{q}{q-k}}(E).
\end{equation}
\noindent {\rm(iii)} The testing inequality 
\begin{equation}\label{TE2}
\int_{B}\Big[{\rm\bf W}_{\frac{2k}{k+1},\, k+1}\om_{B}(x)\Big]^{q} \, dx\leq C \, \om(B)
\end{equation}
holds for all balls $B$ in $\RR^{n}$ .\\
{\rm(iv)} There exists a constant C such that
\begin{equation}
\label{pointwiseHT3}
{\rm\bf W}_{\frac{2k}{k+1},\, k+1}({\rm\bf W}_{\frac{2k}{k+1},\,k+1}\om)^{q}(x)\leq C \, {\rm\bf W}_{\frac{2k}{k+1},\,
k+1}\om(x) < \infty \quad a.e.
\end{equation}
Moreover, there is a constant $C_{0}=C_{0}(n,k,q)$ such that 
if any one of the conditions (\ref{TE1})--(\ref{pointwiseHT3}) holds  with  $C\leq C_{0}$, then equation (\ref{ep-equationHT3}) has a solution $u$ with $\ep=1$ which 
satisfies  the two-sided estimate
\begin{equation*}
\frac 1 K \, {\bf W}_{\frac{2k}{k+1},\, k+1}\om(x)\leq u(x)\leq K \, {\bf W}_{\frac{2k}{k+1},\, k+1}\om(x), 
\qquad x \in \RR^n,
\end{equation*}  
where $K$ depends only on $n, k, q$. Conversely, if there is a solution $u$ to  (\ref{ep-equationHT3})   as in statement {\rm(i)} with 
$\epsilon=1$, then conditions (\ref{TE1})--(\ref{pointwiseHT3}) hold with $C=C_{1}(n,k,q)$. 
\end{theorem}

\begin{corollary} Suppose that $f\in L^{\frac{n(q-k)}{2kq},\, \infty}(\RR^n)$ and 
$d\om=fdx$. If $q>k$ and $\frac{2kq}{q-k}<n$ then the equation 
(\ref{ep-equationHT3}) has a nonnegative solution for some 
$\epsilon>0$.
\end{corollary}

Since  ${\rm Cap}_{I_{\alpha},\, s}(E)=0$  in the case $\alpha \, s\geq n$ for 
all Borel sets $E\subset\RR^n$ (see \cite{AH}), we obtain the following Liouville-type theorems for quasilinear and Hessian differential inequalities. 
\begin{corollary} If $q\leq \frac{n(p-1)}{n-p}$, then the inequality
$-{\rm div}\mathcal{A}(x,\nabla u)\geq u^q$
 admits no nontrivial nonnegative $\mathcal{A}$-superharmonic solutions 
in $\RR^n$. 
Analogously, if $q\leq\frac{nk}{n-2k}$, then the inequality $F_{k}[-u]\geq u^q $ 
admits no nontrivial nonnegative solutions in $\RR^n$. 
\end{corollary}
\begin{remark}{\rm When $1<p<n$ and $q>\frac{n(p-1)}{n-p}$, 
the function $u(x)=c\m{x}^{\frac{-p}{q-p+1}}$ with 
$$c=\Big[\frac{p^{p-1}}{(q-p+1)^{p}}\Big]^{\frac{1}{q-p+1}}
[q(n-p)-n(p-1)]^{\frac{1}{q-p+1}},$$
is a nontrivial admissible (but singular) global solution  of $-\Delta_p u=u^q$ 
(see \cite{SZ}).
Similarly, the function $u(x)=c'\m{x}^{\frac{-2k}{q-k}}$ with
$$c'=\Big[\frac{(n-1)!}{k!(n-k)!}\Big]^{
\frac{1}{q-k}}\Big[\frac{(2k)^{k}}{(q-k)^{k+1}}\Big]^{\frac{1}{q-k}}
[q(n-2k)-nk]^{\frac{1}{q-k}},$$
where $1\leq k <n/2$ and $q>\frac{nk}{n-2k}$, is a singular admissible global solution
of $F_{k}[-u]=u^{q}$ (see \cite{Tso} or \cite{Tru1}, formula (3.2)).
 Thus, we see that 
the exponent $\frac{n(p-1)}{n-p}$ (respectively $\frac{nk}{n-2k}$) is also critical 
for the homogeneous equation $-{\rm div}\mathcal{A}(x,\nabla u)=u^q$ (respectively 
$F_{k}[-u]=u^q$) in $\RR^n$. The situation is different when we restrict ourselves 
only to  {\it locally bounded} solutions in $\RR^n$
(see \cite{GS}, \cite{SZ}).
}\end{remark}

Existence results on a bounded domain $\Om$ analogous to Theorems \ref{T1} and 
\ref{T3} are contained in the following two theorems, where Bessel potentials and the corresponding capacities 
are used in place of respectively Riesz potentials and Riesz capacities.  
\begin{theorem}\label{T2}
Let $\om$ be a measure in ${\mathcal M}_{B}^+(\Om)$ which is compactly  supported 
in $\Om$. Let $p>1$, $q>p-1$, and let $R={\rm diam}(\Om)$. Then the following statements are 
equivalent.\\
{\rm(i)} There exists a nonnegative renormalized solution $u\in 
L^{q}(\Om)$
to the equation 
\begin{equation}
\label{T2ep-equation}
\left\{\begin{array}{c}
-{\rm div}\mathcal{A}(x, \nabla u)=u^{q}+\ep \, \om \quad{\rm in}\quad \Om,\\
u=0 \hspace*{.3in}{\rm on}\quad \partial\Om, 
\end{array}
\right.
\end{equation} 
for some $\ep>0$.\\
{\rm(ii)} For all compact sets $E\subset{\rm supp}\om$, 
\begin{equation*}
\om(E)\leq C \, {\rm Cap}_{{\rm\bf G}_{p},\, \frac{q}{q-p+1}}(E).
\end{equation*}
{\rm(iii)} The testing inequality 
\begin{equation*}
\int_{B}\Big[{\rm\bf W}_{1,\,p}^{2R}\om_{B}(x)\Big]^{q} \, dx\leq C \, \om(B)
\end{equation*}
holds for all balls $B$  such that $B\cap{\rm supp}\om\not =\emptyset$ .\\
{\rm(iv)} There exists a constant C such that
\begin{equation}
\label{T2pointwise2}
{\rm\bf W}_{1,\, p}^{2R}({\rm\bf W}_{1,\, p}^{2R}\om)^{q}(x)\leq C  \, {\rm\bf W}_{1,\, p}^{2R}
\om(x) < \infty \quad {\rm a.e. ~on~} \Om.
\end{equation}
\end{theorem}
\begin{remark}\label{subcritical} {\rm In the case where $\om$ is not compactly supported  in $\Om$, it can 
be easily seen from the proof of this theorem  that any one of the conditions
{\rm (ii)}, {\rm (iii)}, and {\rm (iv)} above is still {\it sufficient} for the solvability
of (\ref{T2ep-equation}) for some $\epsilon>0$. Moreover, in the 
subcritical case $\frac{pq}{q-p+1}>n$, these conditions are redundant since  
the Bessel capacity  ${\rm Cap}_{{\rm\bf G}_{p},\, \frac{q}{q-p+1}}$ of a single
point is  positive (see \cite{AH}). This  ensures that statement {\rm(ii)} of Theorem \ref{T2}
 holds  for 
some constant $C>0$ provided $\om$ is a finite measure.} 
\end{remark}
\begin{corollary} 
Suppose that $f\in L^{\frac{n(q-p+1)}{pq},\, \infty}(\Om)$ and 
$d\om=fdx$. If $q>p-1$ and $\frac{pq}{q-p+1}< n$ then the equation 
(\ref{T2ep-equation}) has a nonnegative renormalized (or equivalently, entropy) 
solution for some $\epsilon>0$.
\end{corollary}
\begin{theorem}\label{maintheorem4T4} 
Let $\Om$ be a uniformly $(k-1)$-convex domain in $\RR^n$, and let $\om\in {\mathcal M}^{+}_
{B}(\Om)$ be compactly supported  in $\Om$. Suppose that $1 \le k \le n$, 
$q>k$, $R={\rm diam}(\Om)$, and $\varphi\in C^{0}(\partial \Om)$, $\varphi \geq 0$. 
Then the following statements
are equivalent.\\
{\rm (i)} There exists a solution $u\geq 0$, $-u\in\Phi^{k}(\Om)\cap L^{q}(\Om)$,
continuous near $\partial\Om$,  to the equation
\begin{equation}\label{khessianep}
\left\{\begin{array}{c}
F_{k}[-u]=u^{q}+\epsilon \, \om \quad {\rm in}\quad \Om, \\
u=\epsilon \, \varphi \quad {\rm on}\quad \partial \Om
\end{array}
\right.
\end{equation} 
for some $\epsilon>0$.\\
{\rm (ii)} For all compact sets $E\subset{\rm supp}\om$,
$$\om(E)\leq C \, {\rm Cap}_{{\rm\bf G}_{2k},\, \frac{q}{q-k}}(E).$$
{\rm (iii)} The testing inequality 
\begin{equation*}
\int_{B}\Big[{\rm\bf W}_{\frac{2k}{k+1},\, k+1}^{2R}\om_{B}(x)\Big]^{q} \, dx\leq C \, \om(B)
\end{equation*}
holds for all balls $B$ such that  $B\cap{\rm supp}\om\not =\emptyset$ .\\
{\rm (iv)} There exists a constant C such that
\begin{equation*}
{\rm\bf W}_{\frac{2k}{k+1},\,k+1}^{2R}({\rm\bf W}_{\frac{2k}{k+1},\,k+1}^{2R}\om)^{q}(x) \leq 
C \, {\rm\bf W}_{\frac{2k}{k+1},\,k+1}^{2R}
\om (x) < \infty \quad {\rm a.e. ~on~} \Om.
\end{equation*}
\end{theorem}
\begin{remark}{\rm As in Remark \ref{subcritical}, suppose that $\om\in\mathcal{M}_{B}
^{+}(\Om)\cap L^{s}(\Om\setminus E)$, for a compact set $E\subset\Om$, where
$s>\frac{n}{2k}$ if $k\leq \frac n 2$, and $s=1$ if $k> \frac n 2$.  Then 
any one of the conditions
{\rm (ii)}, {\rm (iii)}, and {\rm (iv)} in Theorem \ref{maintheorem4T4} is still sufficient for 
the solvability of (\ref{khessianep}) for some $\epsilon>0$. Moreover, in the 
subcritical case $\frac{2kq}{q-k}>n$ these conditions are redundant.
}\end{remark}
\begin{corollary} 
Suppose that $f\in L^{\frac{n(q-k)}{2kq},\, \infty}(\Om)\cap L^{s}(\Om\setminus E)$ for
some $s>n/2k$ and for some compact set $E\subset\Om$. Let 
$d\om=fdx$. If $q>k$ and $\frac{2kq}{q-k}< n$ then the equation 
(\ref{khessianep}) has a nonnegative solution for some 
$\epsilon>0$.
\end{corollary}

Our results on {\it local integral estimates} for quasilinear and Hessian inequalities 
are given in the next two theorems.  We will need the capacity associated with the space 
$W^{\alpha, \, s}$ relative to the domain $\Omega$ defined by
\begin{equation}\label{capomega}
{\rm cap}_{\alpha,\, s}(E,\Om)=\inf\{\norm{f}^{s}_
{W^{\alpha,\,s}(\RR^n)}: \, f\in C_{0}^{\infty}(\Om), f\geq 1 {\rm ~on~} E\}.
\end{equation}

\begin{theorem}\label{localest1}  Let $u$ be a nonnegative $\mathcal{A}$-superharmonic function 
in $\Om$ such that $-{\rm div}\mathcal{A}(x,\nabla u)\geq u^q$, where $p>1$ and 
$q>p-1$. Then there exists a constant $C$ which depends only on $p, q, n$, 
and the structural constants of $\mathcal{A}$  
such that  
$$\int_{B_{R}}u^{q} \, dx\leq C \, R^{n-\frac{pq}{q-p+1}} $$
if $\frac{pq}{q-p+1}<n$, and 
$$\int_{B_{r}}u^{q} \, dx\leq C \, (\log\tfrac{2R}{r})^{\frac{1-p}{q-p+1}} $$
if $\frac{pq}{q-p+1}=n$. Here $0<r\leq R$, $B_r \subset B_R$,  and $B_{R}$ is a ball such that
$B_{2R}\subset\Om$.  

Moreover, if $\frac{pq}{q-p+1}<n$, and  $\Om$ is a bounded $C^{\infty}$-domain then
$$\int_{E}u^q\leq C \, {\rm cap}_{p,\, \frac{q}{q-p+1}}(E,\Om)$$
for any compact set $E\subset\Om$. 
\end{theorem}

\begin{theorem}\label{localest2}  Let $u\geq 0$ be such that  $-u$ is $k$-subharmonic 
and that $F_{k}[-u]\geq u^{q}$ in $\Om$ with $q>k$. 
Then there exists a constant $C$ which depends only on $k, q, n$ such that 
$$\int_{B_{R}}u^{q} \, dx\leq C \, R^{n-\frac{2kq}{q-k}} $$
if $\frac{2kq}{q-k}<n$, and 
$$\int_{B_{r}}u^{q} \, dx\leq C \, (\log\tfrac{2R}{r})^{\frac{-k}{q-k}} $$
if $\frac{2kq}{q-k}=n$. Here $0<r\leq R$, $B_r \subset B_R$,  and $B_{R}$ is a ball such that
$B_{2R}\subset\Om$.  

Moreover, if $\frac{2kq}{q-k}<n$ and $\Om$ is 
a bounded $C^{\infty}$-domain then   
$$\int_{E}u^q\leq C \, {\rm cap}_{2k,\, \frac{q}{q-k}}(E,\Om)$$
for any compact set $E\subset\Om$.
\end{theorem}

As a consequence of Theorems \ref{T2} and \ref{maintheorem4T4}, we deduce the 
following results concerning removable singularities of quasilinear and 
fully nonlinear equations.

\begin{theorem}\label{removeforp} Let $E$ be a compact subset of $\Om$. 
Then any solution $u$ to the problem 
\begin{equation*}
\left\{\begin{array}{c}
 u {\rm ~is~} \mathcal{A}{\text-}{\rm superharmonic~in~}\Om\setminus E,\\
u\in L^{q}_{\rm loc}(\Om\setminus E), \quad u \ge 0,\\ 
-{\rm div}\mathcal{A}(x, \nabla u)=u^{q} \quad {\rm in} \quad\mathcal{D}'(\Om\setminus E),
\end{array}
\right.
\end{equation*} 
is also a solution to 
\begin{equation*}
\left\{\begin{array}{c}
 u {\rm~is~} \mathcal{A}{\text-}{\rm superharmonic~in~}\Om,\\
u\in L^{q}_{\rm loc}(\Om), \quad u \ge 0, \\
-{\rm div}\mathcal{A}(x, \nabla u)=u^{q} \quad {\rm in}\quad \mathcal{D}'(\Om),
\end{array}
\right.
\end{equation*}
if and only if ${\rm Cap}_{{\rm\bf G}_{p},\,\frac{q}{q-p+1}}(E)=0$. 
\end{theorem}
   
\begin{theorem}\label{removefork} Let $E$ be a compact subset of $\Om$.  
 Then any solution $u$ to the problem 
\begin{equation*}
\left\{\begin{array}{c}
-u {\rm~is~} k{\text-}{\rm subharmonic~in~}\Om\setminus E,  \\
u\in L^{q}_{\rm loc}(\Om\setminus E), \quad u \ge 0,\\
F_{k}[-u]=u^{q} \quad {\rm in}\quad \mathcal{D}'(\Om\setminus E),
\end{array}
\right.
\end{equation*} 
is also a solution to 
\begin{equation*}
\left\{\begin{array}{c}
 -u {\rm~is~} k{\text-}{\rm subharmonic~in~}\Om,  \\
u\in L^{q}_{\rm loc}(\Om), \quad u \ge 0,\\
F_{k}[-u]=u^{q} \quad {\rm in}\quad \mathcal{D}'(\Om),
\end{array}
\right.
\end{equation*}
if and only if ${\rm Cap}_{{\rm\bf G}_{2k},\,\frac{q}{q-k}}(E)=0$. 
\end{theorem}

 In \cite{TW3}, Trudinger and Wang introduced the so called {\it $k$-Hessian}
capacity ${\rm cap}_{k}(\cdot,\Om)$ defined by
$${\rm cap}_{k}(E, \Om)=\sup\Big\{ \int_{E}F_{k}[u]: u {\rm ~is~} k{\text -}{\rm subharmonic
~in~} \Om, -1<u<0 \Big\}.$$
 Our next theorem  asserts 
that locally 
the $k$-Hessian capacity is equivalent to the Bessel capacity ${\rm Cap}_{{\rm\bf G}_{\frac{2k}{k+1}},\,k+1}$. 

In what follows, $\mathcal{Q}=\{Q\}$ will 
stand for a Whitney decomposition of $\Om$ into a union of disjoint dyadic cubes (see Sec.~\ref{Om}). 
\begin{theorem}\label{capequiv}
Let $1\leq k< \frac n 2$ be an integer. Then there are constants $M_{1}$, $M_{2}$ such that
\begin{equation*}
M_{1} \, {\rm Cap}_{{\rm\bf G}_{\frac{2k}{k+1}},\,k+1}(E)\leq {\rm cap}_{k}(E,\Om) \leq M_{2} \, 
{\rm Cap}_{{\rm\bf G}_{\frac{2k}{k+1}},\,k+1}(E), 
\end{equation*}
for any compact set $E\subset\overline{Q}$ with $Q\in \mathcal{Q}$. Furthermore,
if $\Om$ is a bounded $C^{\infty}$-domain then
$${\rm cap}_{k}(E,\Om)\leq C \, {\rm cap}_{\frac{2k}{k+1},
\,k+1}(E,\Om),$$
 for any compact set $E\subset\Om$, where ${\rm cap}_{\frac{2k}{k+1},
\,k+1}(E,\Om)$ is defined by (\ref{capomega}) with $\alpha = 2k$ and $s = \frac{2k}{k+1}$. 
\end{theorem}

%\begin{remark}{\rm It follows from (\ref{T1cap}) and (\ref{T2cap}) that 
%$q=\frac{n(p-1)}{n-p}$ is the critical exponent for the 
%solvability of (\ref{Aom}). In fact, if $q\leq\frac{n(p-1)}{n-p}$
%then ${\rm Cap}_{p, \frac{q}{q-p+1}}
%(E)=0$ for all $E\subset\RR^n$ (see \cite{AH, M2}), and hence (\ref{T1ep-equation}) has no global
%solutions on $\RR^n$ provided $\om\not = 0$. On the other hand, if $q<\frac{n(p-1)}{n-p}$ then 
%${\rm cap}_{p, \frac{q}{q-p+1}}(E)>0$ for all non-empty subsets $E$ of $\RR^n$. Therefore, 
%the equation (\ref{T2ep-equation}) always has a solution for some $\epsilon>0$.  
%See also the Introduction in \cite{SZ} for some other related 
%critical exponents for differential equations or inequalities involving p-Laplacian.      
%}\end{remark}

%***********************************************************************
%***********************************************************************
\section{Discrete models of nonlinear equations}
Let $\mathcal{D}$ be the family of all dyadic cubes in $\RR^n$. For $\om\in {\mathcal 
M}^{+}(\RR^n)$ we define its dyadic Riesz and Wolff potentials respectively by
\begin{eqnarray}\label{3.1a}
\mathcal{I}_{\alpha}\om(x)&=&\sum_{Q\in\mathcal{D}}\frac{\om(Q)}{\m{Q}^{1-\alpha/n}}\chi_{Q}(x), 
\end{eqnarray}
\begin{eqnarray}\label{3.1}
\mathcal{W}_{\alpha,\,p}\om(x)&=&\sum_{Q\in\mathcal{D}}\Big[\frac{\om(Q)}{\m{Q}^{1-\alpha p/n}}\Big]
^{\frac{1}{p-1}}\chi_{Q}(x).
\end{eqnarray}
In this section we are concerned with nonlinear inhomogeneous integral equations of the 
type  
\begin{equation}\label{3.2}
u=\mathcal{W}_{\alpha,\, p}(u^q)+f, \hspace*{.3in} u\in L^{q}_{\rm loc}(\RR^n),
\end{equation}
where $f\in L^{q}_{\rm loc}(\RR^n)$, $q>p-1$, and $\mathcal{W}_
{\alpha,\,p}$ is defined as in (\ref{3.1}) with $\alpha>0$ and $p>1$ such that 
$0<\alpha p<n$. \\
\indent It is convenient to introduce a nonlinear operator $\mathcal{N}$ 
associated with the equation (\ref{3.2}) defined by 
\begin{equation}\label{operatorA}
\mathcal{N}f=\mathcal{W}_{\alpha,\,p}
(f^q), \qquad f\in L_{\rm loc}^{q}(\RR^n),
\end{equation}
so that (\ref{3.2}) can be rewritten as
\begin{equation*}
u=\mathcal{N}u + f,\hspace*{.3in} u\in L_{\rm loc}^{q}(\RR^n).
\end{equation*}
It is obvious that $\mathcal{N}$ is monotonic, i.e., $\mathcal{N}f\geq \mathcal{N}g$ 
whenever $f\geq g\geq 0$ a.e., and  $\mathcal{N}(\lambda f)=\lambda^{\frac{q}{p-1}}
\mathcal{N}f$ for all $\lambda\geq 0$. Since 
\begin{equation}\label{c(p)}
(a+b)^{p'-1}\leq \max\{1,2^{p'-2}\}(a^{p'-1}+b^{p'-1})
\end{equation}
for all $a, b \ge 0$, it follows that
\begin{equation}\label{propertyA}
\Big[\mathcal{N}(f+g)\Big]^{1/q}\leq \max\{1,2^{p'-2}\}\Big[(\mathcal{N}f)^{1/q}
+(\mathcal{N}g)^{1/q}\Big].
\end{equation}

\begin{proposition}\label{discrete0t}Let $\mu$ be a measure in $M^{+}(\RR^n)$, 
$\alpha>0$, $p>1$ and $q>p-1$.
Then the following quantities are equivalent:
\begin{eqnarray*}
&(a)& A_{1}(P,\mu)=\sum_{Q\subset P} \Big[\frac{\mu(Q)}{\m{Q}^{1-
\frac{\alpha p}{n}}}\Big]^{\frac{q}{p-1}}\m{Q},\\
&(b)& A_{2}(P,\mu)=\int_{P}\Big[\sum_{Q\subset P} \frac{\mu(Q)^{\frac{1}{p-1}}}
{\m{Q}^{(1-\frac{\alpha p}{n})\frac{1}{p-1}}}\chi_{Q}(x)\Big]^{q}dx,\\
&(c)& A_{3}(P,\mu)=\int_{P}\Big[\sum_{Q\subset P}\frac{\mu(Q)}{\m{Q}^{1-\frac
{\alpha p}{n}}}\chi_{Q}(x)\Big]^{\frac{q}{p-1}}dx,
\end{eqnarray*}
where $P$ is a dyadic cube in $\RR^n$ or $P=\RR^n$ and the constants of equivalence 
do not depend on $P$ and $\mu$.  
\end{proposition}
\begin{proof} The equivalence of $A_{1}$ and $A_{3}$ follows from Wolff's 
inequality (see \cite{HW}, \cite{COV}). Moreover, it has been proven in \cite{COV}
that 
\begin{equation}\label{max}
A_{3}(P,\mu)\asymp \int_{P}\Big[\sup_{x\in Q\subset P}\frac{\mu(Q)}{\m{Q}^
{1-\alpha p/n}}\Big]^{\frac{q}{p-1}}dx.
\end{equation}
Since 
$$\Big[\sup_{x\in Q\subset P}\frac{\mu(Q)}{\m{Q}^
{1-\alpha p/n}}\Big]^{\frac{1}{p-1}} \leq \sum_{Q\subset P} \frac{\mu(Q)^
{\frac{1}{p-1}}}{\m{Q}^{(1-\frac{\alpha p}{n})\frac{1}{p-1}}}\chi_{Q}(x),$$
from (\ref{max}) we obtain 
$A_{3}\leq C A_{2}$. 
In addition, for $p\leq 2$ we clearly have $A_{2}\leq A_{3}\leq C A_{1}$. Therefore, 
it remains to check that, in the case $p>2$, $A_{2}\leq C A_{1}$ for some $C>0$
independent of $P$ and $\mu$.
By Proposition 2.2 in \cite{COV} we have (note that $q>p-1>1$)   
\begin{eqnarray}
\label{A_{2}}
\hspace{1cm} A_{2}(P,\mu)&=&\int_{P}\Big[\sum_{Q\subset P} \frac{\mu(Q)^{\frac
{1}{p-1}}}{\m{Q}^{(1-\frac{\alpha p}{n})\frac{1}{p-1}}}\chi_{Q}(x)\Big]^{q}dx\\
&\leq& C\sum_{Q\subset P}\frac{\mu(Q)^\frac{1}{p-1}}{\m{Q}^{(1-\frac{\alpha p}
{n})\frac{1}{p-1}+q-2}} \Big[\sum_{Q'\subset Q} \frac{\mu(Q')^{\frac{1}{p-1}}}
{\m{Q'}^{(1-\frac{\alpha p}{n})
\frac{1}{p-1}-1}}\Big]^{q-1}.\nonumber
\end{eqnarray}
On the other hand, by H\"older's inequality,
\begin{eqnarray*}
&&\sum_{Q'\subset Q} \frac{\mu(Q')^{\frac{1}{p-1}}}{\m{Q'}^{(1-\frac{\alpha p}
{n})\frac{1}{p-1}-1}}\\
&=&\sum_{Q'\subset Q}\Big(\mu(Q')^\frac{1}{p-1}\m{Q'}^\epsilon\Big)\m{Q'}^
{-(1-\frac{\alpha p}{n})\frac{1}{p-1}+1-\ep}\\
&\leq& \Big(\sum_{Q'\subset Q}\mu(Q')^\frac{r'}{p-1}\m{Q'}^{\ep r'}\Big)^{\frac
{1}{r'}}\Big(\sum_{Q'\subset Q}\m{Q'}^{-r(1-\frac{\alpha p}{n})\frac{1}{p-1}+r-r\ep}\Big)^
\frac{1}{r},
\end{eqnarray*}
where $r'=p-1>1$, $r=\frac{p-1}{p-2}$ and $\ep>0$ is chosen so 
that $-r(1-\frac{\alpha p}{n})\frac{1}{p-1}+r-r\ep>1$, i.e., $0<\ep<\frac
{\alpha p}{(p-1)n}$. Therefore,
\begin{eqnarray*}
\sum_{Q'\subset Q} \frac{\mu(Q')^{\frac{1}{p-1}}}{\m{Q'}^{(1-\frac{\alpha p}{n})
\frac{1}{p-1}-1}}&\leq& C\mu(Q)^{\frac{1}{p-1}}\m{Q}^{\ep}\m{Q}^{-(1-\frac{
\alpha p}{n})\frac{1}{p-1}+1-\ep}\\
&=& C\frac{\mu(Q)^\frac{1}{p-1}}{\m{Q}^{(1-\frac{\alpha p}{n})\frac{1}{p-1}-1}}.
\end{eqnarray*}
Hence, combining this with (\ref{A_{2}}) we obtain
\begin{eqnarray*}
A_{2}(P,\mu)&\leq& C\sum_{Q\subset P}\frac{\mu(Q)^\frac{1}{p-1}}{\m{Q}^{(1-
\frac{\alpha p}{n})
\frac{1}{p-1}+q-2}} \Big[\frac{\mu(Q)^{\frac{1}{p-1}}}{\m{Q}^{(1-\frac{
\alpha p}{n})\frac{1}{p-1}-1}}\Big]^{q-1}\\
&=&C\sum_{Q\subset P} \frac{\mu(Q)^{\frac{q}{p-1}}}{\m{Q}^{(1-\frac{\alpha p}
{n})\frac{q}{p-1}-1}}=CA_{1}(P,\mu).
\end{eqnarray*} 
This completes the proof of the proposition. 
\end{proof}
\begin{theorem}\label{maindiscrete}
Let $\alpha>0$, $p>1$ be such that $0<\alpha p<n$, and let $q>p-1$. Suppose
$f\in L^{q}_{\rm loc}(\RR^n)$ and $d\om=f^q dx$. Then the following statements are 
equivalent.\\  
\noindent{\rm (i)} The equation 
\begin{equation}\label{equ}
u=\mathcal{W}_{\alpha,\, p}(u^q) + \epsilon f
\end{equation}
has a solution $u\in L_{\rm loc}^{q}(\RR^n)$ for some $\epsilon>0$.\\
\noindent{\rm (ii)} The testing inequality 
\begin{equation}\label{testing1}
\int_{P}\Big[\sum_{Q\subset P}\frac{\om(Q)}{\m{Q}^{1-\frac{\alpha p}{n}}
}\chi_{Q}(x)\Big]^{\frac{q}{p-1}}dx\leq C\m{P}_{\om}
\end{equation}
holds for all dyadic cubes P. \\
%\noindent{\rm (iii)} The weak-type inequality
%\begin{equation}\label{weaktype1}
%\norm{\mathcal{I}_{\alpha p}(g)}_{L^{\frac{q}{q-p+1},\,\infty}(d\mu)}\leq C ||
%g||_{L^{\frac{q}{q-p+1}}(dx)}
%\end{equation} 
%holds for all $0\leq g\in L^{\frac{q}{q-p+1}}(\RR^n)$.\\
\noindent{\rm(iii)}  The testing inequality
\begin{equation}\label{testing2}
\int_{P}\Big[\sum_{Q\subset P} \frac{\om(Q)^{\frac{1}{p-1}}}{\m{Q}^
{(1-\frac{\alpha p}{n})\frac{1}{p-1}}}\chi_{Q}(x)\Big]^{q}dx\leq C \m{P}_{\om}
\end{equation}
holds for all dyadic cubes P. \\
\noindent{\rm(iv)} There exists a constant $C$ such that 
\begin{equation}\label{pw}
\mathcal{W}_{\alpha,\,p}[\mathcal{W}_{\alpha,\,p}(f^q)]^q (x) \leq C 
\mathcal{W}_{\alpha,\,p}(f^q)(x)<\infty \quad {\rm a.e..}
\end{equation}
\end{theorem}
\begin{proof}
We show that (iv)$\Longrightarrow$ (i)$\Longrightarrow$ (ii) $\Longrightarrow$ (iii)
$\Longrightarrow$ (iv). Note that by Proposition \ref{discrete0t} 
we  have (ii)$\Longleftrightarrow$ (iii).
% Moreover, it is known that the testing 
%inequality (\ref{testing1}) is also equivalent to the weak-type inequality 
%(\ref{weaktype1}) (see e.g. \cite{NTV}, \cite{VW}).  
Therefore, it is enough to prove that (iv)$\Longrightarrow$ 
(i)$\Longrightarrow$ (iii)$\Longrightarrow$ (iv).\\

\noindent{\it Proof of {\rm(iv)}$\Longrightarrow$ {\rm(i)}.} Note that 
the pointwise condition (\ref{pw}) can be 
rewritten as 
$$\mathcal{N}^2f\leq C \mathcal{N}f<\infty \quad {\rm a.e.,}$$
where $\mathcal{N}$ is the operator defined by (\ref{operatorA}). 
The sufficiency of this condition for the solvability of 
(\ref{equ}) can be proved using simple iterations:
$$u_{n+1}=\mathcal{N}u_{n} + \epsilon f, \hspace*{.3in} n=0,~1,~2,...$$
starting from $u_{0}=0$. Since $\mathcal{N}$ is monotonic it is easy to see that 
$u_{n}$ is increasing and that $\epsilon^{\frac{q}{p-1}}\mathcal{N}f + \epsilon
f\leq u_{n}$ for all $n\geq 2$. 
Let $c(p)=\max\{1,2^{p'-1}\}$, $c_{1}=0$, $c_{2}=[\epsilon^{\frac{1}{p-1}} 
c(p)]^q$ and 
$$c_{n}=\Big[\epsilon^{\frac{1}{p-1}}c(p)(1+C^{1/q})c_{n-1}^{p'-1}
\Big]^q,\hspace*{.3in} n=3,~4,...$$
where $C$ is the constant in $(\ref{pw})$. Here we choose $\epsilon$ so that
$$\epsilon^{\frac{1}{p-1}}c(p)=\Big(\frac{q-p+1}{q}\Big)^{\frac{q-p+1}{q}}
\Big(\frac{p-1}{q}\Big)^{\frac{p-1}{q}}C^{\frac{1-p}{q^2}}.$$
By induction and using (\ref{propertyA}) we have
$$u_{n}\leq c_{n}\mathcal{N}f + \epsilon f,\hspace*{.3in}n=1,~2,~3,....$$  
Note that 
$$x_{0}=\Big[\frac{q}{p-1}\epsilon^{\frac{1}{p-1}}c(p)C^{1/q}\Big]^{\frac{q(p-1
)}{p-1-q}}$$
is the only root of the equation
$$x=\Big[\epsilon^{\frac{1}{p-1}}c(p)(1+C^{1/q}x)\Big]^{q}$$
and thus $\lim_{n\rightarrow\infty} c_{n}=x_{0}$.
Hence there exists a solution 
$$u(x)=\lim_{n\rightarrow \infty}u_{n}(x)$$ to 
equation (\ref{equ}) (with that choice of $\epsilon$) such that
$$\epsilon f+\epsilon^{\frac{q}{p-1}}\mathcal{W}_{\alpha,\,p}(f^q)\leq u(x)\leq 
\epsilon f + x_{0}\mathcal{W}_{\alpha,\,p}(f^q).$$ 

\noindent{\it Proof of {\rm(i)}$\Longrightarrow$ {\rm(iii)}.} Suppose that $u\in L^{q}_{\rm loc}(\RR^n)$
is a solution of (\ref{equ}). Let $P$ be a cube in $\mathcal{D}$ and $d\mu=u^q 
dx$. Since 
$$[u(x)]^q\geq [\mathcal{W}_{\alpha,\,p}(u^q)(x)]^q \quad {\rm a.e.,} $$ 
we have 
$$\int_{P}[\mathcal{W}_{\alpha,\,p}(u^q)(x)]^qdx\leq \int_{P}
[u(x)]^qdx.$$
Thus,
\begin{equation}\label{formu}
\int_{P}\Big[\sum_{Q\subset P} \frac{\mu(Q)^{\frac{1}{p-1}}}{\m{Q}^
{(1-\frac{\alpha p}{n})\frac{1}{p-1}}}\chi_{Q}(x)\Big]^{q}dx\leq C \m{P}_{\mu}, 
\end{equation}
for all $P\in \mathcal{D}.$ By Proposition \ref{discrete0t}, inequality 
(\ref{formu}) is equivalent to
 \begin{equation*}
\int_{P}\Big[\sum_{Q\subset P}\frac{\mu(Q)}{\m{Q}^{1-\frac{\alpha p}{n}}
}\chi_{Q}(x)\Big]^{\frac{q}{p-1}}dx\leq C\m{P}_{\mu}
\end{equation*}
for all $P\in \mathcal{D}$, which in its turn is equivalent to the weak-type inequality 
\begin{equation}\label{weaktype2}
\norm{\mathcal{I}_{\alpha p}(g)}_{L^{\frac{q}{q-p+1},\, \infty}(d\mu)}\leq C \norm{
g}_{L^{\frac{q}{q-p+1}}(dx)},
\end{equation} 
for all $g\in L^{\frac{q}{q-p+1}}(\RR^n)$, $g \ge 0$ (see \cite{NTV}, \cite{VW}). 
Note that by (\ref{equ}),
$$d\mu=u^q dx\geq \epsilon^q f^q \, dx=\epsilon^q \, d\om.$$
We now deduce from (\ref{weaktype2}),  
\begin{equation}\label{forom}
\norm{\mathcal{I}_{\alpha p}(g)}_{L^{\frac{q}{q-p+1},\,\infty}(d\om)}\leq \frac{C}{
\epsilon^{q-p+1}} \norm{
g}_{L^{\frac{q}{q-p+1}}(dx)}
\end{equation} 
Similarly, by duality and Proposition \ref{discrete0t} we see that 
(\ref{forom}) is equivalent to the 
testing inequality (\ref{testing2}). The implication (i)$\Longrightarrow$
(iii) is proved.\\  

\noindent{\it Proof of {\rm(iii)}$\Longrightarrow$ {\rm(iv)}.} We first deduce from the testing
inequality (\ref{testing2}) that 
\begin{equation}\label{frostman}
\m{P}_{\om}\leq C \m{P}^{1-\frac{\alpha p q}{n(q-p+1)}}.
\end{equation}
for all dyadic cubes $P$. In fact, this can be verified by using (\ref{testing2})
and the obvious estimate
$$\int_{P}\Big[\frac{\om(P)}{\m{P}^{1-\alpha p/n}}\Big]^{\frac{q}{p-1}}dx\leq
\int_{P}\Big[\sum_{Q\subset P} \frac{\om(Q)^{\frac{1}{p-1}}}{\m{Q}^
{(1-\frac{\alpha p}{n})\frac{1}{p-1}}}\chi_{Q}(x)\Big]^{q}dx.$$ 

Following \cite{KV}, \cite{V3}, we next introduce a certain decomposition of the dyadic Wolff potential 
$\mathcal{W}_{\alpha,\,p} \mu$.
To each dyadic cube $P\in
\mathcal{D}$, we associate the  ``upper" and 
``lower" parts of $\mathcal{W}_{\alpha,\,p} \mu$ 
 defined respectively by
$$\mathcal{U}_{P}\mu(x)=\sum_{Q\subset P}\Big[ \frac{\mu(Q)}{\m{Q}^{1-\alpha p/n}}
\Big]^{\frac{1}{p-1}}\chi_{Q}(x),$$  
$$\mathcal{V}_{P}\mu(x)=\sum_{Q\supset P}\Big[ \frac{\mu(Q)}{\m{Q}^{1-\alpha p/n}}
\Big]^{\frac{1}{p-1}}\chi_{Q}(x).$$  
Obviously,
  $$\mathcal{U}_{P}\mu(x)\leq \mathcal{W}_{\alpha,\,p}\mu(x),\hspace*{.3in} 
\mathcal{V}_{P}\mu(x)\leq \mathcal{W}_{\alpha,\,p}\mu(x),$$
 and for $x\in P$,
\begin{equation*}  
\mathcal{W}_{\alpha,\,p}\mu(x)= \mathcal{U}_{P}\mu(x)+\mathcal{V}_{P}\mu(x)
-\Big[\frac{\mu(P)}{\m{P}^{1-\alpha p/n}}\Big]^{\frac{1}{p-1}}.
\end{equation*}
Using the notation just introduced, we can rewrite the testing inequality 
(\ref{testing2}) in the form:
\begin{equation}\label{testing2'}
\int_{P}[\mathcal{U}_{P}\om(x)]^q \, dx\leq C \, \m{P}_{\om},     
\end{equation}
for all dyadic cubes $P$. Recall that $d\om=f^q \, dx$. The desired pointwise inequality 
(\ref{pw}) 
can be restated as 
\begin{equation}\label{pwexplicit}
\sum_{P\in\mathcal{D}}\Big[\frac{\int_{P}[\mathcal{W}_{\alpha,\,p}\om(y)]^q \, dy}
{\m{P}^{1-\alpha p/n}}\Big]^{\frac{1}{p-1}}\chi_{P}(x)\leq C \, \mathcal{W}_{\alpha,\,
p}\om(x).
\end{equation}
From the discussion above we have, for $y\in P$,
$$\mathcal{W}_{\alpha,\,p}\om(y)\leq \mathcal{U}_{P}\om(y)+\mathcal{V}_{P}\om(y)
$$ 
while from the testing inequality (\ref{testing2'}),
\begin{equation*}
\sum_{P\in\mathcal{D}}\Big[\frac{\int_{P}[\mathcal{U}_{P}\om(y)]^q \, dy}
{\m{P}^{1-\alpha p/n}}\Big]^{\frac{1}{p-1}}\chi_{P}(x)\leq C \, \mathcal{W}_{\alpha,\,
p}\om(x).
\end{equation*}
Therefore, to prove (\ref{pwexplicit}) it enough to prove  
\begin{equation}\label{pwlower}
\sum_{P\in\mathcal{D}}\Big[\frac{\int_{P}[\mathcal{V}_{P}\om(y)]^q \, dy}
{\m{P}^{1-\alpha p/n}}\Big]^{\frac{1}{p-1}}\chi_{P}(x)\leq C \, \mathcal{W}_{\alpha,\,
p}^{d} \om(x).
\end{equation}
Note that, for $y\in P$, 
$$\mathcal{V}_{P}\om(y)=\sum_{Q\supset P}\Big[ \frac{\om(Q)}{\m{Q}
^{1-\alpha p/n}}\Big]^{\frac{1}{p-1}}={\rm const}.$$
 An application of the elementary inequality
$$\Big(\sum_{k=1}^{\infty} a_{k}\Big)^s\leq s\sum_{k=1}^{\infty}a_{k}
\Big(\sum_{j=k}^{\infty} a_{j}\Big)^{s-1}  $$
where $1\leq s<\infty$ and $0\leq a_{k}<\infty$, then gives 
\begin{eqnarray*}
[\mathcal{V}_{P}\om(y)]^{\frac{q}{p-1}}&\leq& C\sum_{Q\supset P}\Big[ \frac{\om(Q)}{\m{Q}
^{1-\alpha p/n}}\Big]^{\frac{1}{p-1}}\Big\{
\sum_{R\supset Q}\Big[ \frac{\om(R)}{\m{R}
^{1-\alpha p/n}}\Big]^{\frac{1}{p-1}} \Big\}^{\frac{q}{p-1}-1}.
\end{eqnarray*}
Using this inequality we see that the left-hand side of (\ref{pwlower}) 
is bounded from above by a constant multiple  of 
$$\sum_{P\in \mathcal{D}}\m{P}^\frac{\alpha p}{n(p-1)}
\sum_{Q\supset P}\Big[ \frac{\om(Q)}{\m{Q}^{1-\alpha p/n}}\Big]^{\frac{1}
{p-1}}\Big\{\sum_{R\supset Q}\Big[ \frac{\om(R)}{\m{R}^{1-\alpha p/n}}
\Big]^{\frac{1}{p-1}} \Big\}^{\frac{q}{p-1}-1}\chi_{P}(x).$$
Changing the order of summation, we see that it is equal to 
$$\sum_{Q\in \mathcal{D}}\Big[ \frac{\om(Q)}{\m{Q}^{1-\alpha p/n}}\Big]^{\frac{1}
{p-1}}\chi_{Q}(x) \Big\{\sum_{P\subset Q}\m{P}^\frac{\alpha p}{n(p-1)}\chi_{P}(x)
[\mathcal{V}_{Q}\om(x)]^{\frac{q}{p-1}-1}\Big\}.$$
By (\ref{frostman}), the expression in the curly brackets above is uniformly 
bounded. Therefore, the proof of estimate (\ref{pwlower}), and hence of  
(iii)$\Longrightarrow$ (iv), is complete. 
\end{proof}

%**************************************************************
%**************************************************************
\section{$\mathcal{A}$-superharmonic functions}\label{pre}
In this section, we recall for later use some facts on $\mathcal{A}$-superharmonic functions, 
most of which can be found in \cite{HKM}, \cite{KM1}, \cite{KM2}, and   \cite{TW4}. 
Let $\Om$ be an open set in $\RR^n$, and $p>1$. 
We will mainly be interested in the case where $\Om$ is bounded and $1<p\leq n$, or 
$\Om=\RR^n$ and $1<p<n$. 
We assume  that $\mathcal{A}:\RR^n\x \RR^n\ra\RR^n$ is a 
vector valued mapping which satisfies the following structural properties: 
\begin{eqnarray}\label{2.1}
&&\text{the~ mapping~ x}\ra \mathcal{A}(x,\xi) {\rm ~is~ measurable~ for~ all~}\xi\in\RR^n,\\
&&\text{the~ mapping ~}\xi\ra \mathcal{A}(x,\xi){\rm ~is~ continuous~ for~ a.e.~ x~}\in\RR^n, 
\label{2.2}
\end{eqnarray}
and there are constants $0<\alpha\leq\beta<\infty$ such that for a.e. $x$ in $\RR^n$, 
and for all $\xi$ in $\RR^n$, 
\begin{eqnarray}
\label{2.3}
&\mathcal{A}(x,\xi)\cdot\xi\geq \alpha\m{\xi}^p,\quad \m{\mathcal{A}(x,\xi)}\leq\beta\m{\xi}^{p-1},&\\
\label{2.5}
&(\mathcal{A}(x,\xi_{1})-\mathcal{A}(x,\xi_{2}))\cdot(\xi_{1}-\xi_{2})>0, \quad \text{~if~}
\xi_{1}\not = \xi_{2},&\\
\label{2.6}
&\mathcal{A}(x,\lambda\xi)=\lambda\m{\lambda}^{p-2}\mathcal{A}(x,\xi), \quad {\rm ~if~}\lambda
\in\RR\setminus\{0\}.&
\end{eqnarray}
\indent For $u\in W^{1,\,p}_{{\rm loc}}(\Om)$, we define the divergence of 
$\mathcal{A}(x,\nabla u)$ in the sense of distributions, i.e., if $\varphi\in 
C^{\infty}_{0}(\Om)$, then 
\begin{eqnarray*}
{\rm div}\mathcal{A}(x, \nabla u)(\varphi)=-\int_{\Om}\mathcal{A}(x, \nabla u)
\cdot\nabla \varphi \, dx.
\end{eqnarray*}
\indent It is well known that every solution $u\in W^{1,\,p}_{{\rm loc}}(\Om)$ to the 
equation  
\begin{eqnarray}
\label{homo}
-\text{div}\mathcal{A}(x,\nabla u)=0
\end{eqnarray}
has a continuous representative. Such continuous solutions are said to be 
$\mathcal{A}$-$harmonic$ in $\Om$. If $u\in W_{{\rm loc}}^{1,\,p}(\Om)$ and
\begin{eqnarray*}
\int_{\Om}\mathcal{A}(x,\nabla u)\cdot\nabla\varphi \, dx\geq 0, 
\end{eqnarray*}
for all nonnegative $\varphi\in C^{\infty}_{0}(\Om)$, i.e., $-{\rm div}\mathcal{A}(x,\nabla u)
\geq 0$ in the distributional sense, then $u$ is called a $supersolution$
of the equation (\ref{homo}) in $\Om$.\\
\indent A lower semicontinuous function $u:\Om\ra (-\infty, \infty]$ is called 
$\mathcal{A}$-$superharmonic$ if $u$ is not identically infinite in each component
of $\Om$, and if for all open sets  $D$ such that 
${\overline D}\subset\Om$, and all functions $h\in C(\overline{D})$, $\mathcal{A}$-harmonic in $D$, it follows that 
$h\leq u$ on $\partial D$ implies $h\leq u$ in $D$.

It is worth mentioning that $p$-superharmonicity can also be defined equivalently 
 using the language of viscosity solutions (see \cite{JLM}).

\indent  We recall here the fundamental connection between supersolutions of (\ref{homo}) and 
$\mathcal{A}$-superharmonic functions \cite{HKM}. 
\begin{proposition}[\cite{HKM}] \label{pro2.1}{\rm (i)} If $u\in W_{{\rm loc}}^{1,\,p}(\Om)$ 
is such that 
\begin{eqnarray*}
-{\rm div}\mathcal{A}(x,\nabla u)\geq 0,
\end{eqnarray*}
then there is an $\mathcal{A}$-superharmonic function $v$ such that $u=v$ a.e.. Moreover,
\begin{eqnarray}
\label{liminf}
v(x)={\rm ess}\lim_{y\ra x}{\rm inf~} v(y), \qquad  x\in\Om.
\end{eqnarray}
\indent {\rm (ii)} If v is $\mathcal{A}$-superharmonic, then (\ref{liminf}) holds. Moreover,
if $v\in W^{1,\,p}_{{\rm loc}}(\Om)$, then 
\begin{eqnarray*}
-{\rm div}\mathcal{A}(x, \nabla v)\geq 0. 
\end{eqnarray*}
\indent {\rm (iii)} If v is $\mathcal{A}$-superharmonic and locally bounded, then 
$v\in W^{1,\, p}_{{\rm loc}}(\Om)$,    and 
\begin{eqnarray*}
 -{\rm div}\mathcal{A}(x, \nabla v)\geq 0.
\end{eqnarray*}
\end{proposition}
\indent Note that an $\mathcal{A}$-superharmonic function $u$ does not necessarily belong
to $W^{1,\,p}_{{\rm loc}}(\Om)$, but its truncation $\min\{u,k\}$ does, for every integer
$k$, by Proposition \ref{pro2.1}(iii). Using this we set 
\begin{eqnarray*}
Du=\lim_{k\ra\infty} \, \nabla \, [ \, \min\{u,k\}], 
\end{eqnarray*}  
defined a.e. If either $u\in L^{\infty}(\Om)$ or $u\in W^{1,\,1}_{{\rm loc}}(\Om)$, then $Du$ 
coincides with the regular distributional gradient of $u$. In general we have the following
gradient estimates \cite{KM1} (see also \cite{HKM}, \cite{TW4}).
\begin{proposition}[\cite{KM1}]\label{gradient} Suppose u is $\mathcal{A}$-superharmonic in $\Om$ and 
$1\leq q< \frac{n}{n-1}$. Then both $\m{Du}^{p-1}$ and $\mathcal{A}(\cdot,Du)$
belong to $L^{q}_{{\rm loc}}(\Om)$. Moreover, if  $p>2-\frac{1}{n}$, then $Du$ is 
the distributional gradient of u.
\end{proposition}
\indent We can now extend the definition of the divergence of $\mathcal{A}(x, 
\nabla u)$ if $u$ is merely an $\mathcal{A}$-superharmonic function in $\Om$. For such $u$ we set 
\begin{eqnarray*}
-{\rm div}\mathcal{A}(x, \nabla u)(\varphi)=\int_{\Om}\mathcal{A}(x, Du)
\cdot \nabla \varphi \, dx, 
\end{eqnarray*}  
for all $\varphi\in C^{\infty}_{0}(\Om)$. Note that by Proposition \ref{gradient} and the 
dominated convergence theorem, 
\begin{eqnarray*}
-{\rm div}\mathcal{A}(x, \nabla u)(\varphi)=\lim_{k\ra\infty}\int_{\Om}\mathcal
{A}(x, \nabla\min\{u,k\})\cdot \nabla \varphi \, dx\geq 0 
\end{eqnarray*}
whenever $\varphi\in C^{\infty}_{0}(\Om)$ and $\varphi\geq 0$.\\
\indent Since $-{\rm div}\mathcal{A}(x, \nabla u)$ is a nonnegative distribution in $\Om$ for 
an $\mathcal{A}$-superharmonic  $u$, it follows that there is 
a positive (not necessarily finite) Radon measure denoted by $\mu[u]$ such that 
\begin{eqnarray*}
-{\rm div}\mathcal{A}(x, \nabla u)=\mu[u] \quad {\rm in} \quad \Om. 
\end{eqnarray*}
Conversely, given a positive finite measure $\mu$ in a bounded $\Om$, there is an 
$\mathcal{A}$-superharmonic function u such that $-{\rm div}\mathcal{A}(x, \nabla u)=\mu$
in $\Om$ and $\min\{u,k\}\in W^{1,p}_{0}(\Om)$ for all integers $k$. Moreover, if 
$\mu$ is a positive finite measure in $\RR^n$ we can also find a positive 
$\mathcal{A}$-superharmonic function $u$ such that $-{\rm div}\mathcal{A}(x, \nabla u)=\mu$
in $\RR^n$. We refer to \cite{KM1} and \cite{KM2} for details.\\
\indent The following weak continuity result in \cite{TW4} will be used later in Sec. \ref{R^n}
to prove the  existence of $\mathcal{A}$-superharmonic solutions to quasilinear 
equations.    
\begin{theorem}[\cite{TW4}]\label{weakcont} Suppose that $\{u_{n}\}$ is a sequence of nonnegative 
$\mathcal{A}$-superharmonic functions in $\Om$ that converges a.e. to  an 
$\mathcal{A}$-superharmonic function $u$. Then 
the sequence of measures $\{\mu[u_{n}]\}$ converges to $\mu[u]$ weakly, i.e.,
$$\lim_{n\rightarrow\infty}\int_{\Om}\varphi \, d\mu[u_{n}]=\int_{\Om}\varphi \, d\mu[u],$$
for all $\varphi\in C_{0}^{\infty}(\Om)$.
\end{theorem}
\indent In \cite{KM2} (see also \cite[Theorem 3.1]{Mi} and \cite{MZ}) the following 
pointwise potential  estimate for $\mathcal{A}$-superharmonic
functions was established, which serves as a major tool in our study of quasilinear equations of Lane--Emden type.   
\begin{theorem}[\cite{KM2}]\label{potential}
Suppose  $u \ge 0$ is an $\mathcal{A}$-superharmonic function in $B(x,3r)$.
If $\mu=-{\rm div}\mathcal{A}(x, \nabla u)$, then
\begin{eqnarray*}
C_{1} \, {\rm\bf W}_{1,\,p}^{r}\mu(x)\leq u(x)\leq C_{2} \, \inf_{B(x,r)}u + C_{3} \, {\rm\bf W}_{1,\,p}^{2r}\mu(x),
\end{eqnarray*} 
where $C_{1}, C_{2}$ and $C_{3}$ are positive constants which depend only on $n, p$ and 
the structural constants $\alpha$ and $\beta$.
\end{theorem}
\indent A consequence of Theorem \ref{potential} is the following global version of the above potential
pointwise estimate.
\begin{corollary}[\cite{KM2}]\label{globalcor}
Let $u$ be an $\mathcal{A}$-superharmonic function in $\RR^n$ with $\inf_{\RR^n}u=0$. If
$\mu=-{\rm div}\mathcal{A}(x, \nabla u)$, then 
\begin{eqnarray}
\label{global}
\frac{1}{K} \, {\rm \bf W}_{1,\,p}\mu(x)\leq u(x)\leq K \, {\rm \bf W}_{1,\,p}\mu(x), 
\end{eqnarray}
for all $x\in \RR^n$, where $K$ is a positive constant depending only on $n, p$ and the structural constants
$\alpha$ and $\beta$. 
\end{corollary}

%**********************************************************
%**********************************************************
\section{Quasilinear equations on $\RR^n$} \label{R^n}
In this section, we study the solvability problem for the quasilinear equation
\begin{eqnarray}\label{globalA}
-{\rm div}\mathcal{A}(x, \nabla u)=u^{q}+\om
\end{eqnarray}
in the class of nonnegative $\mathcal{A}$-superharmonic functions 
on the entire space $\RR^n$, where $\mathcal{A}(x,\xi)\cdot\xi\approx\m{\xi}^{p}$ is defined
precisely as in Sec.~4. Here we assume $1<p<n$, $q>p-1$, and $\om\in {\mathcal M}^{+}(\RR^n)$. 
In this setting, all solutions are understood in the ``potential-theoretic" sense, i.e., 
$0\leq u\in L^{q}_{{\rm loc}}(\RR^n)$
is a solution to (\ref{globalA}) if $u$ is an $\mathcal{A}$-superharmonic function and
\begin{eqnarray}
\int\lim_{k\ra\infty}\mathcal{A}(x,\nabla \min\{u, k\})\cdot\nabla\varphi \, dx=\int u^{q}\varphi \, dx +
\int\varphi \, d\om, 
\end{eqnarray} 
for all test functions $\varphi\in C^{\infty}_{0}(\RR^n)$.\\
\indent First we prove the continuous counterpart of Proposition \ref{discrete0t}. Here 
we use the well-known argument due to Fefferman and Stein \cite{FS} which is based on 
the averaging over shifts of the dyadic lattice $\mathcal{D}$.

\begin{proposition}\label{cont}
Let $0<r\leq\infty$. Let $\mu\in {\mathcal M}^{+}(\RR^n)$, $\alpha>0$, $p>1$, and $q>p-1$. 
Then the following quantities are equivalent.
\begin{eqnarray*}
&(a)&\norm{{\rm\bf W}_{\alpha p,\,\frac{q}{q-p+1}}^{r}\mu}_{L^{1}(d\mu)}=
\int_{\RR^{n}}\int_{0}^{r}\Big[\frac{\mu(B_{t}(x))}{t^{n-\frac{\alpha pq}
{q-p+1}}}\Big]^{\frac{q}{p-1}-1}\frac{dt}{t}d\mu,\\
&(b)& \norm{{\rm\bf W}_{\alpha,\,p}^{r}\mu}_{L^{q}(dx)}^{q}=\int_{\RR^{n}}\Big
\{\int_{0}^{r}\Big[\frac{\mu(B_{t}(x))}{t^{n-\alpha p}}\Big]
^{\frac{1}{p-1}}\frac{dt}{t}\Big\}^{q}dx,\\
&(c)& \norm{{\rm\bf I}_{\alpha p}^{r}\mu}_{L^{\frac{q}{p-1}}(dx)}^{\frac{q}{p-1}}=
\int_{\RR^{n}}\Big[\int_{0}^{r}\frac{\mu(B_{t}(x))}{t^{n-\alpha p}}\frac{dt}{t}
\Big]^{\frac{q}{p-1}}dx,
\end{eqnarray*}
where the constants of equivalence do not depend on $\mu$ and $r$. 
\end{proposition}
\begin{proof}
We will prove only the equivalence  of $(b)$ and $(c)$, i.e., there are constants 
$C_{1},~C_{2} > 0$ such that  
\begin{equation}\label{seconddir}
C_{1}\norm{{\rm\bf W}_{\alpha,\,p}^{r}\mu}_{L^{q}(dx)}^{q}\leq 
\norm{{\rm\bf I}_{\alpha p}^{r}\mu}_{L^{\frac{q}{p-1}}(dx)}^{\frac{q}{p-1}}
\leq C_{2} \norm{{\rm\bf W}_{\alpha,\,p}^{r}\mu}_{L^{q}(dx)}^{q}.
\end{equation} 
The equivalence of $(a)$ and $(c)$ which follows from Wolff's inequality (see \cite{AH}, \cite{HW}),
can also be deduced by a similar argument. We first restrict ourselves to the 
case $r<\infty$.
Observe that there is a constant $C>0$ such that 
\begin{equation}\label{rand2r}
\norm{{\rm\bf I}_{\alpha p}^{2r}\mu}_{L^{\frac{q}{p-1}}(dx)}^{\frac{q}{p-1}}
\leq C \norm{{\rm\bf I}_{\alpha p}^{r}\mu}_{L^{\frac{q}{p-1}}(dx)}^{\frac{q}{p-1}}.
\end{equation}
In fact, since  
$$\int_{0}^{2r}\frac{\mu(B_{t}(x))}{t^{n-\alpha p}}\frac{dt}{t} 
\leq C \int_{0}^{r}\frac{\mu(B_{t}(x))}{t^{n-\alpha p}}
\frac{dt}{t} +C \frac{\mu(B_{2r}(x))}{r^{n-\alpha p}},$$
(\ref{rand2r}) will follow from the estimate 
\begin{equation}\label{5.4}
\int_{\RR^n}\Big[\frac{\mu(B_{2r}(x))}{r^{n-\alpha p}}\Big]^{\frac{q}{p-1}} dx
\leq C \int_{\RR^n}\Big[\int_{0}^{r}\frac{\mu(B_{t}(x))}{t^{n-\alpha p}}\frac{dt}
{t}\Big]^{\frac{q}{p-1}} dx.
\end{equation}
Note that for a partition of $\RR^n$ into a union of disjoint cubes $\{Q_{j}\}$ such that ${\rm diam}
(Q_{j})=r/4$ we have
\begin{eqnarray*}
\int_{\RR^n}\mu(B_{2r}(x))^{\frac{q}{p-1}}dx &=& \sum_{j}\int_{Q_{j}} \mu(B_{2r}(x))^
{\frac{q}{p-1}}dx\\
&\leq& C\sum_{j}\int_{Q_{j}}\mu(Q_{j})^{\frac{q}{p-1}}dx,
\end{eqnarray*}
where we  have used the fact that the ball $B_{2r}(x)$ is contained in the union of at most
$N$ cubes in $\{Q_{j}\}$ for some constant $N$ depending only on $n$. 
Thus
\begin{eqnarray*}
\int_{\RR^n}\Big[\frac{\mu(B_{2r}(x))}{r^{n-\alpha p}}\Big]^{\frac{q}{p-1}}dx 
&\leq& C\sum_{j}\int_{Q_{j}}\Big[\frac{\mu(B_{r/2}(x))}{r^{n-\alpha p}}\Big]^
{\frac{q}{p-1}}dx\\
&\leq& C\sum_{j}\int_{Q_{j}}\Big[\int_{0}^{r}\frac{\mu(B_{t}(x))}{t^{n-\alpha p}}
\frac{dt}{t}\Big]^
{\frac{q}{p-1}}dx,
\end{eqnarray*}
which gives  (\ref{5.4}).   \\
\indent By arguing as in \cite{COV}, we can find constants $a$, $C$ and $c$ depending only on $n$ such that 
\begin{equation*}
{\rm\bf W}_{\alpha,\,p}^{r}\mu(x)\leq C r^{-n}\int_{\m{t}\leq cr}\sum_{\substack
{Q\in\mathcal{D}_{t}\\\ell(Q)\leq 4r/a}}
\Big[\frac{\mu(Q)}{\m{Q}^{1-\alpha p/n}} \Big]^{\frac{1}{p-1}}\chi_{Q}(x)dt
\end{equation*}   
where $\mathcal{D}_{t}$, $t\in\RR^n$, denotes the lattice $\mathcal{D}+t=\{Q=Q'+t: Q'\in
\mathcal{D}\}$
and $\ell(Q)$ is the side length of $Q$.  Using  Proposition 2.2 in \cite{COV} 
and arguing as in the proof of Theorem \ref{discrete0t} we obtain
\begin{eqnarray*}
&&\int_{\RR^n}\Big\{\sum_{\substack{Q\in\mathcal{D}_{t}\\\ell(Q)\leq 4r/a}}
\Big[\frac{\mu(Q)}{\m{Q}^{1-\alpha p/n}} \Big]^{\frac{1}{p-1}}\chi_{Q}(x)\Big\}^{q}dx\\
&\asymp& \int_{\RR^n}\Big[\sum_{\substack{Q\in\mathcal{D}_{t}\\\ell(Q)\leq 4r/a}}
\frac{\mu(Q)}{\m{Q}^{1-\alpha p/n}} \chi_{Q}(x)\Big]
^{\frac{q}{p-1}}dx
\end{eqnarray*}
where the constant of equivalence is independent of $\mu$, $r$ and $t$. 
The last two estimates together with the integral   Minkowski  
inequality then give
\begin{eqnarray*}
&&||{\rm\bf W}_{\alpha,\,p}^{r}\mu||_{L^{q}(dx)}\\
&\leq&C r^{-n}\int_{\m{t}\leq cr}\Big\{
\int_{\RR^n}\Big(\sum_{\substack{Q\in\mathcal{D}_{t}\\\ell(Q)\leq 4r/a}}
\Big[\frac{\mu(Q)}{\m{Q}^{1-\alpha p/n}} \Big]^{\frac{1}{p-1}}\chi_{Q}(x)\Big)
^{q}dx\Big\}^{\frac{1}{q}}dt\\
&\leq& C r^{-n}\int_{\m{t}\leq cr}\Big[
\int_{\RR^n}\Big(\sum_{\substack{Q\in\mathcal{D}_{t}\\\ell(Q)\leq 4r/a}}
\frac{\mu(Q)}{\m{Q}^{1-\alpha p/n}} \chi_{Q}(x)\Big)
^{\frac{q}{p-1}}dx\Big]^{\frac{1}{q}}dt.
\end{eqnarray*}
Note that 
\begin{eqnarray*}
\sum_{\substack{Q\in\mathcal{D}_{t}\\\ell(Q)\leq 4r/a}}\frac{\mu(Q)}
{\m{Q}^{1-\alpha p/n}}\chi_{Q}(x)
&\leq& C \sum_{2^k\leq 4r/a}\frac{\mu(B(x,\sqrt{n}2^{k}))}{2^{k(n-\alpha p)}}\\
&\leq& C {\rm\bf  I}_{\alpha p}^{8r\sqrt{n}/a }\mu(x).
\end{eqnarray*}
where $C$ is independent of $t$. Thus, in view of (\ref{rand2r}), 
we obtain the lower estimate  in (\ref{seconddir}).\\
\indent Now by letting $R\rightarrow\infty$ in the inequality
\begin{eqnarray*}
||{\rm\bf W}^{R}_{\alpha,\,p}\mu||_{L^{q}(dx)}^{q}
\leq C||{\rm\bf I}_{\alpha p}^{R}\mu||_{L^{\frac{q}{p-1}}(dx)}^
{\frac{q}{p-1}},
\hspace{.3in} 0<R<\infty,
\end{eqnarray*}
we get the lower estimate in (\ref{seconddir}) with $r=\infty$. 
The upper estimate in (\ref{seconddir}) can be deduced in a similar way. 
This completes the proof of Proposition \ref{cont}.
\end{proof}

The next theorem gives a characterization of the existence of nonnegative solutions
to the equation $-{\rm div}\mathcal{A}(x, \nabla u) =\mu$ 
in terms of Wolff's potentials. 
\begin{theorem}\label{deltapu=mu}
Let $\mu$ be a measure in ${\mathcal M}^{+}(\RR^{n})$. Suppose that ${\rm\bf W}_{1,\,p}\mu<\infty$ a.e. 
 Then there is a nonnegative $\mathcal{A}$-superharmonic
function u in $\RR^{n}$ such that  
\begin{eqnarray}\label{linear}
-{\rm div}\mathcal{A}(x, \nabla u) =\mu \quad{\rm in}\quad \RR^n,
\end{eqnarray} 
and 
\begin{eqnarray}\label{AE}
\frac{1}{K}{\rm\bf W}_{1,\, p}\mu(x)\leq u(x)\leq K\,{\rm\bf W}_{1,\,p}\mu(x), 
\end{eqnarray}
for all $x$ in $\RR^{n}$, where $K$ is the constant in (\ref{global}).
Conversely, if $u$ is a nonnegative $\mathcal{A}$-superharmonic
function in $\RR^{n}$ which solves (\ref{linear}) then 
${\rm\bf W}_{1,\,p}\mu<\infty$ a.e. on $\RR^{n}$.  
\end{theorem}  
\begin{proof}
The second statement of the theorem follows immediately from the lower Wolff potential estimate. To prove 
the  first statement, we  let $\mu_{k}=\mu_{B_{k}}$, the restriction of $\mu$ on the ball
$B_{k}$ of radius $k$ and centered at the origin, so that $\mu_{k}\rightarrow\mu$ weakly as measures
and $\mu_{k}$ are finite positive Borel measures on $\RR^{n}$. Thus arguing as in the proof
of Theorem 2.4 in \cite{KM1} we can find 
nonnegative $\mathcal{A}$-superharmonic functions $u_{k}$ in $\RR^{n}$ such that 
$$-{\rm div}\mathcal{A}(x, \nabla u_{k})=\mu_{k}$$
in $\RR^{n}$. By replacing $u_{k}$ with $u_{k}-\inf_{\RR^{n}} u_{k}$ we can assume that 
$\inf_{\RR^{n}}u_{k}=0$ so that $u_{k}(x)\leq K\,{\rm\bf W}_{1,\,p}\mu_{k}(x)\leq 
K\,{\rm\bf W}_{1,\,p}\mu(x)<\infty$ for a.e. $x$ in $\RR^{n}$ 
by Corollary \ref{globalcor}. 
Let $\{u_{k_{j}}\}$ be a subsequence of $\{u_{k}\}$ such that $u_{k_{j}}\rightarrow u$ a.e. on $\RR^{n}$
for some nonnegative $\mathcal{A}$-superharmonic function u in $\RR^{n}$ (see \cite[Theorem 1.17]{KM1}).
Then by Theorem \ref{weakcont} we see that 
$$-{\rm div}\mathcal{A}(x, \nabla u)=\mu$$
in $\RR^{n}$. The estimate (\ref{AE}) then follows from the Wolff potential 
estimate, which  completes the proof of the theorem.
\end{proof}
\begin{theorem}\label{sufficency}
Let $\om \in \mathcal{M}^{+}(\RR^{n})$, $1<p<n$, and $q>p-1$. Assume that 
\begin{equation}
\label{W(Wom)^{q}}
{\rm\bf W}_{1,\,p}({\rm\bf W}_{1,\,p}\om)^{q}\leq C\,{\rm\bf W}_{1,\,p}\om<\infty \quad {\rm a.e.,} 
\end{equation}
 where 
\begin{eqnarray}\label{condonC}
C \leq \Big(\frac{q-p+1}{qK\max\{1,2^{p'-2}\}}\Big)^{q(p'-1)}\Big(\frac{p-1}{q-p+1}\Big), 
\end{eqnarray}
and $K$ is the constant in (\ref{global}).
 Then there is an $\mathcal{A}$-superharmonic function $u\in L^{q}_{\rm loc}(\RR^{n})$ such that 
\begin{eqnarray}
\label{suff}
\left\{\begin{array}{c}
\inf_{x\in \RR^{n}}u(x)=0\\
-{\rm div}\mathcal{A}(x, \nabla u)=u^{q}+\om,
\end{array}
\right.
\end{eqnarray}
and
$$\frac{1}{M}{\rm\bf W}_{1,\,p}\om(x)\leq u(x)\leq M \,{\rm\bf W}_{1,\,p}\om(x),$$
for all $x$ in $\RR^{n}$, where the constant $M$ depends  only $n, p, q$ and the 
structural constants $\alpha$ and $\beta$. 
\end{theorem}
\begin{proof}
Let $\{u_{k}\}_{k\geq 0}$ be a sequence of $\mathcal{A}$-superharmonic functions such that 
$\inf_{\RR^{n}}u_{k}=0$, $u_{k}\in L^{q}_{\rm loc}(\RR^n)$,
\begin{eqnarray*}
\int \mathcal{A}(x, \nabla u_{0})\cdot\nabla\varphi dx&=&\int\varphi d\om,
\end{eqnarray*}
and
\begin{eqnarray}\label{m{Du_{k+1}}^{p-2}}
\int \mathcal{A}(x, \nabla u_{k+1})\cdot\nabla\varphi dx&=&\int u_{k}^{q}\varphi \, dx+\int\varphi \, d\om,
\end{eqnarray} 
for all integer $k\geq 0$ and $\varphi\in C_{0}^{\infty}{(\RR^{n})}$.
 The existence of such a sequence is guaranteed by Theorem \ref{deltapu=mu} and 
condition (\ref{W(Wom)^{q}}). Put $c_{0}=K$, where $K$ is the constant in (\ref{global}). 
By the potential estimate we see that $u_{0}\leq c_{0}{\rm \bf W}_{1,\,p}\om$ and
$u_{k+1}\leq K {\bf W}_{1,\, p}(u_{k}^{q}+\om)$ for all $k\geq 0$. From these estimates and (\ref{c(p)}) 
we get
\begin{eqnarray*}
u_{1}&\leq& K \max\{1,2^{p'-2}\}\Big[{\rm\bf W}_{1,\,p}(u_{0}^{q})+{\rm\bf W}_{1,\,p}\mu\Big]\\
&\leq& K \max\{1,2^{p'-2}\}(c_{0}^{q(p'-1)}C+1){\rm\bf W}_{1,\,p}\mu\\
&=& c_{1}{\rm\bf W}_{1,\,p}\mu, 
\end{eqnarray*} 
where $c_{1}=K \max\{1,2^{p'-2}\}(c_{0}^{q(p'-1)}C+1)$. By induction we can find a sequence $\{c_{n}\}_{n\geq 0}$
of positive numbers such that $u_{n}\leq c_{n}{\rm\bf W}_{1,\,p}\mu$ for all $n\geq 0$  with  $c_{0}=K$ and 
$c_{n+1}=K\max\{1,2^{p'-2}\}(c_{n}^{q(p'-1)}C+1)$
for all $n\geq 0$.
It is then easy to see that $c_{n}\leq \frac{K\max\{1,2^{p'-2}\}q}{q-p+1}$ for all $n\geq 0$ as long as
(\ref{condonC}) is satisfied. Thus
$$u_{k}\leq \frac{K\max\{1,2^{p'-2}\}q}{q-p+1}{\rm\bf W}_{1,\,p}\om$$
for all $k\geq 0$. By Theorem 1.17 in \cite{KM1} we can find 
a subsequence which is also denoted by $\{u_{k}\}_{k\geq 0}$ and an $\mathcal{A}$-superharmonic 
function $u$ such 
that $u_{k}\ra u$ almost everywhere. As $k\ra \infty$ in (\ref{m{Du_{k+1}}^{p-2}}), the left-hand side 
tends to $\int\mathcal{A}(x, \nabla u
)\cdot\nabla\varphi dx$ by weak continuity results in \cite{TW4} while the right hand side  tends 
to $\int u^{q}\varphi dx+\int\varphi d\om $ by the  dominated convergence theorem. 
Therefore $u$ is a solution to 
(\ref{suff}),  which completes the proof of the theorem.  
\end{proof}
\begin{theorem}\label{maintheorem1}
Let $\om$ be a locally finite positive  measure on $\RR^n$, $1<p<n$, 
and  $q>p-1$. Then the 
following statements are equivalent.\\
{\rm(i)} There exists a nonnegative solution $u\in 
L^{q}_{\rm loc}(\RR^n)$ to the equation 
\begin{equation}
\label{ep-equation}
\left\{\begin{array}{c}
\inf_{x\in\RR^{n}}u(x)=0\\ 
-{\rm div}\mathcal{A}(x, \nabla u)=u^{q}+\ep\om 
\end{array}
\right.
\end{equation} 
for some $\ep>0$.\\
{\rm(ii)} The testing inequality 
\begin{equation}
\label{testingI_{p}}
\int_{B}\Big[{\rm\bf I}_{p}\om_{B}(x)\Big]^{\frac{q}{p-1}}dx\leq C\om(B)
\end{equation}
holds for all balls $B$ in $\RR^{n}$.\\
{\rm(iii)} The testing inequality 
\begin{equation}
\label{testingW}
\int_{B}\Big[{\rm\bf W}_{1,\,p}\om_{B}(x)\Big]^{q}dx\leq C \om(B)
\end{equation}
holds for all balls $B$ in $\RR^{n}$ .\\
{\rm(iv)} There exists a constant C such that
\begin{equation}
\label{pointwise}
{\rm\bf W}_{1,\,p}({\rm\bf W}_{1,\,p}\om)^{q}\leq C\, {\rm\bf W}_{1,\,p}\om < \infty \quad {\rm a.e..}
\end{equation}
Moreover, 
if the constant $C$ in (\ref{pointwise}) satisfies
\begin{eqnarray*} 
C\leq\Big(\frac{q-p+1}{qK\max\{1,2^{p'-2}\}}\Big)^{q(p'-1)}\Big(\frac{p-1}{q-p+1}\Big),
\end{eqnarray*} 
 where $K$ is the constant in (\ref{global}),
 then the equation (\ref{ep-equation}) has a solution $u$ with $\ep=1$ which obeys 
the two-sided estimate
\begin{eqnarray*}
C_{1}{\bf W}_{1,\,p}\om(x)\leq u(x)\leq C_{2}{\bf W}_{1,\,p}\om(x)
\end{eqnarray*}  
for all $x\in \RR^n$.
\end{theorem}
\begin{remark}{\rm It is of interest to note that statement (ii) in Theorem 
\ref{maintheorem1} is also equivalent to the following capacitary condition 
(see e.g. \cite{V2}).

\noindent{\it {\rm (v)} There exists a constant $C>0$ such that 
$$\om(E)\leq C\, {\rm Cap}_{{\rm\bf I}_{p},\,\frac{q}{q-p+1}}(E)$$
for all compact sets $E\subset\RR^n$.
}
}\end{remark}
\begin{proof}[Proof of Theorem \ref{maintheorem1}]
We show that (i)$\Ra$(ii)$\Ra$(iii)$\Ra$(iv)$\Ra$(i). Note that 
(\ref{testingI_{p}}) is also equivalent to the testing inequality (see e.g. \cite{VW}):
\begin{eqnarray*}
\int_{\RR^{n}}\Big[{\rm\bf I}_{p}\om_{B}(x)\Big]^{\frac{q}{p-1}} \, dx\leq C \, \om(B). 
\end{eqnarray*}
By applying Proposition (\ref{cont}) we deduce (ii)$\Ra$(iii). The 
implication (iv)$\Ra$(i) clearly follows from Theorem \ref{sufficency}.
Therefore, it remains to check (i)$\Ra$(ii) and (iii)$\Ra$(iv).\\

\noindent{\it Proof of} (i)$\Ra$(ii).
Let $u$ be a nonnegative solution of (\ref{ep-equation}) and let $\mu=u^{q}+
\epsilon\om$.
Then $\mu$ is a positive measure such that $\mu\geq u^{q}$, $\mu\geq \epsilon \om$ and
$u(x)\geq \frac{1}{K} {\rm\bf W}_{1,\,p}\mu(x)$ where $K$ is the constant in (\ref{global}). Therefore,
\begin{eqnarray*}
\int_{P}d\mu&\geq&\int_{P}u^{q} \, dx\geq C\int_{P}({\rm\bf W}_{1,\,p}\mu)^{q} \, dx\\
&\geq& C\int_{P}\Big[\sum_{Q\subset P}\Big(\frac{\mu(Q)}{\m{Q}^
{1-\frac{p}{n}}}\Big)^{\frac{1}{p-1}}\chi_{Q}(x)\Big]^{q} \, dx,  
\end{eqnarray*}  
for all dyadic cubes $P$ in $\RR^{n}$. Using this and Proposition \ref{discrete0t}, we get
\begin{eqnarray*}
\sum_{Q\subset P}\Big[\frac{\mu(Q)}{\m{Q}^{1-\frac{p}{n}}}\Big]^{\frac{q}{p-1}}\m{Q}
\leq C \, \mu(P), \hspace*{.3in} P\in \mathcal{D}. 
\end{eqnarray*}
It is known that the preceding condition  is equivalent to the inequality (see \cite[Sec. 3]{V1}) 
$$\norm{{\rm\bf I}_{p}(f)}_{L^{\frac{q}{q-p+1}}(d\mu)}\leq C \, \norm{f}_
{L^{\frac{q}{q-p+1}}(dx)},$$  
where $C$ does not depend on $f \in L^{\frac{q}{q-p+1}}(dx)$. 
Since $\mu\geq\epsilon \, \om$, from this 
we have 
$$\norm{{\rm\bf I}_{p}(f)}_{L^{\frac{q}{q-p+1}
}(d\om)}\leq \epsilon^{\frac{q-p+1}{-q}}C\norm{f}_{L^{\frac{q}{q-p+1}}(dx)}.$$
Therefore, by duality  we obtain the testing inequality (\ref{testingI_{p}}). 
This completes the proof of  (i)$\Ra$(ii). \\

\noindent {\it Proof of} (iii)$\Ra$(iv).
We first claim that  (\ref{testingW}) yields 
\begin{equation}\label{infinite}
\int_{r}^{\infty}\Big[\frac{\om(B_{t}(x))}{t^{n-p}}\Big]^{\frac{1}{p-1}}\frac{dt}{t}\leq 
C \, r^{\frac{-p}{q-p+1}}, 
\end{equation}
where $C$ is independent of $x$ and $r$. Note that for $y \in B_{t}(x)$ 
and  $\tau\geq 2t$, we have $B_{t}(x)\subset B_{\tau}(y)$. Thus,   
\begin{eqnarray*}
{\rm\bf W}_{1,\,p}\om_{B_{t}(x)}(y)&\geq& \int_{2t}^{\infty}\Big(\frac{\om(B_{\tau}(y)\cap B_{t}(x))}{\tau^{n-p}}\Big)^{\frac{1}{p-1}}\frac{d\tau}{\tau}\\
&\geq& C \, \Big(\frac{\om(B_{t}(x))}{t^{n-p}}\Big)^{\frac{1}{p-1}}. 
\end{eqnarray*}
Combining this with (\ref{testingW}), we obtain 
\begin{equation}\label{measureB}
\om(B_{t}(x))\leq C \, t^{n-\frac{pq}{q-p+1}}, 
\end{equation}
which clearly implies  (\ref{infinite}). 

\indent Next, we introduce a decomposition of the Wolff potential
${\rm\bf W}_{1,\,p}$ into its lower and upper parts defined respectively by
\begin{eqnarray*}
{\rm\bf L}_{r}\mu(x)=\int_{r}^{\infty}\Big[\frac{\mu({B_{t}(x)})}{t^{n-p}}\Big]^{\frac
{1}{p-1}}\frac{dt}{t},~~~~ r>0,~x\in\RR^n, 
\end{eqnarray*}
and
\begin{eqnarray*}
{\rm\bf U}_{r}\mu(x)=\int_{0}^{r}\Big[\frac{\mu({B_{t}(x)})}{t^{n-p}}\Big]^{\frac
{1}{p-1}}\frac{dt}{t},~~~~ r>0,~x\in\RR^n.
\end{eqnarray*}
Let $d\nu=({\rm\bf W}_{1,\,p}\om)^{q}dx$. For each $r>0$ let 
$d\mu_{r}=({\rm\bf U}_{r}\om)^{q}dx$ and $d\lam_{r}=({\rm\bf L}_{r}\om)^{q}dx$. Then 
\begin{eqnarray}
\label{sum}
\nu\leq C(q)(\mu_{r}+\lam_{r})\end{eqnarray}
Let $x\in\RR^{n}$ and $B_{r}=B_{r}(x)$. Since ${\rm\bf W}_{1,\,p}({\rm\bf W}_{1,\,p}\om)^{q}={\rm\bf W}_{1,\,p}\nu$, we have to prove that 
$${\rm\bf W}_{1,\,p}\nu(x)=\int_{0}^{\infty}\Big[\frac{\nu(B_{r})}{r^{n-p}}\Big]^{\frac{1}{p-1}}\frac{dr}{r}
\leq C \, {\rm\bf W}_{1,\,p}\om(x).$$
For $r>0$, $t\leq r$ and $y\in B_{r}$ we have $B_{t}(y)\subset B_{2r}$. Therefore it is easy  to 
see that ${\rm\bf U}_{r}\om={\rm\bf U}_{r}\om_{B_{2r}}$ on $B_{r}$. Using this together with (\ref{testingW}), we have
$$\mu_{r}(B_{r})=\int_{B_{r}}({\rm\bf U}_{r}\om)^{q}dx=\int_{B_{r}}({\rm\bf U}_{r}\om_{B_{2r}})^{q}dx\leq C\om(B_{2r}).$$
Hence, 
\begin{eqnarray}
\label{upper}
\int_{0}^{\infty}\Big[\frac{\mu_{r}(B_{r})}{r^{n-p}}\Big]^{\frac{1}{p-1}}\frac{dr}{r}
&\leq& C\int_{0}^{\infty}\Big[\frac{\om(B_{2r})}{r^{n-p}}\Big]^{\frac{1}{p-1}}\frac{dr}{r}\\
&\leq& C' \, {\rm\bf W}_{1,\,p}\om(x).\nonumber
\end{eqnarray}
On the other hand, for $y\in B_{r}$ and $t\geq r$,  we have $B_{t}(y)\subset B_{2t}$, and 
consequently  
\begin{eqnarray}
\label{L_{r}y<x}
{\rm\bf L}_{r}\om(y)&\leq& \int_{r}^{\infty}\Big[\frac{\om(B_{2t})}{t^{n-p}}\Big]^{\frac{1}{p-1}}\frac{dt}{t}\\
&\leq& C\int_{2r}^{\infty}\Big[\frac{\om(B_{t})}{t^{n-p}}\Big]^{\frac{1}{p-1}}\frac{dt}{t}\nonumber\\
&\leq& C \, {\rm\bf L}_{r}\om(x).\nonumber
\end{eqnarray}
Using (\ref{L_{r}y<x}),we obtain
$$\lam_{r}(B_{r})=\int_{B_{r}}({\rm\bf L}_{r}\om(y))^{q}dy\leq C ({\rm\bf L}_{r}\om(x))^{q}\m{B_{r}}.$$
Thus, 
\begin{eqnarray*}
&&\int_{0}^{\infty}\Big[\frac{\lam_{r}(B_{r})}{r^{n-p}}\Big]^{\frac{1}{p-1}}\frac{dr}{r}\\
&\leq& C'\int_{0}^{\infty}({\rm\bf L}_{r}\om(x))^{\frac{q}{p-1}}\Big(\frac{\m{B_{r}}}{r^{n-p}}\Big)^{\frac{1}{p-1}}\frac{dr}{r}\\
&=& C'\int_{0}^{\infty}\Big[\int_{r}^{\infty}\Big(\frac{\om(B_{t})}{t^{n-p}}\Big)^{\frac{1}{p-1}}\frac{dt}{t}\Big]^{\frac{q}{p-1}}
\Big(\frac{\m{B_{r}}}{r^{n-p}}\Big)^{\frac{1}{p-1}}\frac{dr}{r}\\
&=&C'\frac{q}{p-1}\int_{0}^{\infty}\Big[\int_{0}^{r}\Big(\frac{\m{B_{t}}}{t^{n-p}}\Big)^{\frac{1}{p-1}}
\frac{dt}{t}\Big][{\rm\bf L}_{r}\om(x)]^{\frac{q}{p-1}-1}
\Big[\frac{\om(B_{r})}{r^{n-p}}\Big]^{\frac{1}{p-1}}\frac{dr}{r}, 
\end{eqnarray*}
where we have used integration by parts in the last equality. Using now 
(\ref{infinite}), we get
$$\Big[\int_{0}^{r}\Big(\frac{\m{B_{t}}}{t^{n-p}}\Big)^{\frac{1}{p-1}}
\frac{dt}{t}\Big][{\rm\bf L}_{r}\om(x)]^{\frac{q}{p-1}-1}\leq C.$$
Hence, 
\begin{eqnarray}
\label{lower}
\int_{0}^{\infty}\Big[\frac{\lam_{r}(B_{r})}{r^{n-p}}\Big]^{\frac{1}{p-1}}\frac{dr}{r}
&\leq& C'' \int_{0}^{\infty}\Big[\frac{\om(B_{r})}{r^{n-p}}\Big]^{\frac{1}{p-1}}\frac{dr}{r}\\
&=&C'' \, {\rm\bf W}_{1,\,p}\om(x).\nonumber
\end{eqnarray}
Combining (\ref{sum}), (\ref{upper}) and (\ref{lower}) gives  
$${\rm\bf W}_{1,\,p}\nu(x)=\int_{0}^{\infty}\Big[\frac{\nu(B_{r})}{r^{n-p}}\Big]^{\frac{1}{p-1}}\frac{dr}{r}
\leq C \, {\rm\bf W}_{1,\,p}\om(x),$$
for a suitable constant $C$ independent of $\om$. Thus,  (iii) implies (iv) as claimed. This completes 
the proof of the theorem.
\end{proof}

%************************************************************
%************************************************************
\section{Renormalized solutions of quasilinear Dirichlet problems}\label{Om}
Let $\Om$ be a bounded, open subset of $\RR^n$, $n\geq 2$. 
We denote by $\mathcal{M}_{B}(\Om)$ (respectively $\mathcal{M}_{B}^{+}(\Om)$)
the set of all Radon measures (respectively nonnegative Radon measures) 
with bounded variation in $\Om$. 
Let $\mathcal{A}$ be as in Sec. \ref{pre} and let $1<p<\infty$. In this 
section we consider the Dirichlet problem
\begin{eqnarray}\label{Dirichlet}
\left\{\begin{array}{c}
-{\rm div}\mathcal{A}(x, \nabla u)=u^{q}+\om,\\
u\geq 0 \quad {\rm in}\quad\Om,\\
 \hspace*{.1in}u=0 \quad {\rm on}\quad \partial\Om,  
\end{array}
\right.
\end{eqnarray}
where $\om\in \mathcal{M}_{B}^{+}(\Om)$ and  $q>p-1$. 

It is well known that when the data is not regular enough, a solution of nonlinear 
Leray-Lions type equations does not necessarily belong to the Sobolev space ${\rm W}_{0}^{1,\,p}(\Om)$.
Therefore, we use the framework of renormalized solutions which seems proper for such problems 
with measure data (see, e.g.,  \cite{DMOP}).\\
\indent For a measure $\mu$ in $\mathcal{M}_{B}(\Om)$, its positive and  negative parts
are denoted by $\mu^{+}$ and $\mu^{-}$, respectively. We say that  a sequence of measures
$\{\mu_{n}\}$ in $\mathcal{M}_{B}(\Om)$ converges in the narrow topology to  
$\mu \in \mathcal{M}_{B}(\Om)$ if    
$$\lim_{n\rightarrow\infty}\int_{\Om}\varphi \, d\mu_{n}=\int_{\Om}\varphi \, 
d\mu$$
for every bounded and continuous function $\varphi$ on $\Om$. \\
\indent Denote by $\mathcal{M}_{0}(\Om)$ (respectively $\mathcal{M}_{s}(\Om)$) the set of all measures in  $\mathcal{M}_{B}(\Om)$ 
which are continuous (respectively singular) with 
respect to the capacity ${\rm cap}_{1,\,p}(\cdot,\Om)$. Here ${\rm cap}_{1,\,p}(\cdot,\Om)$ is the
capacity relative to the domain $\Om$ defined by
\begin{equation}\label{cap1p}
{\rm cap}_{1,\,p}(E,\Om)=\inf\Big\{\int_{\Om}\m{\nabla \phi}^{p}dx: \phi\in C_{0}^{\infty}
(\Om), \phi\geq 1 {\rm ~on~} E \Big\}
\end{equation}
for any compact set $E\subset\Om$. Recall that, for every measure $\mu$ in $\mathcal{M}_{B}
(\Om)$, there exists a unique pair of measures ($\mu_{0}, \mu_{s}$) with $\mu_{0}\in \mathcal
{M}_{0}(\Om)$ and  $\mu_{s}\in \mathcal{M}_{s}(\Om)$, such that $\mu=\mu_{0}+\mu_{s}$. If
$\mu$ is nonnegative, then so are $\mu_{0}$ and $\mu_{s}$ (see \cite[Lemma 2.1]{FST}).  

For $k>0$ and for $s\in \RR$ we denote by $T_{k}(s)$ the truncation $T_{k}(s)=\max\{-k,\min
\{k,s\}\}$. Recall also from \cite{BBG} that if $u$ is a measurable function  on $\Om$ which is 
finite almost everywhere and satisfies $T_{k}(u)\in W_{0}^{1,\,p}(\Om)$ for every $k>0$, then there exists
a measurable function $v:\Om\rightarrow \RR^n$ such that 
\begin{equation}
\nabla T_{k}(u)=v\chi_{\{\m{u}<k\}} {\rm ~~~almost ~everywhere~ in~} \Om, {\rm for~ all~} k>0. 
\nonumber
\end{equation}
Moreover, $v$ is unique up to almost everywhere equivalence. We define the gradient $Du$ of 
$u$ as this function $v$, and set $Du=v$.   \\
\indent In \cite{DMOP}, several equivalent definitions of renormalized solutions
are given. In what follows, we will need the following ones.
\begin{definition}\label{rns1}{\rm
Let $\mu$ be in $\mathcal{M}_{B}(\Om)$. Then $u$ is said to be a renormalized 
solution of 
\begin{eqnarray}\label{liDirichlet}
\left\{\begin{array}{c}
-{\rm div}\mathcal{A}(x, \nabla u)=\mu \quad {\rm in}\quad\Om,\\
 \hspace*{-.5in} u=0 \quad {\rm on}\quad \partial\Om,  
\end{array}
\right.
\end{eqnarray}
if the following conditions hold:\\
\noindent (a) The function $u$ is measurable and finite almost everywhere, and $T_{k}(u)$ belongs 
to $W_{0}^{1,\,p}(\Om)$ for every $k>0$.\\
\noindent (b) The gradient $Du$ of $u$ satisfies $\m{Du}^{p-1}\in L^{q}(\Om)$ for all $q<\frac{n}{n-1}$.\\
\noindent (c) If $w$ belongs to $W_{0}^{1,\,p}(\Om)\cap L^{\infty}(\Om)$ and if there exist $k>0$, 
$w^{+\infty}$ and $w^{-\infty}$ in $W^{1,\,r}(\Om)\cap L^{\infty}(\Om)$, with $r>N$, such that 
\begin{equation*}
\left\{\begin{array}{c}
w=w^{+\infty} \hspace*{.2in}{\rm a.e. ~on~ the~ set}~~ \{u>k\},\\
w=w^{-\infty} \hspace*{.2in}{\rm a.e. ~on ~the ~set}~~ \{u<-k\}  
\end{array}
\right.
\end{equation*}
then
\begin{equation}\label{test}
\int_{\Om}\mathcal{A}(x, Du)\cdot\nabla wdx=\int_{\Om}w d\mu_{0}+\int_{\Om}w^{+\infty}
d\mu_{s}^{+}-\int_{\Om}w^{-\infty}d\mu_{s}^{-}.
\end{equation}
}\end{definition}      
\begin{definition}\label{rns}{\rm
Let $\mu$ be in $\mathcal{M}_{B}(\Om)$. Then  $u$ is a renormalized 
solution of (\ref{liDirichlet})
if $u$ satisfies (a) and (b) in Definition \ref{rns1}, and if the following conditions hold:\\
\noindent (c) For every $k>0$ there exist two nonnegative measures in $\mathcal{M}_{0}(\Om)$, $\lambda^{+}_{k}$ and 
$\lambda^{-}_{k}$, concentrated on the sets $\{u=k\}$ and $\{u=-k\}$, respectively, such that 
$\lambda^{+}_{k}\rightarrow \mu_{s}^{+}$ and  
$\lambda^{-}_{k}\rightarrow \mu_{s}^{-}$ in the narrow topology 
of measures.  
\noindent (d) For every $k>0$
\begin{equation}\label{truncate}
\int_{\{\m{u}<k\}}\mathcal{A}(x, Du)\cdot\nabla\varphi dx= \int_{\{\m{u}<k\}}\varphi
d\mu_{0} +  \int_{\Om}\varphi d\lambda_{k}^{+} -  \int_{\Om}\varphi d\lambda_{k}^{-}
\end{equation}
for every $\varphi$ in ${\rm W}_{0}^{1,\,p}(\Om)\cap L^{\infty}(\Om)$.
}\end{definition}
\begin{remark}\label{quasi}{\rm
By  \cite[Remark 2.18]{DMOP}, if $u$ is a renormalized solution of (\ref{liDirichlet})
then (the ${\rm cap}_{1,\,p}$-quasi continuous representative
of) $u$ is finite ${\rm cap}_{1,\,p}$-quasieverywhere.
Therefore, $u$ is finite $\mu_{0}$-almost everywhere.    
}\end{remark}
\begin{remark}\label{approx}{\rm
By (\ref{truncate}), if $u$ is a renormalized solution  of (\ref{liDirichlet}) then 
\begin{equation}\label{4.4}
-{\rm div}\mathcal{A}(x,\nabla T_{k}(u))=\mu_{k} \quad {\rm in}\quad \Om, 
\end{equation}
where 
$$\mu_{k}=\chi_{\{\m{u}<k\}}\mu_{0}+\lambda^{+}_{k}-\lambda^{-}_{k}.$$
Since 
$T_{k}(u)\in W_{0}^{1,\,p}(\Om)$, by (\ref{2.3}) we see that 
$-{\rm div}\mathcal{A}(x, \nabla T_{k}(u))$ and hence $\mu_{k}$ belongs to 
the dual space $W^{-1,\,p'}(\Om)$ of $W_{0}^{1,\,p}(\Om)$. Moreover, by Remark \ref{quasi}, 
$\m{u}<\infty$ $\mu_{0}$-almost everywhere and hence $\chi_{\{\m{u}<k\}} \rightarrow 
\chi_{\Om}$ $\mu_{0}$-almost everywhere as $k\rightarrow\infty$. Therefore, by the  monotone
convergence theorem,  $\mu_{k}$ converges to 
$\mu$ in the narrow topology of measures.
}\end{remark}
\begin{remark}\label{nonnegative}{\rm 
If $\mu \geq 0$, i.e., $\mu\in \mathcal{M}_{B}^{+}(\Om)$, and $u$ is a renormalized solution of 
(\ref{liDirichlet}) then $u$ is nonnegative. To see this, for each $k>0$ we  ``test" (\ref{test}) 
 with $w=-T_{k}(u^-)$ where $u^-=-\min\{u,0\}$, $w^{+\infty}=0$ and $w^{-\infty}=-k$: 
\begin{eqnarray*}
-\int_{\Om}\mathcal{A}(x,Du)\cdot\nabla T_{k}(u^-)dx&=&-\int_{\Om}T_{k}(u^-)d\mu_{0} 
+\int_{\Om}k d\mu_{s}^-\\
&=& -\int_{\Om}T_{k}(u^-)d\mu_{0}\leq 0, 
\end{eqnarray*}
since $\mu_{s}^{-}=0$. Thus using (\ref{2.3}) we get
$$\int_{\Om}\m{\nabla T_{k}(u^-)}^pdx\leq 0$$
for every $k>0$. Therefore $u^- = 0$ a.e., i.e., $u$ is nonnegative. 
}\end{remark}
\begin{remark}\label{rm5.6}{\rm Let $\mu\in {\mathcal M}_{B}^+(\Om)$ and let $u$ be a renormalized solution
of (\ref{liDirichlet}). Since $u^{-}=0$ a.e. 
(by Remark \ref{nonnegative}) and hence $u^-=0$ ${\rm cap}_{1,p}$-quasi everywhere 
(see \cite[Theorem 4.12]{HKM}), in Remark \ref{approx}  we may take $\lambda_{k}^{-}=0$, 
 and thus $\mu_{k}$ is nonnegative.
Hence by (\ref{4.4}) and Proposition \ref{pro2.1}, the functions
$v_{k}$ defined by $v_{k}(x)={\rm ess}\liminf_{y\rightarrow x} T_{k}(u)(y)$ are 
$\mathcal{A}$-superharmonic and increasing. 
Using Lemma 7.3 in \cite{HKM}, it is 
then easy to see that 
$v_{k}\rightarrow v$  as $k\rightarrow\infty$ everywhere in $\Om$  for some 
$\mathcal{A}$-superharmonic function $v$ on $\Om$ such that $v=u$ a.e.. In other words,
$v$ is an  $\mathcal{A}$-superharmonic  representative  of $u$. 
}\end{remark}
\begin{remark}\label{repre}{\rm When we are dealing with pointwise values of a renormalized 
solution u to the problem (\ref{liDirichlet}) with measure data $\mu\geq 0$, we always identify     
u with its $\mathcal{A}$-superharmonic representative mentioned in Remark \ref{rm5.6}.
}\end{remark}
In Theorem \ref{om-estimate} below, we give a global pointwise potential estimate
for renormalized solutions on a bounded domain $\Om$, whose proof is based on its  
local counterpart given in Theorem \ref{potential} and the following 
lemma. 

\begin{lemma}\label{compareuw}
Suppose that $u$ is a renormalized solution of the problem (\ref{liDirichlet}) with data 
$\mu\in {\mathcal M}_{B}^{+}(\Om)$. Let $B=B(x_{0},2{\rm diam}(\Om))$ be a ball centered at $x_{0}\in\Om$. 
Then there exists a nonnegative $\mathcal{A}$-superharmonic function $w$ on $B$ 
such that $u\leq w$ on $\Om$, and  
\begin{eqnarray*}
\left\{\begin{array}{c}
-{\rm div}\mathcal{A}(x, \nabla w)=\mu \quad {\rm in}\quad B,\\
 \norm{w}_{L^{p-1}(B)}\leq C R^{\frac{p}{p-1}}\mu(\Om)^{\frac{1}{p-1}}.  
\end{array}
\right.
\end{eqnarray*}
\end{lemma}
\begin{proof}
Let $u_{k}=\min\{u,k\}$, and let 
$$\mu_{k}=\chi_{\{u<k\}}\mu_{0}+\lambda_{k}^{+}$$
be as in Remark \ref{approx} (note that $\lambda_{k}^{-}=0$ by Remark \ref{rm5.6}). We see that $u_{k}\in W_{0}^{1,p}(\Om)$
is the unique solution of problem (\ref{liDirichlet}) with data $\mu_{k}$. Since 
$\mu_{k}$ is continuous with respect to the capacity ${\rm cap}_{1,\,p}(\cdot,B)$, we have 
a unique renormalized (or entropy) solution $w_{k}$ to the problem 
\begin{eqnarray*}
\left\{\begin{array}{c}
-{\rm div}\mathcal{A}(x, \nabla w_{k})=\mu_{k} \quad {\rm in}\quad B,\\
w_{k}=0 \quad {\rm on}\quad \partial B.  
\end{array}
\right.
\end{eqnarray*}
We now extend $u_{k}$ by zero outside $\Om$, and set 
$$\Psi=\min\{w_{k}-u_{k},0\}=\min\{\min\{w_{k},k\}-u_{k},0\}.$$
Note that $\Psi\in W_{0}^{1,\,p}(\Om)\cap W_{0}^{1,\,p}(B)\cap L^{\infty}(B)$ 
since $\m{\Psi}\leq u_{k}$. Then using $\Psi$ as a test function we have 
\begin{eqnarray*}
0&=&\int_{B}\mathcal{A}(x,\nabla w_{k})\cdot\nabla\Psi dx -\int_{\Om}\mathcal{A}
(x,\nabla u_{k})\cdot\nabla\Psi dx\\
&=& \int_{B\cap\{w_{k}<u_{k}\}}\mathcal{A}(x,\nabla w_{k})\cdot\nabla\Psi dx -
\int_{B\cap\{w_{k}<u_{k}\}}\mathcal{A}(x,\nabla u_{k})\cdot\nabla\Psi dx\\
&=& \int_{B\cap\{w_{k}<u_{k}\}}[\mathcal{A}(x,\nabla w_{k})-\mathcal{A}(x,\nabla 
u_{k})]\cdot(\nabla w_{k}-\nabla u_{k})dx.
\end{eqnarray*} 
Thus $\nabla w_{k}=\nabla u_{k}$ a.e. on the set $B\cap \{w_{k}<
u_{k}\}$ by hypothesis (\ref{2.5}) on $\mathcal{A}$. Hence $\Psi=0$ a.e., i.e., 
\begin{equation}\label{uw}
u_{k}\leq w_{k}\quad{\rm a.e.}
\end{equation}
Since $\mu_{k}$ converges to $\mu$ in the narrow
topology of measures on $\Om$ (and hence also on $B$),
arguing as in the proof of 
\cite[Theorem 2.4]{KM1}, we can find a 
subsequence $\{w_{k_{j}}\}$ of $\{w_{k}\}$ such that $w_{k_{j}}\rightarrow w$ a.e., where 
$w$ is a nonnegative $\mathcal{A}$-superharmonic function on $B$ such that
\begin{eqnarray*}
-{\rm div}\mathcal{A}(x, \nabla w)=\mu \quad {\rm in}\quad B.
\end{eqnarray*}
By (\ref{uw}) we have $u\leq w$ a.e. on $\Om$, and hence $u\leq w$ everywhere on
$\Om$ due to Remark \ref{repre} and Proposition \ref{pro2.1}. Note that for $p<n$
we have
\begin{eqnarray*}
\norm{w_{k}}_{L^{\frac{n(p-1)}{n-p},\,\infty}(B)}\leq C \, \mu_{k}(\Om)^{\frac{1}{p-1}},
\end{eqnarray*}    
for some constant $C$ independent of $R$ and $k$ (see \cite[Theorem 4.1]{DMOP} or 
\cite[Lemma 4.1]{BBG}). Thus 
\begin{eqnarray}\label{w(B)}
\norm{w_{k}}_{L^{p-1}(B)}\leq C R^\frac{p}{p-1}\mu_{k}(\Om)^{\frac{1}{p-1}}.
\end{eqnarray}
The inequality (\ref{w(B)}) also holds for $p\geq n$, see for example \cite[Lemma 2.1]{Gre}. 
Finally, using Fatou's lemma and (\ref{w(B)}), we obtain 
\begin{eqnarray*}
\norm{w}_{L^{p-1}(B)}\leq C R^\frac{p}{p-1}\mu(\Om)^{\frac{1}{p-1}}.
\end{eqnarray*} 
This completes the proof of the lemma. 
\end{proof}

\begin{theorem}\label{om-estimate} Suppose that $u$ is a renormalized solution of the problem (\ref{liDirichlet}) 
with data $\mu\in {\mathcal M}_{B}^{+}(\Om)$. Let $R={\rm diam}(\Om)$. Then there is a constant $K$ independent
of $\mu$ and $R$ such that
\begin{eqnarray}\label{uestimate}
u(x)\leq K \, {\rm\bf W}_{1,\,p}^{2R}\mu(x), 
\end{eqnarray}
for all $x$ in $\Om$.
\end{theorem}
\begin{proof}
Let $w$ and $B$ be as in Lemma \ref{compareuw}. Fix  $x\in\Om$. We denote by 
$d(x)$ the distance from $x$ to the boundary $\partial B$ of $B$. By Theorem \ref{potential},
Lemma \ref{compareuw}, and the fact that $d(x)\geq R$, we have
\begin{eqnarray}\label{westimate}
w(x)&\leq& C\,{\rm\bf W}_{1,\,p}^{2d(x)/3}\mu(x) + C\inf_{B(x,d(x)/3)}w\\
&\leq& C\,{\rm\bf W}_{1,\,p}^{2R}\mu(x) +
C d(x)^{\frac{-n}{p-1}}\norm{w}_{L^{p-1}(B)}\nonumber\\
&\leq& C\,{\rm\bf W}_{1,\,p}^{2R}\mu(x) + C R^{\frac{p-n}{p-1}}\mu(\Om)^{\frac{1}{p-1}}\nonumber\\
&\leq& C\,{\rm\bf W}_{1,\,p}^{2R}\mu(x).\nonumber  
\end{eqnarray}
Therefore, from (\ref{westimate}) and Lemma \ref{compareuw}  we get  the desired 
inequality (\ref{uestimate}).   
\end{proof}

\begin{theorem}\label{sufficency2}
Let $\om \in {\mathcal M}_{B}^{+}(\Om)$. Let $p>1 $ and $q>p-1$. 
Suppose that $R={\rm diam}(\Om)$, and 
\begin{equation}
\label{W(Wom)^{q}2}
{\rm\bf W}_{1,\,p}^{2R}({\rm\bf W}_{1,\,p}^{2R}\om)^{q}\leq C\,{\rm\bf W}_{1,\,p}^{2R}\om<\infty 
\quad  {\rm a.e.,}  
\end{equation}
where 
$$C \leq \Big(\frac{q-p+1}{qK\max\{1,2^{p'-2}\}}\Big)^{q(p'-1)}\Big(\frac{p-1}{q-p+1}\Big),$$
and $K$ is the constant in Theorem \ref{om-estimate}.
 Then there is a renormalized solution $u\in L^{q}(\Om)$ to the  Dirichlet 
problem 
\begin{eqnarray}
\label{suff2}
\left\{\begin{array}{c}
-{\rm div}\mathcal{A}(x, \nabla u)=u^{q}+\om \quad{\rm in}\quad\Om,\\
\hspace*{.2in}u=0 \quad {\rm on}\quad \partial\Om
\end{array}
\right.
\end{eqnarray}
such that
$$ u(x)\leq M \, {\rm\bf W}_{1,\,p}^{2R}\om(x),$$
for all $x$ in $\Om$, where the constant $M$ depends  only $n, p, q$ and the 
structural constants $\alpha$ and $\beta$. 
\end{theorem}
\begin{proof}
Let $\{u_{k}\}_{k\geq 0}$ be a sequence of renormalized solutions defined inductively for the 
following Dirichlet problems:
\begin{eqnarray}
\left\{\begin{array}{c}
-{\rm div}\mathcal{A}(x, \nabla u_{0})=\om \quad{\rm in~}\quad\Om,\\
\hspace*{.2in}u_{0}=0 \quad {\rm on} \quad\partial\Om,
\end{array}
\right.
\end{eqnarray}
and 
\begin{eqnarray}
\label{suff3}
\left\{\begin{array}{c}
-{\rm div}\mathcal{A}(x, \nabla u_{k})=u_{k-1}^{q}+\om \quad{\rm in}\quad\Om,\\
\hspace*{.2in}u_{k}=0 \quad{\rm on}\quad \partial\Om.
\end{array}
\right.
\end{eqnarray}
for $k\geq 1$. By Theorem \ref{om-estimate} we have 
$$u_{0}\leq K\, {\rm W}_{1,\,p}^{2R}\om, \qquad 
u_{k}\leq K\, {\rm W}_{1,\,p}^{2R}(u_{k-1}^{q}+\om).$$
Thus by arguing as in the proof of Theorem \ref{sufficency},  we obtain a constant $M>0$ such that 
\begin{equation}\label{4.42}
u_{n}\leq M\, {\rm W}_{1,\,p}^{2R}\om<\infty \quad {\rm a.e.}
\end{equation}
for all $n\geq 0$. By passing to a subsequence (see \cite[Theorem 1.17]{KM1} 
or \cite[Sec. 5.1]{DMOP}), we can assume 
that  $u_{n}\rightarrow u$ a.e. on $\Om$ for some nonnegative function $u$. Note that by 
(\ref{4.42})
\begin{equation*}
u\leq M\,{\rm W}_{1,\,p}^{2R}\om<\infty\quad a.e.
\end{equation*}
and $u_{n}^{q}\rightarrow u^{q}$ in ${\rm L}^{1}(\Om)$. Finally, in view of 
(\ref{suff3}), the stability result
in \cite[Theorem 3.4]{DMOP}  asserts that $u$ is a renormalized solution of
(\ref{suff2}), which proves the theorem.  
\end{proof}

Let $\mathcal{Q}=\{Q\}$ be a Whitney decomposition of $\Om$, i.e., $\mathcal{Q}$ is a disjoint
subfamily of the family of dyadic cubes in $\RR^n$ such that $\Om=\cup_{Q\in\mathcal{Q}}Q$,
where we can assume that $2^{5}{\rm diam}(Q)\leq{\rm dist(Q, \partial\Om)}\leq 2^{7}
{\rm diam}(Q)$. Let $\{\phi_{Q}\}_{Q\in\mathcal{Q}}$ be a partition of unity associated 
with the Whitney decomposition of $\Om$ above: $0\leq\phi_{Q}\in C_{0}^{\infty}(Q^{*})$,
$\phi_{Q}\geq 1/C(n)$ on $\overline{Q}$, $\sum_{Q}\phi_{Q}=1$ and $\m{D^{\gamma}\phi_{Q}}
\leq A_{\gamma}({\rm diam}(Q))^{-\m{\gamma}}$
for all multi-indices $\gamma$. Here $Q^{*}=(1+\epsilon)Q$, $0<\epsilon<\frac{1}{4}$ and 
$C(n)$ is a positive constant depending only on $n$ such that each point in  $\Om$ is 
contained in at most $C(n)$ of the cubes $Q^{*}$ (see \cite{St1}).  
 
\begin{theorem} \label{LCE} Let $\om$ be a locally finite nonnegative  measure 
on an open (not necessarily bounded) set $\Om$.  Let $p>1$ and $q>p-1$. 
Suppose that 
there exists a nonnegative $\mathcal{A}$-superharmonic function $u$ in $\Om$ such that 
\begin{equation}\label{nonboundary}
-{\rm div\mathcal{A}(x,\nabla u)}=u^q +\om \quad {\rm in}\quad\Om.
\end{equation} 
Then, for each cube $P\in\mathcal{Q}$ and compact set $E\subset \Om$, 
\begin{eqnarray}\label{oneachcube1}
\mu_{P}(E)\leq C \, {\rm Cap}_{{\rm\bf I}_{p},\,\frac{q}{q-p+1}}(E) 
\end{eqnarray}
if $\frac{pq}{q-p+1}<n$, and 
\begin{eqnarray}\label{oneachcube1'}
\mu_{P}(E)\leq C(P) \, {\rm Cap}_{{\rm\bf G}_{p},\,\frac{q}{q-p+1}}(E) 
\end{eqnarray}
if $\frac{pq}{q-p+1}\geq n$. Here 
$d\mu =u^{q}dx+d\om$, and the constant $C$ in (\ref{oneachcube1}) is 
independent of $P$ and $E\subset \Om$, but the constant $C(P)$ in
 (\ref{oneachcube1'}) may depend on the side length of $P$.

Moreover, if $\frac{pq}{q-p+1}<n$ and $\Om$ is a bounded $C^{\infty}$-domain, then 
$$\mu(E)\leq C \, {\rm cap}_{p,\,\frac{q}{q-p+1}}(E,\Om),$$
for all compact sets $E\subset\Om$, where ${\rm cap}_{p,\,\frac{q}{q-p+1}}(E,\Om)$ 
is defined by (\ref{capomega}).
\end{theorem}
\begin{proof}
Let $P$ be a fixed dyadic cube in $\mathcal{Q}$. For a dyadic
cube $P'\subset P$ we have 
$${\rm dist}(P',\partial\Om)\geq{\rm dist}(P,\partial\Om)\geq 2^{5}{\rm diam}(P) \geq 2^{5}
{\rm diam}(P').$$
The lower estimate in Theorem \ref{potential} then yields 
\begin{eqnarray*}
u(x)&\geq& C \, {\rm\bf W}_{1,\, p}^{2^{3}{\rm diam}(P')}\mu(x)\\
&\geq& C \sum_{k=0}^{\infty}\int_{2^{-k+2}{\rm diam}(P')}^{2^{-k+3}{\rm diam}(P')}
\Big[ \frac{\mu(B_{t}(x))}{t^{n-p}
}\Big]^{\frac{1}{p-1}}\frac{dt}{t}\\
&\geq& C \sum_{Q\subset P'}\Big[ \frac{\mu(Q)}{\m{Q}^{1-p/n}}\Big]^{\frac{1}{p-1}}\chi_{Q}(x),
\end{eqnarray*}
for all $x\in P'$. Thus it follows from Proposition \ref{discrete0t} that 
\begin{equation}\label{forP}
\sum_{Q\subset P'}\Big[\frac{\mu(Q)}{\m{Q}^{1-p/n}}\Big]^{\frac{q}{p-1}}\m{Q}
\leq C\int_{P'}u^{q}dx \leq C\mu(P'), \quad\quad  P'\subset P.
\end{equation}
Hence
\begin{equation}\label{estimateformu1}
\mu(P')\leq C\m{P'}^{1-\frac{pq}{n(q-p+1)}},\quad\quad  P'\subset P.
\end{equation}
To get a better estimate for $\mu(P')$ in the case $\frac{pq}{q-p+1}=n$, we observe 
that (\ref{forP}) is a dyadic Carleson condition. Thus by the dyadic  
Carleson imbedding theorem (see, e.g.,  \cite{NTV}, \cite{V1}) we obtain, 
for $\frac{pq}{q-p+1}=n$,
\begin{equation}\label{Carleson}
\sum_{Q\subset P} \mu(Q)^{\frac{q}{p-1}} \Big[\frac{1}{\mu(Q)}\int_{Q}f d\mu\Big]
^{\frac{q}{p-1}}\leq C \int_{P}f^{\frac{q}{p-1}}d \mu,
\end{equation}
where  
$f\in L^{\frac{q}{p-1}}(d\mu_{P})$, $f\geq 0$.  From (\ref{Carleson}) with 
$f=\chi_{P'}$, one gets
\begin{equation}\label{estimateformu2}
\mu(P')\leq C \Big(\log\frac{2^n\m{P}}{\m{P'}}\Big)^{\frac{1-p}{q-p+1}},\quad  P'\subset P,
\end{equation}
if $\frac{pq}{q-p+1}=n$. 
 Now let $P'$ be a dyadic cube in $\RR^n$. From Wolff's inequality for Riesz potentials
(see \cite{HW}) we have
\begin{eqnarray}\label{KS1}
&&\int_{\RR^n}({\rm\bf I}_{p}\mu_{P'\cap P})^{\frac{q}{p-1}}dx\\
&\leq& C 
\sum_{Q\in \mathcal{D}}\Big[\frac{\mu_{P}(P'\cap Q)}{\m{Q}^{1-p/n}}\Big]^{\frac{q}{p-1}}\m{Q}
\nonumber\\
&=& C \sum_{Q\subset P'}\Big[\frac{\mu_{P}(Q)}{\m{Q}^{1-p/n}}\Big]^{\frac{q}{p-1}}\m{Q} 
+ C\sum_{P'\varsubsetneq Q}\Big[\frac{\mu_{P}(P')}{\m{Q}^{1-p/n}}\Big]^{\frac{q}{p-1}}\m{Q}
\nonumber.
\end{eqnarray}
Thus, for $\frac{pq}{q-p+1}<n$, by combining (\ref{forP}) and (\ref{KS1}) we deduce 
\begin{equation}\label{KS2}
\int_{\RR^n} ({\rm\bf I}_{p}\mu_{P'\cap P})^{\frac{q}{p-1}} \, dx\leq C \, 
\mu_{P}(P').
\end{equation} 
In the case $\frac{pq}{q-p+1}\geq n$, a similar argument using (\ref{forP}), (\ref{estimateformu1}), 
(\ref{estimateformu2}) and Wolff's inequality for Bessel potentials:
$$\int_{\RR^n}({\rm\bf G}_{p}\mu_{P'\cap P})^{\frac{q}{p-1}}dx \leq C(P) 
\sum_{Q \in \mathcal{D},\, Q \subset P}\Big[\frac{\mu_{P}(P'\cap Q)}{\m{Q}^{1-p/n}}\Big]^
{\frac{q}{p-1}}\m{Q}, 
$$
(see \cite{AH}), also gives
\begin{equation}\label{KS3}
\int_{\RR^n}({\rm\bf G}_{p}\mu_{P'\cap P})^{\frac{q}{p-1}}dx\leq C(P) 
\mu_{P}(P'),
\end{equation} 
where the constant $C(P)$ may depend on the side-length of $P$. 
Note that (\ref{KS2}), which holds for all dyadic
cubes $P'\subset \RR^n$, is the well-known Kerman-Sawyer condition (see \cite{KS}). Therefore,
$$\norm{{\rm \bf I}_{p}(f)}^{\frac{q}{q-p+1}}_{L^{\frac{q}{q-p+1}}(d\mu_{P})}
\leq C \norm{f}^{\frac{q}{q-p+1}}_{L^{\frac{q}{q-p+1}}(dx)}$$ 
for all $f\in L^{\frac{q}{q-p+1}}(\RR^n)$  which is equivalent to the capacitary
condition: 
$$\mu_{P}(E)\leq C \, {\rm Cap}_{{\rm\bf I}_{p},\,\frac{q}{q-p+1}}(E)$$
for all compact sets $E\subset\RR^n$. Thus we obtain (\ref{oneachcube1}). The inequality
(\ref{oneachcube1'}) is proved in the same way using (\ref{KS3}). 
From (\ref{oneachcube1}) and the definition of ${\rm cap}_{p,\,\frac{q}{q-p+1}}(\cdot,\Om)$, we see that, 
for each cube $P\in\mathcal{Q}$,
$$\mu_{P}(E)\leq C {\rm cap}_{p,\,\frac{q}{q-p+1}}(E \cap P,\Om)$$
for all compact sets $E\subset\Om$. Thus
\begin{eqnarray*}
\mu(E)&\leq& \sum_{P\in\mathcal{Q}}\mu_{P}(E)\\
&\leq& C \sum_{P\in\mathcal{Q}}{\rm cap}_{p,\,
\frac{q}{q-p+1}}(E\cap P,\Om)\\
&\leq& C \, {\rm cap}_{p,\,
\frac{q}{q-p+1}}(E, \Om),
\end{eqnarray*}
where the last inequality follows from the quasiadditivity of the capacity $
{\rm cap}_{p,\frac{q}{q-p+1}}(\cdot,\Om)$ which is considered in the next theorem.
\end{proof}
\begin{remark} \label{loglocal}{\rm Let $B_{R}$ be a ball such that $B_{2R}\subset
\Om$. It is easy to see that there exists a constant $c>0$ such that $\ell(P)\geq c R$ 
for any Whitney cube P that intersects $B_{R}$. On the other hand, if $B_{r}$ is a ball in 
$B_{R}$ then we can find at most $N$ dyadic cubes $P_{i}$ with 
$cr/4\leq\ell(P_{i})<cr/2$ that cover $B_{r}$, where $N$ depends only on $n$. Thus if $\frac{pq}
{q-p+1}=n$ then from 
(\ref{estimateformu2}) we see that 
$$\mu(B_{r})\leq C (\log\tfrac{2R}{r})^{\frac{1-p}{q-p+1}}$$  
for all balls $B_{r}\subset B_{R}$. Here the constant $C$ is independent of $R$ and $r$.
}\end{remark}
\begin{theorem}\label{quasiadditivity} Suppose that $\Om$ is a $C^{\infty}$-domain in $\RR^n$. Then 
there exists a constant $C>0$ such that 
$$\sum_{Q\in\mathcal{Q}}{\rm cap}_{p,\,
\frac{q}{q-p+1}}(E\cap Q,\Om)\leq C {\rm cap}_{p,\,
\frac{q}{q-p+1}}(E,\Om)$$
for all compact sets $E\subset\Om$.
\end{theorem}
\begin{proof}
Obviously, we may assume that ${\rm cap}_{p,\,
\frac{q}{q-p+1}}(E,\Om)>0$. Then by definition there exists $f\in C_{0}^{\infty}(\Om)$,
$f\geq 1$ on $E$ such that 
$$2 \, {\rm cap}_{p,\,\frac{q}{q-p+1}}(E,\Om)\geq\norm{f}^{\frac{q}{q-p+1}}_{W^{p,\,
\frac{q}{q-p+1}}(\RR^n)}.$$
By the refined localization principle on the smooth domain $\Om$ for the function space 
$W^{p,\,\frac{q}{q-p+1}}$ we have 
$$\norm{f}^{\frac{q}{q-p+1}}_{W^{p,\,\frac{q}{q-p+1}}(\RR^n)}\geq C \sum_{Q\in\mathcal{Q}}
\norm{f\phi_{Q}}^
{\frac{q}{q-p+1}}_{W^{p,\,\frac{q}{q-p+1}}(\RR^n)},$$ (see e.g. \cite[Theorem 5.14]{Tri}). Thus
\begin{equation}\label{sumphicube}
\sum_{Q\in\mathcal{Q}}\norm{f\phi_{Q}}^{\frac{q}{q-p+1}}_{W^{p,\,\frac{q}{q-p+1}}(\RR^n)}\leq 
C{\rm cap}_{p,\,\frac{q}{q-p+1}}(E,\Om).
\end{equation}
Note that for $x\in E\cap \overline{Q}$,
$$f\phi_{Q}\geq \phi_{Q}\geq 1/C(n).$$
Hence by definition we have 
$${\rm cap}_{p,\,\frac{q}{q-p+1}}(E\cap \overline{Q},\Om)\leq C \norm{f\phi_{Q}}^{\frac{q}{q-p+1}}_{W^{p,\,
\frac{q}{q-p+1}}(\RR^n)}.$$
From this and (\ref{sumphicube}) we deduce the desired inequality.
\end{proof}
\begin{theorem}\label{maintheorem3}
Let $\om$ be a  measure in ${\mathcal M}_{B}^+(\Om)$ with compact support in $\Om$. 
Let $p>1$, $q>p-1$ and $R={\rm diam}(\Om)$. Then the 
following statements are equivalent.\\
{\rm (i)} There exists a nonnegative renormalized solution $u\in 
L^{q}(\Om)$
 to the equation 
\begin{equation}
\label{ep-equation2}
\left\{\begin{array}{c}
-{\rm div}\mathcal{A}(x, \nabla u)=u^{q}+\ep\om \quad{\rm in}\quad \Om,\\
u=0 \quad{\rm on}\quad \partial\Om 
\end{array}
\right.
\end{equation} 
for some $\ep>0$.\\
{\rm (ii)} The testing inequality 
\begin{equation}
\label{testingI2}
\int_{B}({\rm\bf G}_{p}\om_{B})^{\frac{q}{p-1}} \, dx\leq C \, \om(B)
\end{equation}
holds for all balls $B$ such that $B\cap {\rm supp}\om\not=\emptyset$.\\
{\rm (iii)} For all compact sets $E\subset{\rm supp}\om$,
$$\om(E)\leq C \, {\rm Cap}_{{\rm\bf G}_{p},\,\frac{q}{q-p+1}}(E).$$
{\rm (iv)} The testing inequality 
\begin{equation}
\label{testingW2}
\int_{B}\Big[{\rm\bf W}_{1,\,p}^{2R}\om_{B}(x)\Big]^{q} \, dx\leq C \, \om(B)
\end{equation}
holds for all balls $B$ such that $B\cap {\rm supp}\om\not=\emptyset$ .\\
{\rm (v)} There exists a constant C such that
\begin{equation}
\label{pointwise2}
\hspace*{.3in}{\rm\bf W}_{1,\,p}^{2R}({\rm\bf W}_{1,\,p}^{2R}\om)^{q}\leq C \,  {\rm\bf W}_{1,\,p}^{2R}
\om < \infty \quad {\rm a.e. ~on~} \Om.
\end{equation}
 Moreover, if the constant $C$ in (\ref{pointwise2}) satisfies
\begin{eqnarray*} 
C\leq\Big(\frac{q-p+1}{qK\max\{1,2^{p'-2}\}}\Big)^{q(p'-1)}\Big(\frac{p-1}{q-p+1}\Big), 
\end{eqnarray*} 
where $K$ is the constant in Theorem \ref{om-estimate}, then the equation 
(\ref{ep-equation2}) 
has a solution $u$ with $\ep=1$ which obeys the estimate
\begin{eqnarray*}
u(x)\leq M \,{\bf W}_{1,\,p}^{2R}
\om(x)
\end{eqnarray*}  
for all $x\in \Om$.
\end{theorem}

\begin{proof} It is well known that  statements (ii) and (iii) above are equivalent 
(see e.g. \cite{V2}). Thus it remains 
to prove that (i)$\Longrightarrow$ (ii)$\Longrightarrow$ (iv) $\Longrightarrow$(v) $\Longrightarrow$(i). 
Since $\om$ is compactly supported in $\Om$, using Theorem  \ref{LCE} we have 
(i)$\Longrightarrow$ (iii) $\Longrightarrow$(ii).
 As before, the testing inequality (\ref{testingI2})
is also equivalent to the Kerman--Sawyer condition
\begin{equation}\label{KSG}
\int_{\RR^n}\Big[{\rm\bf G}_{p}\om_{B}(x)\Big]^{\frac{q}{p-1}} \, dx\leq C \, \om(B),
\end{equation}
(see \cite{KS}, \cite{V2}). Note that 
\begin{equation}\label{G&I}
\int_{\RR^n}\Big[{\rm\bf G}_{p}\mu(x)\Big]^{\frac{q}{p-1}} \, dx\asymp 
\int_{\RR^n}\Big[\int_{0}^{2R}\frac{\mu(B_{t}(x))}{t^{n-p}}\Big]^
{\frac{q}{p-1}} dx,
\end{equation}
where the constants of equivalence are independent of the measure $\mu$, 
(see \cite{HW}, \cite{AH}).
From (\ref{KSG}), (\ref{G&I}), and Proposition \ref{cont} we  deduce the
implication  (ii)$\Longrightarrow$(iv). Note that Theorem \ref{sufficency2} gives
(v)$\Longrightarrow$(i).
Thus it remains to show that (iv)$\Longrightarrow$(v).
In fact, the proof of this implication is similar to the proof of (iii)$\Longrightarrow$(iv) in Theorem 
\ref{maintheorem1}. We will only sketch some crucial steps here. We define the
lower and upper parts of the truncated Wolff potential ${\rm\bf W}_{1,\,p}^{2R}$ respectively
by
\begin{eqnarray*}
{\rm\bf L}_{r}^{2R}\mu(x)=\int_{r}^{2R}\Big[\frac{\mu({B_{t}(x)})}{t^{n-p}}\Big]^{\frac
{1}{p-1}}\frac{dt}{t},~~~~ 0<r<2R,~x\in\RR^n
\end{eqnarray*}
and
\begin{eqnarray*}
{\rm\bf U}_{r}^{2R}\mu(x)=\int_{0}^{r}\Big[\frac{\mu({B_{t}(x)})}{t^{n-p}}\Big]^{\frac
{1}{p-1}}\frac{dt}{t},~~~~ 0<r<2R,~x\in\RR^n.
\end{eqnarray*}
Since $R={\rm diam}(\Om)$ and $\om\in M_{B}^{+}(\Om)$, 
to prove (\ref{pointwise2}), it is enough to verify that, for $x\in\Om$,
\begin{eqnarray}\label{mur}
\int_{0}^{2R}\Big[\frac{\mu_{r}(B_{r}(x))}{r^{n-p}}\Big]
^{\frac{1}{p-1}}\frac{dr}{r}
\leq C \, {\rm\bf W}_{1,\,p}^{2R}\om(x),
\end{eqnarray}
and
\begin{eqnarray}\label{lambdar}
\int_{0}^{2R}\Big[\frac{\lambda_{r}(B_{r}(x))}
{r^{n-p}}\Big]^{\frac{1}{p-1}}\frac{dr}{r}
\leq C \,  {\rm\bf W}_{1,\,p}^{2R}\om(x), 
\end{eqnarray}
where $d\mu_{r}=({\rm\bf U}_{r}^{2R}\om)^q dx$,  $d\lambda_{r}=({\rm\bf L}_
{r}^{2R}\om)^q dx$ and $0<r<2R$. The proof of (\ref{mur}) is the same as before. 
For the proof of (\ref{lambdar}), we need an estimate similar to (\ref{infinite}).
Namely,
\begin{equation}\label{infinite2}
\int_{r}^{4R}\Big[\frac{\om(B_{t}(x))}{t^{n-p}}\Big]^{\frac{1}{p-1}}\frac{dt}{t}
\leq C(R,\om(\Om)) \, r^\frac{-p}{q-p+1}
\end{equation} 
for all $0<r\leq 4R$ and $x\in\Om$. 
In fact, note that for $0<t<R/2$ and $y\in B_{t}(x)$, 
\begin{eqnarray*}
{\rm\bf W}_{1,\,p}^{2R}\om_{B_{t}(x)}(y)&\geq& \int_{2t}^{2R}\Big[
\frac{\om(B_{\tau}(y)\cap B_{t}(x))}{\tau^{n-p}}\Big]^{\frac{1}{p-1}}\frac
{d\tau}{\tau}\\
&\geq& C(n,p)\Big[\frac{\om(B_{t}(x))}{t^{n-p}}\Big]^{\frac{1}{p-1}}. 
\end{eqnarray*}    
As before, from this inequality and (\ref{testingW2}) one gets
\begin{equation}\label{ball}
\om(B_{t}(x))\leq C t^{n-\frac{pq}{q-p+1}},\hspace*{.3in} 0<t<R/2.
\end{equation}
To prove (\ref{infinite2}), we can assume that $0<r<R/2$ and write the left-hand
side of (\ref{infinite2}) as
\begin{equation}\label{split}
\int_{r}^{R/2}\Big[\frac{\om(B_{t}(x))}{t^{n-p}}\Big]^{\frac{1}{p-1}}\frac{dt}{t}
+\int_{R/2}^{4R}\Big[\frac{\om(B_{t}(x))}{t^{n-p}}\Big]^{\frac{1}{p-1}}\frac{dt}{t}.
\end{equation}
Applying (\ref{ball}) to the first term of (\ref{split}) and using the fact that 
$\om\in {\mathcal M}_{B}^{+}(\Om)$ in the second term of (\ref{split}), we finally obtain
(\ref{infinite2}). This completes the proof of (iv)$\Longrightarrow$(v), and so Theorem \ref{maintheorem3} is proved. 
\end{proof}
\begin{remark}\label{generalom}{\rm From the proof of  Theorem \ref{maintheorem3} we see that 
if $\om$ is not assumed to be compactly supported in 
$\Om$, then any one of the conditions (ii)--(v) is still sufficient for the solvability 
of (\ref{ep-equation2}) for some $\epsilon>0$.
}\end{remark}

\begin{theorem}\label{removeforpH} Let $E$ be a relatively closed subset of $\Om$. Suppose that 
${\rm Cap}_{{\rm\bf G}_{p},\,\frac{q}{q-p+1}}(E)=0$. Then any solution $u$ of
\begin{equation}
\label{onomminuse}
\left\{\begin{array}{c}
 u {\rm~is~} \mathcal{A}{\text-}{\rm superharmonic~in~}\Om\setminus E,\\
u\in L^{q}_{\rm loc}(\Om\setminus E), \quad u \ge 0,\\
-{\rm div}\mathcal{A}(x, \nabla u)=u^{q} \quad {\rm in}\quad 
\mathcal{D}'(\Om\setminus E),  
\end{array}
\right.
\end{equation} 
is also a solution of 
\begin{equation}
\label{onom}
\left\{\begin{array}{c}
 u {\rm~is~} \mathcal{A}{\text-}{\rm superharmonic~in~}\Om,\\
u\in L^{q}_{\rm loc}(\Om), \quad u \ge 0,\\
-{\rm div}\mathcal{A}(x, \nabla u)=u^{q} \quad {\rm in}\quad \mathcal{D}'(\Om).
\end{array}
\right.
\end{equation}
Conversely, if $E$ is a compact set in $\Om$ such that any solution of 
(\ref{onomminuse}) is also a solution of (\ref{onom})
then ${\rm Cap}_{{\rm\bf G}_{p},\,\frac{q}{q-p+1}}(E)=0$. 
\end{theorem}
\begin{proof}
Let us prove the first part of the theorem. Since 
$${\rm Cap}_{{\rm\bf G}_{p},\,
\frac{q}{q-p+1}}(E)=0,$$
 we have 
${\rm cap}_{1,\,p}(E,\Om)=0$ where the capacity ${\rm cap}_{1,\,p}(\cdot,\Om)$ is defined
by (\ref{cap1p}) (see \cite{HKM}).
Thus $u$ can be extended so that  it is a nonnegative $\mathcal{A}$-superharmonic
function in $\Om$ (see \cite{HKM}). Let $\mu[u]$ be the Radon measure on $\Om$ 
associated with $u$, and let $\varphi$ be an arbitrary nonnegative function in $C_{0}^{\infty}(\Om)$.
As in \cite[Lemme 2.2]{BP}, we can find a sequence $\{\varphi_{n}\}$ of nonnegative
functions in $C_{0}^{\infty}(\Om\setminus E)$ such that 
\begin{equation}0\leq \varphi_{n}\leq \varphi; \qquad
\varphi_{n}\rightarrow\varphi \quad {\rm Cap}_{{\rm \bf G}_{p},\,\frac{q}{q-p+1}}{\text-}{\rm quasi~everywhere}. 
\end{equation} 
By Fatou's lemma we have 
\begin{eqnarray*}
\int_{\Om} u^{q} \, \varphi \, dx &\leq& \liminf_{n\rightarrow\infty}\int_{\Om} u^{q} \, 
\varphi_{n} \, dx\\
&=& \liminf_{n\rightarrow\infty}\int_{\Om}\varphi_{n} \, d\mu[u]\\
&\leq& \int_{\Om} \, \varphi \, d\mu[u] <\infty.
\end{eqnarray*}
Therefore $u\in L^{q}_{\rm loc}(\Om)$, and $\mu[u]\geq u^{q}$ in 
$\mathcal{D}'(\Om)$. It is then easy to see that 
$$-{\rm div}\mathcal{A}(x,\nabla u)=u^{q} +\mu^{E} \quad {\rm in} \quad \mathcal{D}'(\Om)$$
for some nonnegative measure $\mu^{E}$ such that $\mu^{E}(A)=0$ for any Borel set
$A\subset \Om\setminus E$. Moreover, by Theorem \ref{maintheorem3} and Remark 
\ref{generalom} we have, for any compact set $K\subset E$,  
$$\mu^{E}(K)\leq C(K) \, {\rm Cap}_{{\rm\bf G}_{p},\,\frac{q}{q-p+1}}(K)=0.$$ 
Thus $\mu^{E}=0$ and $u$ solves
(\ref{onom}). \\
\indent The second part of the theorem is proved in the same way as in 
the linear case $(p=2)$ using the existence results in Theorem \ref{maintheorem3}. 
We refer to  \cite{AP} for details.
\end{proof}

%****************************************************************************
%****************************************************************************

\section{Hessian equations}\label{hessianequation}
In this section, we study a fully nonlinear counterpart of the theory 
presented in the previous sections. Here the notion of 
$k$-subharmonic ($k$-convex) functions associated with the fully nonlinear $k$-Hessian operator $F_{k}$, 
$k=1,..., n$, introduced by Trudinger and Wang in \cite{TW1}--\cite{TW3} 
will play a role similar to that of $\mathcal{A}$-superharmonic
functions in the quasilinear theory.\\
\indent Let $\Om$ be an open set in $\RR^n$, $n\geq 2$. For $k=1,..., n$ and 
$u\in C^{2}(\Om)$, the $k$-Hessian operator $F_{k}$ is defined by 
\begin{eqnarray*}
F_{k}[u]=S_{k}(\lambda(D^{2}u)),
\end{eqnarray*}
where $\lambda(D^{2}u)=(\lambda_{1},...,\lambda_{n})$ denotes the eigenvalues of 
the Hessian matrix of second partial derivatives $D^{2}u$, and $S_{k}$ is the $k^{th}$ symmetric
function on $\RR^n$ given by 
\begin{eqnarray*}
S_{k}(\lambda)=\sum_{1\leq i_{1}<\cdots<i_{k}\leq n}\lambda_{i_{1}}\cdots\lambda_{i_{k}}.
\end{eqnarray*}
 Thus $F_{1}[u]=\Delta u$ and $F_{n}[u]=\det D^{2}u$. Alternatively, we may also write 
\begin{eqnarray*}
F_{k}[u]=[D^{2}u]_{k},
\end{eqnarray*}
where for an $n\times n$ matrix $A$, $[A]_{k}$ is the $k$-trace of $A$, i.e., the
sum of its $k\times k$ principal 
minors.  Several equivalent definitions of $k$-subharmonicity
were given in \cite{TW2}, one of which involves the language of viscosity 
solutions: An upper-semicontinuous function $u: 
\Om\rightarrow [-\infty, \infty)$ is said to be  $k$-subharmonic in $\Om$, $1\leq k\leq n$, 
if $F_{k}[q]\geq 0$ for any quadratic polynomial $q$ such that  $u-q$ has 
a local finite maximum in $\Om$. Equivalently, an upper-semicontinuous function 
$u: \Om\rightarrow [-\infty,\infty)$ is $k$-subharmonic in $\Om$ if, for every 
open set $\Om'\Subset\Om$ and for every function $v\in C^{2}_{\rm loc}(\Om')\cap C^{0}(\overline{
\Om'})$   
satisfying $F_{k}[v]\geq 0$ in $\Om'$, the following implication holds:
\begin{eqnarray*}
u\leq v {\rm ~on~} \partial\Om' \Longrightarrow u\leq v {\rm~in~} \Om',
\end{eqnarray*}
(see \cite[Lemma 2.1]{TW2}). Note that a function $u\in C^{2}_{{\rm loc}}(\Om)$ is $k$-subharmonic 
if and only if 
\begin{eqnarray*}
F_{j}[u]\geq 0 {~\rm in~} \Om {~\rm for ~ all~} j=1,\dots, k.
\end{eqnarray*}
We denote by $\Phi^{k}(\Om)$ the class of all $k$-subharmonic 
functions in $\Om$ which are not identically equal to $-\infty$ in each component of  $\Om$.
It was proven in \cite{TW2} that $\Phi^{n}(\Om)\subset\Phi^{n-1}(\Om)\cdots\subset
\Phi^{1}(\Om)$ where $\Phi^{1}(\Om)$ coincides with the set of all proper 
classical subharmonic functions in $\Om$,  and $\Phi^{n}(\Om)$ is the set of functions 
convex on each component of $\Om$. \\
\indent The following  weak convergence result proved in  \cite{TW2} is fundamental 
to potential theory associated with $k$-Hessian operators. 
\begin{theorem}[TW2]\label{weakcontH} For each $u\in\Phi^{k}(\Om)$, there exists a nonnegative Borel measure 
$\mu_{k}[u]$ in $\Om$ such that \\
{\rm (i)} $\mu_{k}[u]=F_{k}[u]$ for $u\in C^{2}(\Om)$, and\\
{\rm (ii)} if $\{u_{m}\}$ is a sequence in $\Phi^{k}(\Om)$ converging in $L^{1}_{\rm loc}(\Om)$
to a function $u\in\Phi^{k}(\Om)$, then the sequence of the corresponding measures
$\{\mu_{k}[u_{m}]\}$ converges weakly   to $\mu_{k}[u]$.
\end{theorem}
The measure $\mu_{k}[u]$ in the theorem above is called the $k$-Hessian measure associated 
with $u$. Due to (i) in Theorem \ref{weakcontH} we sometimes write $F_{k}[u]$ in place of
$\mu_{k}[u]$ even in the case where $u\in\Phi^{k}(\Om)$ does not belong to $C^{2}(\Om)$. The $k$-Hessian measure
is an important tool in potential theory for $\Phi^{k}(\Om)$. It was used by D. A. Labutin to derive 
pointwise estimates for functions 
in $\Phi^{k}(\Om)$ in terms of the Wolff potential, which is an analogue of the Wolff potential 
estimates for $\mathcal{A}$-superharmonic functions in Theorem \ref{potential}. 
\begin{theorem} [\cite{L}]\label{localH} Let $u\geq 0$ be such that $-u\in \Phi^{k}(B(x,3r))$, 
where $1 \le k \le n$. If $\mu=\mu_{k}[-u]$ then 
\begin{eqnarray*}
C_{1}{\rm\bf W}_{\frac{2k}{k+1},\,k+1}^{r/8}\mu(x)\leq u(x)\leq  C_{2} \inf_{B(x,r)}u +
C_{3} {\rm\bf W}_{\frac{2k}{k+1},\,k+1}^{2r}\mu(x),
\end{eqnarray*}
where the constants $C_{1}$, $C_{2}$ and $C_{3}$ are independent of $x$, $u$, and $r$.
\end{theorem} 
The following global estimate is deduced from the preceding theorem as in the quasilinear case. 
\begin{corollary}\label{RRHE} Let $u\geq 0$ be such that $-u\in \Phi^{k}(\RR^n)$, 
where $1\leq k<\frac n 2$. If $\mu=\mu_{k}[-u]$
and $\inf_{\RR^n}u=0$ then for all $x\in\RR^n$,
\begin{eqnarray*}
\frac{1}{K} \, {\rm\bf W}_{\frac{2k}{k+1},\,k+1}\mu(x)\leq u(x)\leq K \, {\rm\bf W}_{
\frac{2k}{k+1},\,k+1}\mu(x),
\end{eqnarray*}
for some constant $K$ independent of $x$ and $u$. 
\end{corollary}

Let $\Om$ be a bounded  uniformly $(k-1)$-convex domain in $\RR^n$, that is, 
$\partial\Om\in C^{2}$ and $H_{j}(\partial\Om)>0$, $j=1, . . .,k-1$, where 
$H_{j}(\partial\Om)$ denotes the $j$-mean curvature of the boundary $\partial\Om$.   
We consider the following fully nonlinear 
problem:
\begin{eqnarray}\label{eqH}
\left\{\begin{array}{c}
F_{k}[-u]=u^{q}+\om \quad {\rm in}\quad \Om,\\
\hspace*{-.14in} u\geq 0 \quad {\rm in}\quad \Om,\\
 \hspace*{.02in}u=\varphi \quad {\rm on}\quad \partial \Om,
\end{array}
\right.
\end{eqnarray}
in the class of functions $u$ such that $-u$ is $k$-subharmonic in 
$\Om$. Here $\om$ is a Borel 
measure compactly supported in $\Om$, and the boundary condition 
in (\ref{eqH}) is understood in the classical sense.  Characterizations of the 
existence of $u\in\Phi^{k}(\Om)$ continuous near $\partial\Om$ which solves (\ref{eqH})
can be obtained using the iteration scheme in the proof Theorem 
\ref{sufficency2} together with the argument in the proof Theorem \ref{maintheorem3}. 
To do so we need an analogue of the global upper potential 
estimates on a bounded domain  given in Theorem 
\ref{om-estimate} for quasilinear operators. 
   
\begin{theorem}\label{globalH} Let $\mu$ be a nonnegative Borel measure compactly 
supported in a bounded domain 
$\Om\subset\RR^n$. Suppose that $u\geq 0$, $-u\in
\Phi^{k}(\Om)$ such that $u$ is continuous near $\partial\Om$ and solves
\begin{eqnarray*}
\left\{\begin{array}{c}
\mu_{k}[-u]=\mu+f \quad {\rm in}\quad \Om,\\
u=\varphi \quad {\rm on}\quad \partial \Om,
\end{array}
\right.
\end{eqnarray*}
where  $0\leq\varphi\in C^{0}(\partial\Om)$ and $0\leq f\in L^{s}(\Om)$ with $s>\frac{n}{2k}$ if $1\leq k\leq \frac n 2$,  
and $s=1$ if $\frac n 2<k\leq n$. Then 
for all $x\in\Om$,
\begin{eqnarray*} u(x)\leq K \, \Big[{\rm\bf W}_{\frac{2k}{k+1},\,k+1}^{2R}(\mu+f)(x)+
\max_{\partial\Om}\varphi \Big],
\end{eqnarray*}
where $R={\rm diam}(\Om)$ and $K$ is a constant independent of $x$, $u$, and $\Om$.
\end{theorem}
\begin{proof}
Suppose that the support of $\mu$ is contained in $\Om'$ for some open set
$\Om'\Subset\Om$.
Let $M=\sup_{\overline\Om\setminus \Om'} u$ and $u_{m}=\min\{u, m\}$ for $m>M$. Then 
$-u_{m}\in\Phi^{k}(\Om)$, continuous near $\partial\Om$, solves
\begin{eqnarray*}
\left\{\begin{array}{c}
\mu_{k}[-u_{m}]=\mu_{m} \quad {\rm in}\quad \Om,\\
u_{m}=\varphi \quad {\rm on}\quad \partial \Om,
\end{array}
\right.
\end{eqnarray*}
for some nonnegative Borel measure $\mu_{m}$ in $\Om$. Since $u_{m}\rightarrow u$
in $L^{1}_{\rm loc}(\Om)$, by Theorem \ref{weakcontH} we have 
\begin{eqnarray}\label{weak}
\mu_{m}\rightarrow \mu+f {\rm ~ weakly~ as ~ measures~in~} \Om.
\end{eqnarray}  
Note that $u_{m}=u$ in $\overline\Om\setminus\Om'$ since $m>M$. Thus $\mu_{m}=\mu_{
k}[u]=f$ in $\Om\setminus\overline{\Om'}$ for all $m>M$. Using this and (\ref{weak})
it is easy to see that 
\begin{eqnarray*}
\int_{\Om}\phi d\mu_{m}\rightarrow\int_{\Om}\phi d\mu + \int_{\Om}\phi f dx
\end{eqnarray*}
as $m\rightarrow\infty$ for all $\phi\in C^{0}(\overline\Om)$, i.e., 
\begin{eqnarray*}
\mu_{m}\rightarrow \mu+f {\rm ~ in~ the~ narrow ~topology~ of~ measures}.
\end{eqnarray*}  
We now take a ball $B\supset \Om$ with $B=B(x_{0}, 2R)$, $x_{0}\in\Om$ and 
consider the solutions $w_{m}\geq 0$, $-w_{m}\in\Phi^{k}(\Om)$, continuous near 
$\partial \Om$, of  
\begin{eqnarray*}
\left\{\begin{array}{c}
\mu_{k}[-w_{m}]=\mu_{m}  \quad{\rm in}\quad B,\\
w_{m}=\max_{\partial\Om}\varphi  \quad {\rm on}\quad \partial B,
\end{array}
\right.
\end{eqnarray*}
where $m>M$. Since $u_{m}$ is bounded in $\Om$ the measure $\mu_{m}$ is continuous 
with respect to the capacity ${\rm cap}_{k}(\cdot,\Om)$, and hence with respect to the capacity
${\rm cap}_{k}(\cdot,B)$ (see \cite{TW3}). Here ${\rm cap}_{k}(\cdot,\Om)$ is the $k$-Hessian
capacity defined by 
\begin{equation}\label{khessiancapacity}
{\rm cap}_{k}(E,\Om)=\sup \, \{\mu_{k}[u](E): u\in\Phi^{k}(\Om), -1<u<0\}
\end{equation}
for a compact set $E\subset\Om$.
By a comparison principle (see \cite[Theorem 4.1]{TW3}), we have $w_{m}\geq 
\max_{\partial\Om}\varphi$ in $B$, and hence
$w_{m}\geq u_{m}$ on $\partial\Om$. Thus, applying the comparison principle again, we have
\begin{eqnarray}\label{wmum}
w_{m}\geq u_{m} \quad {\rm in}\quad \Om.
\end{eqnarray}
Since $\mu_{m}\rightarrow\mu+f$ in the narrow topology of measures in $\Om$, we see 
that $\mu_{m}\rightarrow \mu+f $ weakly as measures in $B$. Therefore, arguing as 
in \cite[Sec. 6]{TW2} we can find a subsequence $\{w_{m_{j}}\}$ such that  
$w_{m_{j}}\rightarrow w$ a.e. for some $w\geq 0$, $-w\in\Phi^{k}(B)$ such that
$w$ is continuous near $\partial B$ and 
\begin{eqnarray*}
\left\{\begin{array}{c}
\mu_{k}[-w]=\mu + f \quad {\rm in}\quad B,\\
w=\max_{\partial\Om}\varphi \quad{\rm on}\quad \partial B.
\end{array}
\right.
\end{eqnarray*}
Note that from (\ref{wmum}), $w\geq u$ a.e. on $\Om$ and hence 
$w\geq u$ everywhere on $\Om$. Using this and Theorem \ref{localH} applied to the function
$w$ on $B(x, d(x))$, where $d(x)={\rm dist}(x,\partial B)$ we have, for $x\in \Om$ and
$d\nu=d\mu+fdx$,
\begin{eqnarray}\label{final}
u(x)&\leq& C\,{\rm\bf W}_{\frac{2k}{k+1},\,k+1}^{2R}(\nu)(x)+C \inf_{B(x,d(x)/3)}w\\
&\leq& C\,{\rm\bf W}_{\frac{2k}{k+1},\,k+1}^{2R}(\nu)(x)+C\, d(x)^{-n} \int_{B(x,d(x)/3)}w dy\nonumber\\
&\leq& C\Big({\rm\bf W}_{\frac{2k}{k+1},\,k+1}^{2R}(\nu)(x) + \max_{\partial\Om}\varphi+ R^{2-n/k}\nu(\Om)^{1/k}\Big),\nonumber
\end{eqnarray}
where the last inequality in (\ref{final}) follows from the estimate (6.3) in \cite{TW2}.
The proof of the Theorem \ref{globalH} is then complete by noting that 
$$\int_{R}^{2R}\Big[\frac{\nu(B_{t}(x))}{t^{n-2k}}\Big]^{\frac{1}{k}}\frac{dt}{t}
\geq C R^{2-n/k}\nu(\Om)^{\frac{1}{k}}$$
for all $x\in \Om$.
\end{proof}
 The next theorem 
is a criterion for the existence of global solutions to fully nonlinear equations with 
general measure data, which is an analogue of Theorem \ref{deltapu=mu}.
\begin{theorem} Suppose that $\mu$ is a measure in $\mathcal{M}^{+}(\RR^n)$ such  that
${\rm\bf W}_{\frac{2k}{k+1},\,k+1}\mu < \infty$ a.e. on $\RR^n$.  Then there exists  $u\geq 0$, 
$-u\in\Phi^{k}(\RR^n)$  such that 
\begin{eqnarray}\label{eqnRH}
F_{k}[-u]=\mu \quad {\rm in}\quad \RR^n, 
\end{eqnarray}
and
\begin{eqnarray}\label{WE}
\frac{1}{K} \, {\rm\bf W}_{\frac{k+1}{2k},\,k+1}\mu \leq u\leq K \, {\rm\bf W}_{\frac{k+1}{2k},\,k+1}\mu.
\end{eqnarray}
Conversely, if $u\geq 0$, $-u\in\Phi^{k}(\RR^n)$ solves (\ref{eqnRH}), then  ${\rm\bf W}_{\frac{2k}{k+1},\,k+1}\mu < \infty$ 
a.e. on $\RR^n$.
\end{theorem}
\begin{proof} The second part of the theorem is trivial in view of Corollary \ref{RRHE}.
To prove the first part we denote by $B_{m}$ the open ball in $\RR^n$ centered at the 
origin with radius $m$, $m=1,2,\dots$. Let  $u_{m}\geq 0$, $-u_{m}\in\Phi^{k}(B_{m+1})$,
continuous near $\partial B_{m+1}$ be a solution to the following Dirichlet problem
\begin{eqnarray*}
\left\{\begin{array}{c}
F_{k}[-u_{m}]=\mu_{B_{m}}  ~{\rm in}~ B_{m+1},\\
u_{m}=0 ~{\rm on}~ \partial B_{m+1}.
\end{array}
\right.
\end{eqnarray*}
By Theorem \ref{globalH} we have 
\begin{equation}\label{bound}
u_{m}\leq C\, {\rm\bf W}_{\frac{k+1}{2k},\,k+1}\mu<\infty \quad{\rm a.e.,}
\end{equation}
where $C$ is independent of $m$. Thus by passing to a subsequence we may assume that 
$u_{m}$ converges to $u$ a.e. for some $u\geq 0$ such that $-u\in\Phi^{k}(\RR^n)$. Since
${\rm\bf W}_{\frac{k+1}{2k},\,k+1}\mu\in L^{1}_{{\rm loc}}(\RR^n)$, 
the weak continuity
result (Theorem \ref{weakcontH}), (\ref{bound}) and Corollary \ref{RRHE} then imply that 
$u$ is a solution of (\ref{eqnRH})
which satisfies (\ref{WE}). 
\end{proof}

We are now in a position to establish the main results of this section.

\begin{theorem}\label{maintheoremH}
Let $\om$ be a measure in $\mathcal{M}^{+}(\RR^n)$, $1\leq k<n/2$, 
and  $q>k$. Then the 
following statements are equivalent.\\
{\rm(i)} There exists a nonnegative solution $u$  to the equation 
\begin{equation}
\label{ep-equationH}
\left\{\begin{array}{c}
\displaystyle{\inf_{x\in\RR^{n}}} \, u(x)=0, \\ 
F_{k}[-u]=u^{q}+\ep \, \om  \quad {\rm in}\quad \RR^n, 
\end{array}
\right.
\end{equation} 
such that $-u\in\Phi^{k}(\Om)\cap 
L^{q}_{\rm loc}(\RR^n)$, for some $\ep>0$.\\
{\rm(ii)} The testing inequality 
\begin{equation}
\label{testingI_{2k}}
\int_{B}\Big[{\rm\bf I}_{2k}\om_{B}(x)\Big]^{\frac{q}{k}}dx\leq C \, \om(B)
\end{equation}
holds for all balls $B$ in $\RR^{n}$.\\
{\rm(iii)} For all compact sets $E\subset\RR^n$, 
\begin{equation}\label{capH}
\om(E)\leq C \, {\rm Cap}_{{\rm\bf I}_{2k}, \frac{q}{q-k}}(E).
\end{equation}
\noindent {\rm(iv)} The testing inequality 
\begin{equation}
\label{testingWH}
\int_{B}\Big[{\rm\bf W}_{\frac{2k}{k+1},\,k+1}\om_{B}(x)\Big]^{q}dx\leq C \,  \om(B)
\end{equation}
holds for all balls $B$ in $\RR^{n}$.\\
{\rm(v)} There exists a constant C such that
\begin{equation}
\label{pointwiseH}
{\rm\bf W}_{\frac{2k}{k+1},\,k+1}({\rm\bf W}_{\frac{2k}{k+1},\,k+1}\om)^{q}\leq C \, {\rm\bf W}_{\frac{2k}{k+1},\,
k+1}\om < \infty \quad {\rm a.e.}
\end{equation}
Moreover, if the constant $C$ in (\ref{pointwiseH}) satisfies
$$C\leq \Big(\frac{q-k}{qK} \Big)^{q/k} \frac{k}{q-k},$$
where $K$ is the constant in Corollary \ref{RRHE}, then 
the equation (\ref{ep-equationH}) has a solution $u\geq 0$, 
$-u\in\Phi^{k}(\RR^n)$ with $\ep=1$ which obeys 
the two-sided estimate
\begin{eqnarray*}
C_{1} \, {\bf W}_{\frac{2k}{k+1},\,k+1}\om(x)\leq u(x)\leq C_{2} \, {\bf W}_{\frac{2k}{k+1},\,k+1}\om(x)
\end{eqnarray*}  
for all $x\in \RR^n$.\\
\end{theorem}
\begin{theorem} Let $\Om$ be a bounded uniformly $(k-1)$-convex domain in $\RR^n$.
Suppose that $\om\in\mathcal{M}_{B}^{+}(\Om)$ such that
 $\om =\mu +f$, where $\mu\in \mathcal{M}^{+}_{B}(\Om)$ with ${\rm supp}\mu
\Subset \Om$ and $0\leq f\in L^{s}(\Om)$ with
$s>n/2k$ if $1\leq k\leq n/2$ and $s=1$ if $n/2<k\leq n$. Let  $q>k$, $R={\rm diam}(\Om)$ and 
$0\leq\varphi\in C^{0}(\partial \Om)$. Assume that 
\begin{equation}\label{A}
{\rm\bf W}_{\frac{2k}{k+1},\,k+1}^{2R}({\rm\bf W}_{\frac{2k}{k+1},\,k+1}^{2R}\om)^{q}\leq A{\rm\bf W}
_{\frac{2k}{k+1},\,k+1}^{2R}\om
\end{equation}
and 
\begin{equation}\label{B}
{\rm\bf W}_{\frac{2k}{k+1},\,k+1}^{2R}\Big[{\rm\bf W}_{\frac{2k}{k+1},\,k+1}^{2R}
(\max_{\partial\Om}\varphi)^{q}\Big]^{q}
\leq B{\rm\bf W}_{\frac{2k}{k+1},\,k+1}^{2R}
(\max_{\partial\Om}\varphi)^{q},
\end{equation}
where 
\begin{eqnarray}\label{Acont}
A\leq \Big(\frac{q-k}{3^{\frac{q-1}{q}}qK}\Big)^{q/k}\Big(\frac{k}{q-k}\Big)
\end{eqnarray}
and
\begin{eqnarray}\label{Bcont}
B\leq \Big(\frac{q-k}{3^\frac{q-1}{q}qK^\frac{q}{k}} \Big)^{q/k}\Big(\frac{k}{q-k}\Big).
\end{eqnarray}
Here $K$ is the constant in Theorem \ref{globalH}. 
Then there exists a function $u\geq 0$, $-u\in \Phi^{k}(\Om)\cap L^{q}(\Om)$, continuous near 
$\partial \Om$ such that  
\begin{equation}
\label{HE}
\left\{\begin{array}{c}
F_{k}[-u]=u^{q}+\om \quad {\rm in}\quad \Om, \\
u=\varphi \quad{\rm on}\quad \partial \Om.
\end{array}
\right.
\end{equation} 
Moreover, there is a constant $C=C(n,k,q)$ such that 
$$u\leq C\Big\{{\rm\bf W}_{\frac{2k}{k+1},\,k+1}\om+ 
{\rm\bf W}_{\frac{2k}{k+1},\, k+1}(\max_{\partial\Om}\varphi)^{q}+\max_{
\partial\Om}\varphi\Big\}.$$ 
\end{theorem}
\begin{remark}{\rm Condition (\ref{B}) is redundant if $\varphi$ is small 
enough.}
\end{remark}
\begin{proof}
Let $\{u_{m}\}_{m\geq 0}$ be a sequence of nonnegative
functions on $\Om$ defined inductively by the following Dirichlet 
problems:
\begin{equation*}
\left\{\begin{array}{c}
F_{k}[-u_{0}]=\om \quad {\rm in}\quad \Om, \\
u_{0}=\varphi \quad {\rm on}\quad \partial \Om, 
\end{array}
\right.
\end{equation*} 
and
\begin{equation}\label{iterate}
\left\{\begin{array}{c}
F_{k}[-u_{m}]=u_{m-1}^{q}+\om \quad {\rm in}\quad \Om, \\
u_{m}=\varphi \quad {\rm on}\quad \partial \Om,
\end{array}
\right.
\end{equation} 
for $m\geq 1$. Here for each $m\geq 0$, $-u_{m}$ is $k$-subharmmonic and is continuous near 
$\partial\Om$. By Theorem \ref{globalH} we have
\begin{eqnarray*}
u_{0}&\leq& K\, {\rm\bf W}_{\frac{2k}{k+1},\,k+1}^{2R}\om +K \max_{\partial\Om}
\varphi\\
&=&a_{0}{\rm\bf W}_{\frac{2k}{k+1},\,k+1}^{2R}\om + b_{0}{\rm\bf W}_
{\frac{2k}{k+1},\,k+1}^{2R}(\max_{\partial\Om}\varphi)^{q} + K\max_{\partial\Om}
\varphi, 
\end{eqnarray*}
where $a_{0}=K$ and $b_{0}=0$. Thus 
\begin{eqnarray*}
u_{1}&\leq& K\,{\rm\bf W}_{\frac{2k}{k+1},\,k+1}^{2R}(u_{0}^{q}+\om)+K\max_{\partial\Om}
\varphi\\
&\leq& K\Big\{(3^{q-1}a_{0}^{q})^{\frac{1}{k}}{\rm\bf W}_{\frac{2k}{k+1},\,k+1}^{2R}
({\rm\bf W}_{\frac{2k}{k+1},\,k+1}^{2R}\om)^{q}+\\
&&(3^{q-1}b_{0}^{q})^{\frac{1}{k}}{\rm\bf W}_{\frac{2k}{k+1},\,k+1}^{2R}\Big[{\rm\bf W}_
{\frac{2k}{k+1},\,k+1}^{2R}(\max_{\partial\Om}\varphi)^{q}\Big]^{q}+\\
&& K^{\frac{q}{k}}{\rm\bf W}_{\frac{2k}{k+1},\,k+1}^{2R}(\max_{\partial\Om}\varphi)^{q}
+{\rm\bf W}_{\frac{2k}{k+1},\,k+1}^{2R}\om  \Big\} +K\max_{\partial\Om}\varphi.
\end{eqnarray*}
 Then by (\ref{A}) and (\ref{B}),
\begin{eqnarray*}
u_{1}&\leq& K[(3^{q-1}a_{0}^{q})^{\frac{1}{k}}A+1]{\rm\bf W}_{\frac{2k}{k+1},\,k+1}^{2R}\om+\\
&&K[(3^{q-1}b_{0}^{q})^{\frac{1}{k}}B+K^{\frac{q}{k}}]{\rm\bf W}_{\frac{2k}{k+1},\,k+1}^{2R}
(\max_{\partial\Om}\varphi)^{q}+  K\max_{\partial\Om}\varphi\\
&=& a_{1}{\rm\bf W}_{\frac{2k}{k+1},\,k+1}^{2R}\om +
b_{1}{\rm\bf W}_{\frac{2k}{k+1},\,k+1}^{2R}(\max_{\partial\Om}\varphi)^{q}+K\max_{\partial\Om}\varphi,
\end{eqnarray*}
where $$a_{1}=K[(3^{q-1}a_{0}^{q})^{\frac{1}{k}}A+1], \qquad b_{1}=K[(3^{q-1}b_{0}^{q})^{\frac{1}{k}}B+
K^{\frac{q}{k}}].$$
By induction we have 
$$u_{m}\leq a_{m}{\rm\bf W}_{\frac{2k}{k+1},\,k+1}^{2R}\om +
b_{m}{\rm\bf W}_{\frac{2k}{k+1},\,k+1}^{2R}(\max_{\partial\Om}\varphi)^{q}+K
\max_{\partial\Om}\varphi,$$
where 
$$a_{m+1}=K[(3^{q-1}a_{m}^{q})^{\frac{1}{k}}A+1], \qquad 
b_{m+1}=K[(3^{q-1}b_{m}^{q})^{\frac{1}{k}}B+K^{\frac{q}{k}}], $$
for all $m\geq 0$. It is then easy to see that 
$$a_{m}\leq \frac{Kq}{q-k}, \qquad b_{m}\leq \frac{K^{\frac{q}{k}+1}q}{q-k},$$
provided (\ref{Acont}) and (\ref{Bcont}) are satisfied. Thus 
\begin{eqnarray}\label{boundforun}
u_{m}&\leq& \frac{Kq}{q-k}{\rm\bf W}_{\frac{2k}{k+1},\,k+1}^{2R}\om +\\
&&+\frac{K^{\frac{q}{k}+1}q}{q-k}{\rm\bf W}_{\frac{2k}{k+1},\,k+1}^{2R}(\max_{\partial\Om}
\varphi)^{q}+K\max_{\partial\Om}\varphi.\nonumber
\end{eqnarray}
Using (\ref{A}), (\ref{B}), (\ref{boundforun}), and passing to a subsequence, we can find a 
function $u\geq 0$ such that $-u$ is $k$-subharmonic and $u_{m}^q\rightarrow u^q$ in 
$L^{1}(\Om)$. Thus in view of (\ref{iterate}) and Theorem \ref{weakcontH} we see that $u$ is a desired 
solution of (\ref{HE}). 
\end{proof}

\begin{theorem} Let $\om$ be a locally finite nonnegative measure on an open
(not necessarily bounded) set $\Om$. 
Let $q>k$, where $1\le k \le n$. Suppose that 
$u\geq 0$,  $-u\in \Ph^{k}(\Om)$ such that $u$ is a solution of
$$F_{k}[-u]=u^q +\om \quad\quad {\rm in} \quad\Om.$$ 
Then for each cube $P\in\mathcal{Q}$, where $\mathcal{Q}=\{Q\}$ is a Whitney 
decomposition of $\Om$ as before (see Sec. \ref{Om}), we have  
\begin{equation}\label{LE1}
\mu_{P}(E)\leq C \, {\rm Cap}_{{\rm\bf I}_{2k},\,\frac{q}{q-k}}(E), 
\end{equation}
if $\, \frac{2kq}{q-k}<n$, and 
\begin{equation}\label{LE2}
\mu_{P}(E)\leq C(P) \, {\rm Cap}_{{\rm\bf G}_{2k},\,\frac{q}{q-k}}(E), 
\end{equation}
if $\, \frac{2kq}{q-k}\geq n$, for all compact sets $E\subset \Om$.
Here $d\mu=u^{q}dx + d\om$, and the constant $C$ in (\ref{LE1}) does not depend on $P$ and 
$ E\subset\Om$; however, the constant $C(P)$ in (\ref{LE2}) may depend
on the side length of $P$.

Moreover, if $\frac{2kq}{q-k}<n$, and $\Om$ is a bounded $C^{\infty}$-domain then 
$$\mu(E)\leq C \, {\rm cap}_{2k,\,\frac{q}{q-k}}(E,\Om)$$
for all compact sets $E\subset\Om$, where 
${\rm cap}_{2k,\,\frac{q}{q-k}}(E,\Om)$ is defined by (\ref{capomega}). 
\end{theorem}
\begin{remark}{\rm Let $B_{R}$ be a ball such that  $B_{2R}\subset\Om$. If $\frac{2kq}{q-k}=n$ then as in 
Remark \ref{loglocal} we have
$$\mu(B_{r})\leq C (\log\tfrac{2R}{r}) ^{\frac{-k}{q-k}}$$
for all balls $B_{r}\subset B_{R}$.
}\end{remark}

\begin{theorem}\label{maintheorem4} 
Let $\om$ be a compactly supported measure in ${\mathcal M}^{+}_{B}(\Om)$,
where $\Om$ is a bounded uniformly $(k-1)$-convex domain in $\RR^n$ $(1 \le k \le n)$. Let  
$q>k$, $R={\rm diam}(\Om)$, and $\varphi\in C^{0}(\partial \Om)$, $\varphi \ge 0$. Then the following statements
are equivalent.\\
{\rm (i)} There exists a solution $u\geq 0$, $-u\in\Phi^{k}(\Om)\cap L^{q}(\Om)$,
continuous near $\partial\Om$, to the equation
\begin{equation}\label{khessianepH}
\left\{\begin{array}{c}
F_{k}[-u]=u^{q}+\epsilon\om \quad {\rm in}\quad \Om, \\
u=\epsilon\varphi \quad {\rm on}\quad \partial \Om,
\end{array}
\right.
\end{equation} 
for some $\epsilon>0$.\\
{\rm (ii)} The testing inequality 
\begin{equation}
\label{testingI2H}
\int_{B}({\rm\bf G}_{2k}\om_{B})^{\frac{q}{k}} \, dx\leq C \, \om(B)
\end{equation}
holds for all balls $B$ such that $B\cap {\rm supp}\om\not=\emptyset$.\\
{\rm (iii)} For all compact sets $E\subset{\rm supp}\om$,
$$\om(E)\leq C \, {\rm Cap}_{{\rm\bf G}_{2k},\,\frac{q}{q-k}}(E).$$
{\rm (iv)} The testing inequality 
\begin{equation*}
\int_{B}\Big[{\rm\bf W}_{\frac{2k}{k+1},\,k+1}^{2R}\om_{B}(x)\Big]^{q} \, dx\leq C \, \om(B)
\end{equation*}
holds for all balls $B$ such that $B\cap {\rm supp}\om\not=\emptyset$.\\
{\rm (v)} There exists a constant C such that
\begin{equation*}
\hspace*{.3in}{\rm\bf W}_{\frac{2k}{k+1},\,k+1}^{2R}({\rm\bf W}_{\frac{2k}{k+1},\,k+1}^{2R}\om)^{q}\leq 
C \, {\rm\bf W}_{\frac{2k}{k+1},\,k+1}^{2R}
\om < \infty \quad {\rm a.e. ~on~} \Om.
\end{equation*}
\end{theorem}
\begin{remark}{\rm As in Remark \ref{generalom}, if $\om=\mu+f$, where ${\rm supp}\mu
\Subset\Om$, and $0\leq f\in L^{s}(\Om)$ with $s>\frac{n}{2k}$ if $k\leq n/2$, 
and $s=1$ if $k>n/2$, then any one of the 
conditions (ii)--(v) is still sufficient for the solvability of the equation 
(\ref{khessianepH}) for some $\epsilon>0$.
}\end{remark}
\begin{theorem} Let $E$ be a relatively closed subset of $\Om$. Suppose that 
${\rm Cap}_{{\rm\bf G}_{2k},\,\frac{q}{q-k}}(E)=0$. Then any solution $u$ of
\begin{equation}
\label{onomminuseH}
\left\{\begin{array}{c}
  -u \in \Phi^{k}(\Om\setminus E)\cap L^{q}_{\rm loc}(\Om\setminus E), \quad u\geq 0,\\
F_{k}[-u]=u^{q}\quad {\rm in}\quad \mathcal{D}'(\Om\setminus E),
\end{array}
\right.
\end{equation} 
is also a solution of 
\begin{equation}
\label{onomH}
\left\{\begin{array}{c}
 -u \in \Phi^{k}(\Om)\cap L^{q}_{\rm loc}(\Om), \quad u\geq 0,\\
F_{k}[-u]=u^{q} \quad {\rm in}\quad \mathcal{D}'(\Om).
\end{array}
\right.
\end{equation}
Conversely, if $E$ is a compact set in $\Om$ such that any solution of 
(\ref{onomminuseH}) is also a solution of (\ref{onomH}), 
then ${\rm Cap}_{{\rm\bf G}_{2k},\,\frac{q}{q-k}}(E)=0$. 
\end{theorem}
\begin{proof} To prove this theorem, we proceed as in the proof of Theorem \ref{removeforp}.
For the first statement, note that if ${\rm Cap}_{{\rm\bf G}_{2k},\,\frac{q}{q-k}}(E)=0$ then 
${\rm Cap}_{{\rm\bf G}_{\frac{2k}{k+1}},\,k+1}(E)=0$ and $k<\frac n 2$ (see \cite{AH}),  
which implies that 
$${\rm cap}_{k}(E,B)=0$$ 
for a ball $B\Supset\Om\supset E$  due to Theorem \ref{handb} below. 
Here ${\rm cap}_{k}(\cdot,\Om)$ is the (relative) $k$-Hessian capacity associated 
with the domain $\Om$  (see (\ref{khessiancapacity})).
Thus by \cite[Theorem 4.2]{L}, $E$ is a $k$-polar set, i.e., $(-\infty)$-set of a 
$k$-subharmonic function in $\RR^n$. It is then easy to see that the 
function $\tilde{u}$ defined by 
\begin{equation}
\tilde{u}(x)=\left\{\begin{array}{c}
u(x),\quad x\in\Om\setminus E,\\
\displaystyle{\limsup_{\substack{y\rightarrow x,\,y\not\in E}}} \, u(y),\quad x\in E,
\end{array}
\right.
\end{equation}
belongs to $\Phi^{k}(\Om)$, and $-\tilde{u}$ is an extension of $u$. The rest of the 
proof is then the same as before.
\end{proof}
\begin{theorem}\label{handb}
Let $1\leq k<\frac n 2$ be an integer. Then 
\begin{equation}\label{comparecapacity}
M_{1} \, {\rm Cap}_{{\rm\bf G}_{\frac{2k}{k+1}}, \, k+1}(E)\leq {\rm cap}_{k}(E,\Om) \leq M_{2} \, 
{\rm Cap}_{{\rm\bf G}_{\frac{2k}{k+1}}, \, k+1}(E)
\end{equation}
for any compact set $E\subset\overline{Q}$ with $Q\in\mathcal{Q}$, where the constants
$M_{1}$, $M_{2}$ are independent of $E$ and $Q$. 
\end{theorem}
\begin{proof} Let $R$ be the diameter of $\Om$. From Wolff's inequality it follows that 
${\rm Cap}_{{\rm\bf G}_{\frac{2k}{k+1}}, \, k+1}(E)$ is equivalent to
$$\sup \, \{\mu(E): \, \, \mu\in M^{+}(E), \quad 
{\rm {\bf W}}^{4R}_{\frac{2k}{k+1}, \,  k+1}\mu\leq 1 {\rm ~on~} {\rm supp}\mu\},$$
for any compact set $E\subset\Om$ (see \cite[Proposition 5]{HW}). To prove the 
left-hand inequality
in (\ref{comparecapacity}), let
$\mu\in M^{+}(E)$ such that  
${\rm\bf W}^{4R}_{\frac{2k}{k+1}, \, k+1}\mu\leq 1$ on ${\rm supp}\mu$, and let 
$u\in\Phi^{k}(B)$ 
be a nonpositive solution of  
\begin{equation*}
\left\{\begin{array}{c}
F_{k}[u]=\mu \quad {\rm in}\quad B\\
u=0 \quad {\rm on}\quad \partial B, 
\end{array}
\right.
\end{equation*}
where $B$ is a ball of radius $R$ containing $\Om$. By Theorem \ref{globalH} and the 
boundedness principle for nonlinear potentials (see \cite{AH}), we have
$$\m{u}\leq C \, {\rm\bf W}^{4R}_{\frac{2k}{k+1},\,k+1}\mu\leq C .$$
Thus 
$$\mu(E)=\mu_{k}[u](E)\leq C \, {\rm cap}_{k}(E,\Om),$$
which shows that  
$${\rm Cap}_{{\rm\bf G}_{\frac{2k}{k+1}},\,k+1}(E)\leq C \, {\rm cap}_{k}(E,\Om).$$
To prove the upper estimate in (\ref{comparecapacity}), we let $Q\in\mathcal{Q}$, and fix a compact 
set $E\subset\overline{Q}$. 
Note that for $\mu\in M^{+}(E)$ and $x\in E$ we have
$${\rm\bf W}^{4R}_{\frac{2k}{k+1},\,k+1}\mu(x)={\rm\bf W}^{2{\rm diam}(Q)}_{\frac{2k}{k+1},\,k+1}\mu(x)
+\int_{2{\rm diam}(Q)}^{4R}\Big[\frac{\mu(E)}{t^{n-2k}}\Big]^{\frac{1}{k}}\frac{dt}{t}.$$
Thus, for $k<\frac n 2$,
\begin{equation}\label{4R}
{\rm\bf W}^{4R}_{\frac{2k}{k+1},\,k+1}\mu(x)\leq C\,{\rm\bf W}^{2{\rm diam}(Q)}_{
\frac{2k}{k+1},\,k+1}\mu(x),\quad\quad \forall x\in E.
\end{equation}
Now for $u\in\Phi^{k}(\Om)$ such that $-1<u<0$ by Theorem \ref{localH} we 
obtain
$${\rm\bf W}^{2{\rm diam}(Q)}_{\frac{2k}{k+1},\,k+1}\mu_{E}(x)\leq 
{\rm\bf W}^{2{\rm diam}(Q)}_{\frac{2k}{k+1},\,k+1}\mu(x)\leq C \, \m{u(x)}\leq C,$$
for all $x\in E$, where $\mu=\mu_{k}[u]$. 
Thus, we deduce from (\ref{4R}) that
$${\rm\bf W}^{4R}_{\frac{2k}{k+1},\,k+1}\mu_{E}(x)\leq C,\quad \quad\forall x\in E,$$
 which implies 
\begin{equation}\label{finale}
\mu(E)\leq C \, {\rm Cap}_{{\rm\bf G}_{\frac{2k}{k+1}},\,k+1}(E).
\end{equation}
 Finally, the definition of  ${\rm cap}_{k}(\cdot,\Om)$ and (\ref{finale}) then give 
$${\rm cap}_{k}(E,\Om)\leq C{\rm Cap}_{{\rm\bf G}_{\frac{2k}{k+1}},\,k+1}(E),$$
which completes the proof of the theorem.
\end{proof}
\begin{remark}{\rm If $\Om$ is a $C^{\infty}$-domain in $\RR^n$, and $1\leq k<\frac n 2$,
 then by the quasiadditivity of the capacity 
${\rm cap}_{\frac{2k}{k+1},
\,k+1}(\cdot,\Om)$ (see Theorem \ref{quasiadditivity}) we have the following upper estimate for the 
$k$-Hessian capacity ${\rm cap}_{k}
(\cdot,\Om)$: There exists a constant $C>0$ such that for any compact set $E\subset\Om$,  
$${\rm cap}_{k}(E,\Om)\leq C \, {\rm cap}_{\frac{2k}{k+1},\,
k+1}(E,\Om).$$
}\end{remark}

\end{document}